\documentstyle[leqno,amsfonts,12pt]{article}
\newsymbol \restriction 1316
\begin{document}

\begin{center}
{\Large {\bf{ National Academy of Sciences of Ukraine } }}
\end{center}
 
\begin{center}
{\Large {\bf{ Institute of Mathematics } }}
\end{center}

\bigskip
\bigskip
\bigskip
\bigskip
\bigskip
\begin{flushright}
{ \textsc{Preprint 98.7}	}
\end{flushright}

\bigskip
\bigskip
\bigskip
\bigskip
\bigskip
\bigskip
\bigskip
\begin{tabular}{|lcl|}  \hline

& & \\   & &       \\ 
 &
{\Large {\bf{ Polynomial Approximation in 
$\mathbf{L_p ( R , d  \mu  )} $.  }} {\bf{I}}  }
& \\  & & 
\\ 		& & 
\\ 
  & {\textsc{Andrew G. Bakan}} & 
 \\ & &  
 \\  & & 
 \\  & & \\
 & & \\ 
 &October \ 1998 \ \hfill 98/10 &  	\\
 & & \\
  \hline
 \end{tabular}

\vfill 
\begin{center}
{\textsc{Kyiv \ \ \  \ 1998}}
\end{center}

\newpage
 \thispagestyle{empty}
\bigskip
\bigskip
\bigskip
\bigskip
\noindent
Polynomial Approximation in 
$ L_p ( R , d  \mu  ) $. \ I \  \ 
/ \  Andrew G. Bakan . -  Kyiv, \\ 1998. -  45 p. 
- (Prepr. / Nat.Acad.Sci.Ukr. Inst.Mathematics; 98:7)

\bigskip
\bigskip
\bigskip
\bigskip

\begin{flushleft}
For arbitrary $ w : {{R}} \to [0,1] \ $ the general form of 
the continuous linear functionals \\ on the space $C_w^0 \ $ 
of all functions $f \ $ continuous on the real line,   
$\lim\limits_{| x | \to + \infty } w(x)f(x) = 0 \ $, 
equipped  with seminorm  \ \ $ ||f||_w := \sup\limits_{x \in 
{{R}} } w(x) |f(x)| \ $, \ is found. The weighted analog of 
the Weierstrass  polynomial  approximation theorem and 
a new version of M.G. Krein's \\ theorem  about partial fraction
decomposition of the reciprocal of an entire function \\
are established. New descriptions of the Hamburger   
and Krein classes of entire  functions are obtained. 
Preprint includes the final representation 
of all those measures $\mu \ $ for \\ which algebraic 
polynomials are dense in $ L_p ( R , d  \mu  ) $. 
\end{flushleft}
                                    
\bigskip
\bigskip
\bigskip
\bigskip

\vfill
									
\noindent
Department of Functions Theory \\
Kiev Inst. of Mathematics \\
National Acad. Sci. Ukraine \\
Tereschenkivska str.3, Kyiv - 4 \\
252601, Ukraine \\
e-mail: mathkiev@imath.kiev.ua

\bigskip
\bigskip
\bigskip
\bigskip

\noindent
The print is approved by the Academic Council of Institute 
of Mathematics, \\
National Academy of Sciences of Ukraine 

\bigskip
\bigskip
\bigskip
\bigskip

\noindent
@ Andrew G. Bakan, 1998 {\hfill e-mail: andrew@bakan.kiev.ua}\\

\newpage
 {\LARGE{ {\bf{Polynomial Approximation in 
$\mathbf{L_p ( R , d  \mu  )} $. }}}} 

\medskip
\noindent
{\Large {\bf{Part I. General Results and Representation 
Theorem. } }}

\bigskip
\begin{center}
{\normalsize{ \textsc{Andrew G. Bakan}}}
\end{center}

\begin{center}
{ \small{ \textrm{October 1998}} }
\end{center}

\bigskip
\bigskip
\begin{center}
{\bf {\large{CONTENTS}} }
\end{center}

 Introduction \dotfill{2} 

\medskip 
 \begin{center}
 PART I. GENERAL RESULTS AND REPRESENTATION THEOREM.
 \end{center}

Chapter I. \  Banach space associated with seminormed space 
$ C^0_w  \ $ \dotfill{4}

\medskip
1.1. \ Notations \dotfill{4}

1.2. \ Background \dotfill{5}

1.3. \ Associated normed space $ N^0_w  \ $ \dotfill{5}

1.4. \ Banach space $ B^0_w  \ $ \dotfill{6}

1.5. \ Supplement to S.N. Mergelyan's Theorem  \dotfill{8}

1.6. \ General form of the functionals in 
       $ \left( C^0_w \right)^* \ $ \dotfill{9}

1.7. Proofs	\dotfill{9}

\medskip
Chapter II. \  Hamburger criterion of the polynomial density
in $C^0_w \ $ and $ L_p (\mu )  \ $  \dotfill{14}

\medskip
2.1. \ Notations  \dotfill{14}

2.2. \ Background \dotfill{14}

2.3. \ Main Theorem \dotfill{15}

2.4. \ Proof of Theorem 2.1 \dotfill{16}

 \medskip
Chapter III. \  Hamburger and Krein classes of entire 
functions  \dotfill{18}

\medskip
3.1. \ Notations and Definitions \dotfill{18}

3.2. \ Normal polynomial families in $ L_p (\mu )  \ $
\dotfill{20}

3.3. \ New version of M.G. Krein's Theorem and its 
corollaries\dotfill{23}

3.4. \ Strictly normal polynomial families
\dotfill{29}

3.5. \ Proof of Theorem 3.2. \dotfill{34}

\medskip
Chapter IV. \  Criterion of the polynomial density
in $ L_p (\mu )  \ $  \dotfill{36}

\medskip
4.1. \ Representation Theorem \dotfill{36}

4.2. \ Preliminary Lemmas \dotfill{38}

4.3. \ Proof of Theorem 4.1. \dotfill{40}

 \medskip
References  \dotfill{44}

\newpage
 
 \begin{center}
   {\large{  INTRODUCTION }}
  \end{center}
 
This paper is devoted to the weighted polynomial approximation
problem on the real line.

Let $w(x)\ $ be a nonnegative function of real values $x \ $, 
such that
for each $n = 0, 1, 2, \ldots $, \ $x^n w(x) \ $ is bounded.
In 1924 S.Bernstein [10] asked for conditions on $w \ $ such 
that
the algebraic polynomials ${\mathcal{P}} \ $ are dense in the 
space $C^0_w \ $ of all functions $f \ $ continuous on 
${\mathbb{R}} \ $, 
satisfying $w(x) f(x) \to 0 \ $ as $|x| \to + \infty \ $,
where  $C^0_w \ $ is equipped with the seminorm 
$||f||_w := \sup\limits_{x \in {\mathbb{R}}} w(x)|f(x)| \ $ 
(for a more explicit survey see [1, 30, 32, 40, 41]).

In 1937 S. Isumi and T. Kawata [20] showed that if 
functions $w(x) \ $ 
and $- \log w(e^x ) \ $ are even 
and convex on the real line, respectively, then algebraic 
polynomials 
${\mathcal{P}} \ $ are dense in the space $C^0_w \ $ if and 
only if 
$$ \int\limits_{{\mathbb{R} } }^{} \frac{\log w(x)}{1 + x^2 }
\ d  x \ \ = \ \ - \infty \ . \eqno{(1)} $$

In 1947 N. Akhiezer and S. Bernstein (see  [32, 1]) proved 
that 
a necessary and sufficient condition for the density of
${\mathcal{P}} \ $  in  $C^0_w \ $ is that
$$ \sup_{P \in {\mathfrak{M}}_w } \ \ 
\int\limits_{{\mathbb{R} } }^{} \frac{\log | P(x) | }{1 + 
x^2 } \ d  x \ \ = \ \ + \infty \  , \eqno{(2)} $$

\noindent
where ${\mathfrak{M}}_w := \left\{ P \in {\mathcal{P}} \ | \ 
w(x)|P(x)| \leq  1 + |x| \ \ \forall \ x \in {\mathbb{R} } \ 
\right\} \ $. It was shown in 1956 by \\ S. Mergelyan [32] that 
condition (2) is equivalent to
$$  \ \ \int\limits_{{\mathbb{R} } }^{} \frac{ \log \left[ 
\sup\limits_{P \in {\mathfrak{M}}_w } | P(x) | \right] }{1 + 
x^2 } \ d  x \ \ = \ \ + \infty \  . \eqno{(3)} $$

In 1959 L. de Branges [12] obtained a remarkable theorem for 
functions $w \ $
which are positive and continuous on the real line. He 
proved that ${\mathcal{P}} \ $ is dense in  $C^0_w \ $ if and 
only if for any real entire function $F \ $ of 
exponential type all whose 
zeros $\Lambda_F \  $ are real and simple and
which satisfies:
$$  \ \ \int\limits_{{\mathbb{R} } }^{} \frac{ \log^+   
| F(x) |  }{1 + x^2 } \ d  x \ \ < \ \ + \infty \  , 
\eqno{(4)} $$

\noindent
where $\log^+ x := \max \{ 0, \log x \} \ $, $x \geq 0 \ $, the 
following relation holds:
$$ \sum_{\lambda \in \Lambda_F  }^{} \frac{1}{ w(\lambda )
|F^{\prime } (\lambda ) | } \ = \ \infty  . \eqno{(5)} $$

In 1989 B.Ja.Levin [30] extended conditions (2) and (3) to 
all
spaces $L_p ({\mathbb{R}} ,  d \mu ) \ $, $1 \leq p < \infty 
\ $, where $\mu \ $ is a positive Borel measure on the real 
line with finite moments of any order:
$$ \int\limits_{{\mathbb{R}} }^{} |x|^n \ d  \mu (x) \ < \ 
\infty \ \ \forall \ n = 0, 1, 2, \ldots \ , $$

\noindent
and unbounded support. He proved that each of the conditions 
(2) and (3) represents a necessary and sufficient condition 
for polynomials to be dense in $L_p ({\mathbb{R}} ,  
d \mu ) \ $ where  ${\mathfrak{M}}_w  \ $ is replaced
by
$$ {\mathfrak{M}}_p  \ := \ \left\{ P \in {\mathcal{P}} \ | \ 
\ \ \int\limits_{{\mathbb{R} } }^{} \frac{  | P(x) |^p  }{( 1 
+ |x| )^p } \ d  \mu (x) \ \leq \ 1  \ \right\} \ . $$

\noindent
It should be noted here that the condition (3) for $p = 2 \ $
coincides with M. Riesz's theorem (1922) in classical
moment theory [35; 36; 2, Th. 2.4.1].

In 1996 M. Sodin and P. Yuditskii [41] found a simpler 
proof of de Branges theorem and proved its validity assuming 
only  the upper semicontinuity $w \ $ on ${\mathbb{R}} \ $. 
Moreover, in de Branges condition (5), they have 
replaced the function $F \ $ by an arbitrary real entire function 
$B \ $ of minimal exponential type with only simple real zeros 
$ \Lambda_B \ \subseteq  \left\{ x \in
{\mathbb{R}} \ | \ w(x) > 0 \ \right\} \ $. In 1998 M. Sodin 
and A. Borichev [11] established a criterion similar to (5) for 
polynomial density in all spaces $L_p ({\mathbb{R}} ,  
d \mu ) \ $, $1 \leq p < \infty \ $, under the condition that 
measure $\mu \ $ is discrete and for some positive number 
$a \ $:
$$ \sum\limits_{\lambda \in \ {\mathrm{supp}} \, \mu }^{} 
\frac{1}{(1 + |\lambda |)^a } \ < \ \infty \ . $$

In the first part of that paper we will extend de Branges 
condition (5) to all spaces $L_p ({\mathbb{R}} ,  
d \mu ) \ $, $1 \leq p < \infty \ $, without any additional 
assumption about measure $\mu \ $, and in the second part,  
obtain a new analytical proof of these conditions, showing 
their real nature from the point of view of extremal problems 
theory.

In the first Chapter for an arbitrary function 
$w :{ \mathbb{R}}  \to [0,1] \ $,  we give a complete 
description of the Banach space $B^0_w \ $ associated with 
the seminormed space $C^0_w \ $ (Theorem 1.1). This 
description, under the condition that ${\mathcal{P}} \ $ is dense 
in  $C^0_w \ $,  makes it possible  in Theorem 1.2 to 
characterize all functions $f : \left\{ \ x \in 
{\mathbb{R}} \ | \ w(x) > 0 \  \right\} \to   
{\mathbb{R}} \ $ which can be approximated in the seminorm 
$|| \cdot ||_w \ $ by polynomials. That is why Theorem 1.2  
represents a supplement to S. Mergelyan's theorem 
[32, Th.7 ] in those cases when polynomials are dense in 
the space $C^0_w \ $. Besides that, the weighted analog of the 
Weierstrass polynomial approximation theorem is derived from 
Theorem 1.2 when the set $\left\{ \ x \in 
{\mathbb{R}} \ | \ w(x) > 0 \ \right\} \ $ is bounded.

In the Chapter II, Hamburger criterion of polynomial density,
known in the classical theory of moments, has 
been extended  to all spaces $L_p ({\mathbb{R}} ,  
d \mu ) \ $, $1 \leq p < \infty \ $, and $C^0_w \ $.

Chapter III contains a new version of M. Krein's 
theorem about the partial fraction decomposition of the 
reciprocal of an entire function (Theorem 3.1). Its 
Corollaries 3.1 and 3.2 give a new characterization 
of the Hamburger and Krein classes of entire functions. 
Strictly normal polynomial families are introduced in 
section 3.4 and sufficient conditions to have such 
property are found in the Theorem 3.3.

Chapter IV includes the main result of this paper 
(Theorem 4.1), which allows us to formulate conditions similar to 
(5) in all spaces $L_p ({\mathbb{R}} ,  
d \mu ) \ $, $1 \leq p < \infty \ $.

\bigskip
{\textbf{Acknowledgements.}} \ The author thanks 
Professors  Christian Berg, Matts Esse'n, 
Iossif V. Ostrovskii and Mikhail Sodin for the discussions 
and informational help which initiated that investigation.

This work was done in the framework of the INTAS 
research network 96-0858 "New methods on theory of entire 
and subharmonic functions and their applications to 
probability theory."

\newpage
\begin{center}
{\bf{ \large{CHAPTER I. }}}  {{ \large{Banach space 
associated with seminormed space} $ C^0_w  \ $ }}
\end{center}

\bigskip
\bigskip
{\bf 1.1. Notations. \ } 
{\small  \ Everywhere below in this Chapter, only real linear 
spaces and spaces of real-valued functions are considered. 
It is worth to remind [13, 1.10.2] that the pair 
$X = \left({\mathcal{L}} (X), \left\| \cdot \right\|_X \right)
\ $ is called a seminormed space if 
${\mathcal{L}} (X) \ $ is a linear space and 
$\left\| \cdot \right\|_X \ $ is defined on 
${\mathcal{L}} (X) \ $ seminorm. We will write $X \ $ instead
of $ {\mathcal{L}} (X) \ $, i.e. $ X = \left( X , \, 
\left\|\cdot \right\|_X \right) \ $. 
Denote by $X^* \ $ the Banach space [13, 1.10.6] of all linear 
continuous functionals  $L \ $ on the seminormed space 
$X \ $, equipped with norm  $ \left\| L \right\| \ := 
\sup \{ \left| L(x) \right|  \  \left| \ x \in X , \  
\left\| x \right\|_X \leq 1  \right. \} \ $. For two 
seminormed spaces $X \ $ and $Y \ $ notation $X \equiv Y \ $
indicates that $X \ $ and $Y \ $ coincide
identically, i.e. $X = Y \ $ and $\left\| x \right\|_X =
\left\| x \right\|_Y \ $ $ \forall x \in X \ $. For two 
normed spaces $X \ $ and $Y \ $ notation $X \cong Y \ $ 
means that $X \ $ and $Y \ $ are {\it{ isometric}}, i.e. 
there exists such linear transformation $U : X \to Y \ $, 
that: a) $U(X) = Y \ $; b) $\left\| U(x) \right\|_Y = 
\left\| x \right\|_X \ $ $ \forall x \in X \ $. If $ X = 
\left( X , \, \left\|\cdot \right\|_X \right) \ $ is a 
seminormed space then the normed factor space 
$ X\setminus N_X = \left( X\setminus N_X , \, \left\|\cdot 
\right\|_{X \setminus N_X } \right) \ $ whose elements are 
classes $\pi (x) := x + N_X \ $, $\left\| \pi (x)  
\right\|_{X \setminus N_X } :=  \left\| x \right\|_X \ $ 
$ \forall x \in X \ $ and $N_X := \left\{ \ x \in X \ | \ 
\left\| x \right\|_X = 0 \right\} \ $ is said to be  
[13, 1.10.2] a normed space {\it associated } with seminormed 
space $X  \ $.

Let $A \subseteq {\mathbb{R}} \ $. The closure of $A \ $ is 
denoted by $\overline{A} \ $, and $\chi_{A} (x) := 
\left\{ \begin{array}{ll}
1 , & x \in A ;  \\
0 , & x \notin A .
\end{array}\right. \ $ Let $C (A) \ $ denote the
linear space of all continuous on $A \ $ functions 
$f : A \to {\mathbb{R}} \ $;  $C^0 (A) \ := \left(  C^0 (A) , 
\, \left\| \cdot \right\|_{ C (A)  } \right) \ $ - Banach space 
of such  bounded on $A \ $, i.e. $\left\| f \right\|_{ 
C (A)  } := \sup\limits_{x \in A } |f(x) | < \infty \ $, 
functions $f \in C (A) \ $ that  $\lim\limits_{x \in A, \ |x| 
\to \infty } f(x) = 0 \ $, if $A \ $ is unbounded;   
${\mathbb{Z}}_0 \ $ - the 
set of all nonnegative integers. Function $f \in C 
({\mathbb{R}}) \ $ is called {\it compactly supported } if 
it's equal to zero outside of some compact subset of the real 
line.

For $ A\subseteq  B \subseteq {\mathbb{R}} \ $ and $h : 
B \to {\mathbb{R}} \ $ symbol $h\!\restriction_A    \ $  
denotes the function $h\!\restriction_A : A \to {\mathbb{R}} 
\ $, $h\!\restriction_A (x) = h(x) \ $ $\forall \ x \in A 
\  $.  For every  $n \in {\mathbb{Z}}_0 \ $, let 
${\mathcal{P}}_n := {\mathcal{P}}_n [{\mathbb{R }}] \ $ and 
$ {\mathcal{P}}_n [{\mathbb{C }}] \ $ denote the sets of all 
algebraic polynomials of degree at most $n \ $ with real and 
complex coefficients, respectively, and let also
${\mathcal{P}} :=  \bigcup_{n \in 
{\mathbb{Z}}_0 }^{}{{\mathcal{P}}_n } \ $,
 ${\mathcal{P}} \left[{\mathbb{C}}\right] :=  
 \bigcup_{n \in {\mathbb{Z}}_0 }^{}{{\mathcal{P}}_n 
 [{\mathbb{C }}] } \ $. 
If $\varphi : {\mathbb{R}} \to {\mathbb{R}} \ $ then the 
function $M_{\varphi } (x) := \lim\limits_{\delta 
\downarrow 0 } \sup\limits_{y \in (x-\delta ,x+\delta ) } 
\varphi (y) \ $ is called an upper Bair function of 
$\varphi \ $, and  $\delta \downarrow 0 \ $ means 
$\delta \to 0 \ $ and $\delta > 0 \ $.

To shorten expressions the following notations will be used:
$$ I:= [-1,1]; \ I_0 := (-1,1); \ 
J:= (- \infty , -1)\cup (1, + \infty), \ 
I_R := R\cdot I, \ J_R := R\cdot J, \ R>0 \ . $$

\noindent
For nonnegative function $F : {\mathbb{R}} \to {\mathbb{R}}^+ 
:= [0, + \infty ) \ $ let  $S_F := \left\{ \ x \in 
{\mathbb{R}} \  | \ F(x) > 0 \right\} \ $.

Let ${\mathfrak{B}} ({\mathbb{R}}) \ $ denote the family of 
Borel subsets of ${\mathbb{R}} \ $, ${\mathcal{M }} 
({\mathbb{R }} ) \ $ - linear space of finite Borel measures 
on ${\mathbb{R}} \ $ and $L_p (\mu ) := L_p ({\mathbb{R}} ,  
d \mu ) \ $, $\left\|f\right\|^p_{L_p (\mu ) } := 
\int\limits_{{\mathbb{R}}}^{} |f(x)|^p \ d\mu (x) \ $,
$1 \leq p < \infty \ $. It should be reminded 
that every measure $\mu \in {\mathcal{M }} ({\mathbb{R }} ) 
\ $ is regural [6, VI, Def.8.2, Ex.8.16] and therefore for any 
positive $\mu \in {\mathcal{M }} ({\mathbb{R }} ) \ $ and 
arbitrary $A \in {\mathfrak{B}} ({\mathbb{R}}) \ $ there 
exists [6, VI, (8.14)] such sequence of compactly supported 
continuous functions
$\psi_n [{A, \mu }] : {\mathbb{R}} \to [0,1] \  $, 
$n \geq 1 \ $, that
$$ \lim\limits_{ n \to \infty} \left\| \ \chi_A - \psi_n 
[{A, \mu }] \ \right\|_{L_p (\mu )} \ = \ 0 \ \ \forall \ 1 
\leq p < \infty \   . \eqno{(1.1.1)} $$

\noindent
For every $\mu \in {\mathcal{M }} ({\mathbb{R }} ) \ $ 
Hahn expansion of the space $\left({\mathbb{R}} , 
{\mathfrak{B}} ({\mathbb{R}}) \right) \ $ with respect to 
the measure $\mu \ $ will be denoted by ${\mathbb{R}} = 
{\mathbb{R}}^+_{\mu} \sqcup  {\mathbb{R}}^-_{\mu } \ $, 
where $A \sqcup B \ $ denotes union of disjoint sets $A \ $ 
and $B \ $   
[6, I, Th.16.2]. For the expansion of the measure 
$\mu \in {\mathcal{M }} ({\mathbb{R }} ) \ $ in the sense of 
Jordan we will use the following notations: $\mu = \mu_+ - 
\mu_- \ $, $\mu_+ (A) := \mu (A \cap  {\mathbb{R}}^+_{\mu}  ) 
\ $, $\mu_- (A) := \mu (A \cap  {\mathbb{R}}^-_{\mu}  ) \ $, 
$\forall A \in {\mathfrak{B}} ({\mathbb{R}}) \ $, and 
$\left\|\mu \right\| := \left|\mu \right| ({\mathbb{R}}) \ $,
where $\left|\mu \right| := \mu_+ + \mu_- \ $ [6, I.16]. }

\bigskip
\bigskip
{\bf 1.2.  Background. \ } For arbitrary $w: {\mathbb{R}} 
\to [0,1] \ $ consider the seminormed space

$$ C^0_w := \left( \left\{\ f \in C ({\mathbb{R}}) \ | \ 
\lim\limits_{ |x| \to \infty } w(x) f (x) \ = 0 \  
\right\} \ , \ \ \left\|\cdot \right\|_{w} \ \right)
 , \ \eqno{(1.2.1 )}$$

\noindent
where $\left\|f \right\|_{w} := \left\|w \cdot  f  
\right\|_{C(S_w  )} \ $ $\forall \ f \in C^0_w \ $. 
An obvious inequality 
$\left\|f \right\|_{w} \leq \left\|  f  \right\|_{
C ({\mathbb{R}} )} \ $ $\forall \ f \in C^0 ({\mathbb{R}} 
) \ $ implies inclusions $C^0 ({\mathbb{R}} ) 
\subseteq C^0_w \subset C ({\mathbb{R}} ) \  $ and validity 
of the continuous embedding $C^0 ({\mathbb{R}} ) 
\hookrightarrow C^0_w \ $ [13, 0.2.9].
Since [13, IV, Ex.4.45] for every $L \in C^0 ({\mathbb{R}} )^* 
\ $ there exists such $\mu \in {\mathcal{M}} ({\mathbb{R}} 
) \ $ that $L(f) = \int_{{\mathbb{R}}}^{} f(x) \ d \mu(x) \ $
$\forall \ f \in C^0 ({\mathbb{R}}) \ $, indicated embedding 
means that for any element $L \ $ of the Banach space 
$\left(C^0_w \right)^* \ $ with norm 
$\left\| L \right\|_w \ := \sup \{ \left| L(f) \right| 
\  \left| \ f \in C^0_w , \, \left\| f \right\|_w \leq 1  
\right. \} \ $ there exists such 
$\mu_L \in {\mathcal{M}} ({\mathbb{R}} ) \ $ 
that
$$ L(f) = \int_{{\mathbb{R}}}^{} f(x) \ d \mu_L (x) \ 
 \ \forall \ f \in C^0 ({\mathbb{R}}) \ . \eqno{(1.2.2)} $$

In this Chapter, we describe the Banach space $B^0_w \ $ 
being a completion of the normed space $N^0_w \ $ 
associated [13, I.10.2] with seminormed space $ C^0_w \ $, 
give a supplement to \\ S. Mergelyan's Theorem [32, Th.7], 
formulate the weighted analog of the Weierstrass polynomial 
approximation theorem and establish a general form of any 
functional in $\left(C^0_w \right)^* \ $, i.e. in view of 
(1.2.2) find a complete description of the subspace 
$\left\{\mu_L \right\}_{L \in \left(C^0_w \right)^*  } 
\subseteq {\mathcal{M}} ({\mathbb{R}} ) \  $.

\medskip
{\textsc{Remark 1.1. ( Mergelyan's regularity. ) } }
Studying the polynomial approximation problem in $C^0_w \ $  
S. Mergelyan suggested [32] to change the weight function 
$w \ $ by its upper Bair function $M_w \ $. Let us clarify 
what does that suggestion mean in terms of the seminormed 
spaces.
It is known [33] that an upper Bair function $M_w \ $ be an 
upper semicontinuous function [17] and the following 
relations hold:
$$ 0 \leq  w(x) \leq \ M_w (x) \leq 1 \ \ \forall \ x \in 
{\mathbb{R}} \ ; \ \ \ S_w \subseteq S_{M_w } \subseteq 
\overline{S}_w  =\overline{S}_{M_w }  \ , \eqno{(1.2.3)}$$

\noindent
Besides that for any open set $G \subseteq {\mathbb{R}} \ $:
$$ \left\| f \cdot {\chi }_G \right\|_w \ = \       
\left\| f \cdot  \chi_G  \right\|_{M_w } \ \ \forall \ f 
\in C(\overline{S}_w  ) \ . \eqno{(1.2.4)} $$

\noindent
Therefore the seminormed spaces $C^0_w \ $ and $C^0_{M_w } \ $
coincide identically, i.e. $C^0_w  \equiv C^0_{M_w } \ $.
In spite of the available possibility to consider everywhere
below only upper semicontinuous functions $w \ $, i.e. $w = M_w \ $,
we will not do so and will examine a general case  
$w \neq  M_w \ $ using notation:
$$ h:= M_w \ . $$

\bigskip
\bigskip
{\bf 1.3.  Associated normed space $ N^0_w  \ $. \ }
One can easily conclude from known criterion [13, 1.10.1]
of the separability of locally convex spaces and from the 
continuity of functions in $C^0_w  \ $ that seminormed space 
$C^0_w  \ $ is a normed one if and only if $\overline{S}_w = 
{\mathbb{R}} \ $.

Denote by $N^0_w := C^0_w \setminus N_{C^0_w } \ $ (see 1.1)
the normed space associated with $C^0_w  \ $. Introduce
the normed space
$$ C^0_w (\overline{S}_w  ) := \left(\left\{ \ f \in 
C(\overline{S}_w  ) \ | \ \lim\limits_{x \in \overline{S}_w  ,
\ |x| \to \infty  } w(x) f(x) \ = \ 0 \ \right\}, \ 
\left\|\cdot \right\|_w \ \right) \ \eqno{(1.3.1)}$$

\noindent
and corresponding two normed ones of the restrictions:
$ C^0_w (\overline{S}_w  )\!\restriction_{S_h } := 
\left( \ \left\{ \ \  f\!\restriction_{S_h }   \ \  |  
\ \ \ f \in \right.\right. \ $ \\ $\in \left.\left. C^0_w  
(\overline{S}_w  )  \right\},  \left\|\cdot 
\right\|_w  \right) \ $ , 
$  C^0_w (\overline{S}_w  )\!\restriction_{S_w } := 
\left(\left\{  f\!\restriction_{S_w }   |  
f \in C^0_w (\overline{S}_w  )  \right\},  
\left\|\cdot \right\|_w  \right)  $. Due to (1.2.1),
(1.2.4) $C^0_w (\overline{S}_w  ) \equiv  C^0_h 
(\overline{S}_h  ) \ $ and  hence, 
$C^0_w (\overline{S}_w  )\!\restriction_{S_h } \equiv 
C^0_h (\overline{S}_h  )\!\restriction_{S_h } \ $. It's 
evident, that transformations $f \to f\!\restriction_{
S_h } \  $, $f \to f\!\restriction_{S_w } \  $ $\forall 
\ f \in  C^0_w (\overline{S}_w  ) \ $ determine isometric 
relations $C^0_w (\overline{S}_w  ) \simeq C^0_h 
(\overline{S}_h  )\!\restriction_{S_h } \ $ and $C^0_w 
(\overline{S}_w  ) \simeq C^0_w (\overline{S}_w  )
\!\restriction_{S_w } \ $, respectively. Besides that 
defined by formula $V (\pi (f)) = 
f\!\restriction_{\overline{S}_w  } \ $ $\forall \ f \in  
C^0_w \ $ transformation
$$ V \ : \ N^0_w \ \to \  C^0_w (\overline{S}_w  ) \ , 
\eqno{(1.3.2)} $$

\noindent
determines linear isometry (see 1.1 and [3, IV.1.3 ]) of the
spaces $N^0_w \ $ and $C^0_w (\overline{S}_w  ) \ $. Really,
equality $\left\|V (\pi (f)) \right\|_w = \left\| \pi (f) 
\right\|_{ N^0_w } \equiv \left\|f\right\|_w \ $ is obvious
and relation $V(N^0_w  ) = C^0_w (\overline{S}_w  ) \ $ 
follows from known continuity [33, IV.4, Lemma 2] of the 
linear extension $f_c \ $ to the whole real line of some 
continuous on the closed set $F \subset {\mathbb{R}} \ $ 
function $f \ $ . In addition if [4, IV.5, Th.21] interval 
of the kind $(- \infty , b) \ $ or $(a, + \infty ) \ $ is a 
part of ${\mathbb{R}} \setminus F \ $ then we will regard:
$$\left\{\begin{array}{ll}
f_c (b - \theta ) := \theta f(b)  & \forall  \ \theta \in 
(0,1) ; \\
f_c (b - \lambda ) := 0   & \forall \ \lambda \ \geq 1 \ ; 
  \end{array}\right.
\left\{\begin{array}{ll}
f_c (a + \theta ) := \theta f(a)  & \forall  \ \theta \in 
(0,1) ;  \\
 f_c (a + \lambda ) := 0  & \forall \ \lambda \ \geq 1 \ ;
  \end{array}\right. \eqno{(1.3.3)} $$
  
\noindent  
respectively. Therefore
$$ N^0_w \cong C^0_w (\overline{S}_w  ) \equiv  C^0_h 
(\overline{S}_h  ) \cong C^0_h (\overline{S}_h  
)\!\restriction_{S_h } \cong C^0_w (\overline{S}_w 
)\!\restriction_{S_w } \ , \eqno{(1.3.4)} $$

\noindent
i.e. associated with $C^0_w \ $ normed space $N^0_w \ $ can 
be identified with arbitrary indicated in (1.3.4) isometric 
normed spaces.

\bigskip
\bigskip
{\bf 1.4.   Banach space $ B^0_w \ $ . \ }

\bigskip
{\textsc{Definition 1.4.1. \ }} 
{\textit{Let }} $w : {\mathbb{R}} \to [0,1] \ $, $h:= M_w \ $
  {\textit{ is an upper Bair function of $w \ $ and }} $S_h := 
 \left\{   x \in {\mathbb{R}} \ | \ h(x) > 0 \ \right\} \ $.
{\textit{The space  }} $B^0_w \ $ {\textit{ is called a 
Banach space associated with the seminormed space  }} 
$C^0_w \ $ {\textit{if  }} $B^0_w \ $ {\textit{ is 
equipped 
with norm  }} $\left\|f\right\|_h := \sup\limits_{ x \in S_h }
h(x) |f(x)| \ $ {\textit{and consists of all functions }}
$f : S_h \to {\mathbb{R}} \ $, {\textit{which satisfy the 
following three conditions:}}

\medskip
(1.4.1) {\textit{function}} $f \ $ {\textit{is  continuous on
the set}}
$E_{1 \left/ \right. \delta } (h) := h^{-1} ([\delta , 1 ] )= 
\ $ \\     
$ =  \left\{ x \in {\mathbb{R}} \ | \ h(x) \geq \delta 
\right\}  \ $
{\textit{for any }} $\delta \in (0,1] \ $;

\medskip
(1.4.2) $\lim\limits_{ h(x) \to 0 } h(x) f(x) = 0 \ $, i.e.
$$ \forall \varepsilon > 0 \ \exists \ \delta > 0 \ : \ 
\left\{x \in {\mathbb{R}} \ | \  0 < h(x) < \delta \ \right\}
\subseteq \left\{x \in S_h  \ | \   h(x) |f(x) | < \varepsilon
\ \right\} \ ; $$  

\medskip
(1.4.3) $h(x) f(x) \to 0 \ $, $x \in S_h \ $, $|x| \to \infty 
\ $.

\bigskip
\bigskip
If $\lim\limits_{ |x| \to \infty } h(x) =
0 \ $, then property (1.4.2) implies (1.4.3) and in that case
property (1.4.3) can be excluded from Definition 1.4.1. 

Verify now correctness of the Definition 1.4.1. 

Since $h \ $ is an upper semicontinuous function then [17]
all sets $E_{1 \left/ \right. \delta } (h) \ $ for $\delta 
\in (0,1] \ $, are closed and so for any $f \in B^0_w \ $
and $\varepsilon = 1 \ $ one can find such $R > 0 \ $ in 
(1.4.3) and $\delta > 0 \ $ in (1.4.2) that $\left\|h \cdot 
f \right\|_{C (J_R ) } \leq 1 \ $, $\left\|h \cdot f 
\right\|_{C ( h^{-1} (0, \delta  )  ) } \leq 1 \ $, getting on
the supplement compact $S_h \setminus \left[J_R \cup  h^{-1} 
(0, \delta  ) \right] = I_R \cap E_{1 \left/ \right. \delta }
(h) \ $ the uniformly boundedness $f \ $ by property (1.4.1).
That's why $\left\|f\right\|_h < \infty \ $ $ \forall \ f \in 
B^0_w \ $.

If now $\left\{ f_n\right\}_{n \geq 1 } \subset B^0_w \ $ is 
a fundamental sequence in $B^0_w \ $ then by known scheme 
[4, V.5] one can easily obtain an existence of such $F : 
S_h \to {\mathbb{R}} \ $ that $\lim_{n \to \infty } 
\left\|f_n - F  \right\|_h = 0 \ $. Function $F \ $ 
obviously satisfies conditions (1.4.2) and (1.4.3). 
Since for every $m \geq 1 \ $
and $x \in E_m (h) \ $ : $h(x) \geq \frac{1}{m} \ $, then 
$\left\|f_n - F \right\|_h \geq \frac{1}{m} \left\|f_n - 
F \right\|_{C(E_m (h) ) } \ $, i.e. the sequence 
$f_{n}\!\restriction_{E_m (h) } \in 
C(E_m (h) ) \ $, $n \geq 1 \ $, uniformly on the set 
$E_m (h) \ $ 
converges to $F\!\restriction_{E_m (h) } \ $ and 
therefore [4, IV.2] $F\!\restriction_{E_m (h) } \in C(E_m (h )
) \ $. That's why $F \in B^0_w \ $ and introduced in 
Definition 1.4.1 normed space be in fact Banach one. 

It should be noted at last that as well as $C^0_w \ $ 
(see Remark 1.1 ) Banach space $B^0_w \ $ posseses the 
property
$B^0_w \equiv B^0_{M_w } \ $. Let now formulate the basic
result of that section.

\bigskip
\bigskip
{\textsc{Theorem 1.1. \ }} 
{\textit{ Let $w: {\mathbb{R}} \to [0,1] \ $, $h := M_w \  $
is an upper Bair function of $w \ $ and 
$S_h := \left\{   x \in {\mathbb{R}} \ | \ h(x) > 0 \ 
\right\} \ $. Linear operator }}
$$ T : C^0_w \ \to \ B^0_w \ , $$

\noindent
{\textit{defined by formula }}
$$ Tf = f\!\restriction_{S_h } \ \ \forall \ f \in C^0_w 
\eqno{(1.4.4)} $$

\noindent
{\textit{isometrically and tightly embeds seminormed space 
$C^0_w \ $
into the Banach space $B^0_w \ $, i.e. }}

\medskip
(1.4.5a) $ \left\|Tf \right\|_h = \left\|f \right\|_w \ $
$ \forall \ f \in C^0_w \ $ ; 

\medskip
(1.4.5b) $ T(C^0_w ) \ $ {\textit{is a dense subspace of the
Banach space $B^0_w \ $. }}

\medskip
\noindent
{\textit{In addition $ T(C^0_w ) \ $ coincides with the 
subspace of those functions $f \in B^0_w \ $, which 
can be extended to the continuous on 
$\overline{S}_h \ $ function. }}

\bigskip
\bigskip
Validity  of the following implication (see 1.7.2):
$$ \begin{array}{ll}
\exists \left\{x_n \right\}_{n \geq 0 } \subseteq S_h \ : & 
\lim\limits_{n \to \infty } x_n = x_0 \in {\mathbb{R}} , \ 
\lim\limits_{x_n \to x_0 } h(x_n ) = 0 \ \Rightarrow \ \\
\Rightarrow \ \exists F \in B^0_w \ :  & \ \left\|F 
\right\|_{C (S_h \cap (x_0 - \delta ,  x_0 + \delta ) ) }
\ = \ + \infty \ \ \forall \ \delta \ > \ 0 \ 
  \end{array} \eqno{(1.4.6)} $$

\noindent
allows us to characterize the following partial cases of  
Theorem 1.1.

\bigskip
\bigskip
{\textsc{Corollary 1.1. \ }}
{\textit{Let $w: {\mathbb{R}} \to [0,1] \ $, $h := M_w \  $. 
Then}

\medskip
(1.4.7) $C^0_w \ $ {\textit{is a normed space if and only if}}
$\overline{S}_h  =  {\mathbb{R}} \ $;

\medskip
(1.4.8) {\textit{The following statements are equivalent:
$$\begin{array}{ll}
(1.4.8a) \ C^0_w \  \ \mbox{is a Banach space }  \ ; &
(1.4.8b) \ B^0_w = C^0_w \ ; \\
(1.4.8c) \ S_h = {\mathbb{R}} \ \mbox{and } \ B^0_w 
\subseteq C({\mathbb{R}} )\ ;  &(1.4.8d)\  
\inf\limits_{x \in [-R, R] } 
h(x) \ > \ 0 \ \  \ \ \forall \ R > 0 \ ; 
  \end{array} $$ }}

 \medskip
(1.4.9) {\textit{The following statements are equivalent:
$$\begin{array}{ll}
(1.4.9a) \ N^0_w \  \ \mbox{is a Banach space }  \ ; 
&(1.4.9b)
\ B^0_w = C^0_w (\overline{S}_w) \ ; \\
(1.4.9c) \ S_h = \overline{S}_h \ \mbox{and } \ B^0_w 
\subseteq C( \overline{S}_h )\ ;  &(1.4.9d) \  
\inf\limits_{x \in S_h \cap [-R, R] }  h(x) \ > \ 0 \ \ \  
\forall \ R > 0 \    ; \end{array} $$ 
  }}

\noindent
{\textit{where}} $\inf\limits_{} \emptyset := + \infty \ $.

\bigskip
\bigskip
The following application of the Theorem 1.1 gives some 
explanation why everywhere above we have not assumed the upper 
semicontinuity of $w \ $.

\bigskip
\bigskip
{\textsc{Corollary 1.2. \ }}
{\textit{Let $w: {\mathbb{R}} \to [0,1] \ $ and 
${\mathcal{M}} \ $ is some dense subset of the seminormed 
space $C^0_w \ $. Function $f : S_w \to {\mathbb{R}} \ $ 
can be approximated by elements of ${\mathcal{M}} \ $, i.e.
$$\forall \varepsilon > 0 \ \exists \ m_{\varepsilon } \in 
{\mathcal{M}} \ : \ \ w(x) \left|f(x) - m_{\varepsilon  }  
(x) \right| < \varepsilon \ \ \forall \ x \in S_w \ 
\eqno{(1.4.10)} $$
}}

\noindent
{\textit{if and only if \ \ $\exists g \in B^0_w : \ 
g\!\restriction_{S_w } = f \ $. }}

\bigskip
\bigskip
\bigskip
{\bf 1.5.  Supplement to S.N. Mergelyan's Theorem  . \ }
In [32, Th.7] S. Mergelyan proved that for the weight function 
$ w : {\mathbb{R}} \to [0,1] \ $ satisfying  condition

$$ \left\| x^n w \right\|_{C({\mathbb{R}} )} < + \infty \ 
\ \  \forall \ n \in {\mathbb{Z}}_0 \ , \eqno{(1.5.1)} $$

\noindent
either algebraic polynomials ${\mathcal{P}} \ $ are dense in 
$C_w^0 \ $ or they can approximate only those functions $f : 
S_w \ \to \ {\mathbb{R}} \ $ which can be extended 
from their domain of definition $S_w  \ $ into the whole 
complex plane as an entire function 
of minimal exponential type. I. Hachatryan  [14] gived the 
description of the indicated in the Mergelyan's theorem class
of entire functions. Corollary 1.2 implies the following 
supplement to the Mergelyan's theorem when algebraic 
polynomials ${\mathcal{P}} \ $ are dense in $C^0_w \ $.

\bigskip
\bigskip
{\textsc{Theorem 1.2. }} \ {\textit{Let $ w : 
{\mathbb{R}} \to [0,1] \ $ satisfies conditions (1.5.1),
$M_w \ $ be an upper Bair function of $w \ $ and algebraic 
polynomials ${\mathcal{P}} \ $ are dense in $C^0_w \ $.    
   Then the function $f : S_w \ \to \ {\mathbb{R}}
\ $ can be approximated by polynomials, i.e. 
$$ \exists \left\{P_n \right\}_{n \geq 1 } \ \subset \ 
{\mathcal{P}} \ : 
\lim_{n \to \infty } \sup\limits_{x \in S_w } w(x)\left| P_n 
(x) - f(x) 
\right| \ = \ 0 \ , $$ }}

\noindent
{\textit{ if and only if  that function can be extended into 
the set $S_{M_w} \ $ as a function \\ $f : S_{M_w } \ \to \ 
{\mathbb{R}} \ $, satisfying the 
following conditions: }}

\medskip
\noindent
(1.5.2) {\textit{for every $m \geq 1 \ $  $ f \ $ is a 
continuous function
on the closed set }} \\ 
$\left\{x \in {\mathbb{R}} \ | \ M_w ( x )  \geq \frac{1}{m} 
\right\} \ ; $

\medskip
\noindent
(1.5.3) $ \lim_{M_w (x) \to 0 } \ M_w (x) \cdot   f ( x ) \ 
= \ 0 \ $, i.e. 
$\forall \ \varepsilon > 0 \ \exists \ \delta \ > 0 \ : \ $
$$\left\{x \in {\mathbb{R}} \ | \ 0  < \ M_w ( x ) < \delta    
\right\} \ \subseteq \ \left\{x \in S_{M_w } \ | \  
M_w ( x ) \cdot | f(x) | \ < \varepsilon  \right\} \ .  $$

\bigskip
If $S_w \ $ is a bounded set then   
 conditions (1.5.1) are obviously true and by
Weierstrass approximation theorem algebraic polynomials  
${\mathcal{P}} \ $  are dense in  $C^0_w \ $. That's why for 
arbitrary weight $ w : {\mathbb{R}} \to [0,1] \ $ with 
bounded set $S_w \ $ conditions (1.5.2) and (1.5.3) give 
the {\it{weighted analog of the Weierstrass polynomial 
approximation theorem}}.
It is interesting to note also that for the weight 
$w(x)=  \sqrt{1- x^2 } \cdot \chi_{[-1, 1] } (x) \ $
conditions (1.5.2) and (1.5.3) are equivalent to
$f \in C ((-1, 1)) \ $ and  $\lim\limits_{|x| \to 1 } 
\sqrt{1 - x ^2 } f(x) \   = 0 \ $. This fact is known and can
be found in [28] where according to these two 
conditions the subspaces of known spaces $B^r \ $ were 
introduced.  

 \bigskip
\bigskip
{\bf 1.6. General form of the functionals in  $ \left( C^0_w 
\right)^* \ $ . \ } 
To prove the main theorem of that section the following 
version of known M. Krein's lemma will be necessary [37].

\bigskip
\bigskip
{\textsc{Lemma 1.1. }} \ {\textit{Let $\left(X, p \right) \ $
be a seminormed space and $K \subset X \ $ is a normal cone, 
i.e. convex set $K \ $ satisfies: $\lambda \cdot K  \subseteq 
K \ $ $\forall \lambda \geq 0 \ $ and $p(x) \leq 
p( x + y ) \ $ $\forall \ x,y \in K \ $. If $X^* \ $ is a 
Banach space conjugate to $X \ $ then for $K^* := \left\{ 
x^* \in X^* \ | \ x^* (x) \geq 0 \ \ \forall x \in K \ \right\} 
\ $ the following equality holds: }}
$$ K^* \ - \  K^*  \ =  \ X^* $$

\bigskip
It should be noted that in [5] a necessary and sufficient 
condition for the validity of more general equality $(K_1 
\cap K_2 )^* = K_1^* + K_2^* \ $ has been established and a 
notion of normal pair of cones $(K_1 , K_2 ) \ $ of 
transfinite order $\alpha \ $ has been introduced. Now we 
can formulate the main theorem of this section.

\bigskip
\bigskip
{\textsc{Theorem 1.3. }} 
{\textit{Let  $w : {\mathbb{R}} \to [0,1] \ $, $ M_w \ $
  {\textit{is an upper Bair function of}} $w \ $ and $S_{M_w } 
 := \left\{   x \in {\mathbb{R}} \ | \ M_w (x) > 0 \ \right\} 
 \ $.  If $L \ $ is a linear continuous functional on the 
 seminormed space $C^0_w \ $ then there exists such
 measure $\mu \in {\mathcal{M } } ({\mathbb{R}} ) \ $,
 that $|\mu | ({\mathbb{R}} \setminus  S_{M_w } ) = 0  \ $
 and }}
 $$ L( f ) = \int\limits_{{\mathbb{R } } }^{} M_w (x) f(x)\ 
 d \mu (x) \ \  \forall \ f \in C^0_w \ . \eqno{(1.6.1)} $$
 
 \noindent
 {\textit{ For arbitrary $\mu \in {\mathcal{M } } 
 ({\mathbb{R}} ) \ $ defined by formula (1.6.1) functional 
 $L \ $ is linear and continuous on the 
 seminormed space $C^0_w \ $ with $\left\|L \right\| = |\mu | 
 (S_{M_w } ) \ $ (see 1.1). }}
  
 \bigskip
 \bigskip
\bigskip
{\bf 1.7. Proofs. \ }

\medskip
{\textsc{ 1.7.1. Proof of Theorem 1.1. \ }}
{\normalsize{   Equality (1.4.5a)
follows from (1.2.4). Prove now that $T(C^0_w ) \subseteq 
B^0_w \ $, where $T(C^0_w ) \equiv \left\{ f\!\restriction_{
S_h } \ | \ f \in C^0_w \ \right\} =: C^0_w \!\restriction_{
S_h } \ $. If $g = f\!\restriction_{S_h } \ $, $f \in 
C^0_w \ $, then conditions (1.4.1) and (1.4.3) for function 
$g \ $ are realized. Let us prove that $g \ $ satisfies 
property (1.4.2). For given $\varepsilon > 0 \ $
one can find by (1.4.3) such $R(\varepsilon ) > 0 \ $ that: 
$\left\|h g \right\|_{C(S_h \cap J_{R(\varepsilon ) } ) }  
< \varepsilon \ $. Denote $C(\varepsilon ):= \left\| g 
\right\|_{C(S_h \cap I_{R(\varepsilon ) } ) } \leq \left\|f 
\right\|_{C(\overline{S}_h \cap I_{R(\varepsilon ) } ) } < 
\infty \ $. Then the number $\delta (\varepsilon ) := 
\varepsilon / C(\varepsilon ) > 0 \ $ in view of 
$\left\|h g \right\|_{C(h^{-1} (0, \delta (\varepsilon ) ) 
\cap I_{R(\varepsilon ) } ) } \leq \delta (\varepsilon ) 
C(\varepsilon ) = \varepsilon \ $ will be required for the 
validity of (1.4.2), i.e. $g \in B^0_w \ $. In addition 
proved in 1.3 equality $C^0_w \!\restriction_{\overline{S}_h }
= C^0_w (\overline{S}_h ) \ $ yields
$$ T (C^0_w ) \equiv C^0_w \!\restriction_{S_h }
= C^0_w (\overline{S}_h )\!\restriction_{S_h } \subseteq 
B^0_w \ , \eqno{(1.7.1.1)}$$

\noindent
and to finish the proof it is remained to show that $C^0_w 
\!\restriction_{S_h } \ $ is dense in $B^0_w \ $.

Consider an arbitrary $f \in B^0_w \ $, $\varepsilon > 0 \ $
and prove that there exists such $f_{\varepsilon } \in C^0_w 
\ $ that $\left\|f - f_{\varepsilon } \right\|_h \leq 2 
\varepsilon  \ $. Property (1.4.2) admits to find such 
positive integer $m \geq 1 \ $ that
$$ \left\{ \begin{array}{ll}
\left\{ \ x \in {\mathbb{R}} \ | \ 0 < h(x) < \frac{1}{m} \ 
\right\}   
 & \subseteq \left\{ x \in S_h \ | \ h(x)|f(x)| < \varepsilon
  \right\} \ ;  \\
 m > \frac{\left\|f\right\|_h }{\varepsilon } \ , \ \ E_m 
 \neq \emptyset \ , & 
   \end{array}\right. \eqno{(1.7.1.2)} $$

\noindent
where (see (1.4.1)) $E_{p } := E_{p} (h) \ $ $\forall \ p 
\geq 1 \ $.

Since $f \ $ is a continuous function on the set $E_{m^2 } \ $,
then by means of indicated in 1.2 method we extend $f \ $ into 
the whole real line obtained $f_{\varepsilon } \in C^0_w \ $. 
Let us prove that
$$ \left\|f - f_{\varepsilon } \right\|_h \leq 2 \varepsilon \ .
\eqno{(1.7.1.3)} $$

Since for $x \in S_h \setminus E_{m^2 } \ $:
$0 < h(x) < \frac{1}{m^2 } \leq \frac{1}{m} \ $, then by 
(1.7.1.2): $\left\|f - f_{\varepsilon } \right\|_h = 
\left\|h (f - f_{\varepsilon } ) \right\|_{C(S_h \setminus 
E_{m^2 } ) } \leq \varepsilon + \left\|h f_{\varepsilon } 
\right\|_{C( S_h \setminus E_{m^2 } ) } \ $ and therefore for the
validity of (1.7.1.3) it is sufficient to prove that
$$ \left\|h f_{\varepsilon } \right\|_{C( {\mathbb{R}} 
\setminus  E_{m^2 } )} \ \leq \ \varepsilon \ . 
\eqno{(1.7.1.4)} $$

\noindent
Inequality (1.7.1.4) is trivial if $E_{m^2} = {\mathbb{R}} 
\ $. But if ${\mathbb{R}} \setminus E_{m^2 } \neq \emptyset 
\ $  then ${\mathbb{R}} \setminus E_{m } \neq \emptyset \ $ 
and by (1.7.1.2) ${\mathbb{R}} \setminus E_{m} \neq 
{\mathbb{R}} \ $. Consider an arbitrary $x \in  {\mathbb{R}}
\setminus E_{m^2 } \subseteq  {\mathbb{R}} \setminus E_{m} \ $. 
Such $x \ $ belongs to one of the forming [4, IV.5, 
Th.21]
open set ${\mathbb{R}} \setminus E_{m} \ $  interval $(a, b) 
\ $ and moreover $x \in (a_k , b_k ) \subseteq (a,b) \setminus
E_{m^2 } \ $, where $(a,b) \setminus E_{m^2 } = \bigsqcup_{k 
= 1}^{Q } (a_k , b_k ) \ $, $a_k , b_k \in E_{m^2 } \cup 
\left\{\pm \infty \right\} \ $, $1 \leq Q \leq \infty \ $.

Assume that $(a_k , b_k ) \ $ is a bounded interval. Then 
$\exists \theta \in (0,1) : $ $x = \theta a_k + (1 - \theta )
b_k $, and
$$ h(x) |f_{\varepsilon } (x) | \ \leq \ \theta h(x) |f(a_k )
| + (1 - \theta ) h(x) |f(b_k ) | \ . \eqno{(1.7.1.5)} $$

\noindent
If $a_k > a \ $, then $a_k \in {\mathbb{R}} \setminus E_{m} 
\ $, $h(x) < \frac{1}{m^2 } \leq  h(a_k ) < \frac{1}{m} \ $ 
and by (1.7.1.2) $h(x) |f(a_k )| \leq h(a_k ) |f(a_k ) | 
\leq \varepsilon \ $. If $a_k = a \ $ then $a_k \in E_m \ $ 
and $h(x) < \frac{1}{m^2 } \leq  \frac{h(a_k )}{ m } \ $, 
whence in view of (1.7.1.2): $h(x) |f(a_k )| \leq 
\frac{1}{m} h(a_k ) |f(a_k ) | \leq \frac{\left\|f\right\|_h 
}{m } \leq \varepsilon \ $. Performing the same estimate for 
$h(x) |f(b_k ) | \ $ we obtain from (1.7.1.5) $h(x) 
|f_{\varepsilon } (x) | \leq \varepsilon \ $.

Assume now that $(a_k , b_k ) = (a_k , + \infty) \ $. Then 
by (1.3.3) for $x = a_k + \theta \ $, $\theta \in (0,1) \ $:
$h(x)|f_{\varepsilon } (x)| = \theta h(x) |f(a_k ) | \ $, and
for $x \geq a + 1 \ $ : $f_{\varepsilon } (x) = 0 \ $. Since
$E_{m} \neq \emptyset \ $, then $(a, b ) \neq  
{\mathbb{R}} \ $ and hence, in that case $(a, b) = 
(a, + \infty ) \ $. Estimated $h(x) |f (a_k ) | \ $ as well 
as it has been done above for the case of bounded interval 
$(a_k , b_k ) \ $ we will get  $h(x) |f_{\varepsilon } 
(x) | \leq \varepsilon \ $ again.

The case $(a_k , b_k ) = (- \infty , b_k ) \ $ can be 
considered in just the same way. That's why inequality 
(1.7.1.4) together with Theorem 1.1 is proved.

}}

\bigskip
\bigskip
{\textsc{ 1.7.2. Proof of implication (1.4.6). \ }}
{\normalsize{  Without loss of generality one can consider 
that the sequences $\frac{1}{\lambda_n  } := h(x_n ) > 0 \ $,
and $|x_n - x_0 | \ $, $n \in {\mathbb{Z}}_0 \ $, are 
decreasing. Let as in 1.7.1 $E_{\lambda } := E_{\lambda } (h)
\ $ for
$\lambda \in [1, + \infty ) \ $. Since $x_{n+1} \in 
{\mathbb{R}} \setminus E_{\lambda_{n} } \ $ $\forall n \in 
{\mathbb{Z}}_0 \ $, it is possible to find such sequence
$\left\{\delta_n \right\}_{n \geq 1 } \ $ of positive real 
numbers that $x_{n +1 } + \delta_{n + 1 } I_0 \ \subset 
{\mathbb{R}} \setminus E_{\lambda_{n} } \  $ $\forall n 
\in {\mathbb{Z}}_0 \ $ and the sets $\left\{x_n + \delta_n I  
\right\}_{ n \geq 1 } \ $ are disjoint. Set
$$ F(x) := \sum_{k \geq 1 }^{} \sqrt{\lambda_{k-1} } 
\alpha_k (x) \ , \ \  \alpha_k (x) = \left(1 - \left|
\frac{x - x_k }{\delta_k }\right| \right) \chi_I \left(
\frac{x - x_k }{\delta_k } \right) \ , k \geq 1 \ . 
\eqno{(1.7.2.1)}$$ 

\noindent
Equalities $F(x_k ) = \sqrt{\lambda_{k - 1 } } \ $ $\forall k
\geq 1 \ $ imply $\left\|F\right\|_{C(S_h \cap (x_0 - \delta ,
x_0 + \delta ))} = + \infty \ $ $\forall \delta > 0 \ $.
It remains to prove that $F \in B^0_w \ $. Since 
$x_{n + 1 + p } + \delta_{n + 1 + p } I_0 \subseteq 
{\mathbb{R}} \setminus E_{\lambda_{n + p} }  \subseteq 
{\mathbb{R}} \setminus E_{\lambda_{n } } \ $ $\forall n, p 
\in {\mathbb{Z}}_0 \ $ then for every $n \geq 1 \ $:
$$  F(x)\!\restriction_{E_{\lambda_{n }}} = 
\sum_{k = 1 }^{n}  \sqrt{\lambda_{k-1} } \alpha_k (x)
\!\restriction_{ E_{\lambda_{n }}} \ \in \ C( E_{\lambda_n }) 
\ ,  $$
 
\noindent
and consequently, property (1.4.1) for $F \ $ is fulfilled.
Validity of (1.4.3) is obviuos. Let us show that $F \ $ 
satisfies (1.4.2). Really for any $n \geq 1 \ $ inequalities
$0 < h(x) < \frac{1}{\lambda_n } \ $ yield
$$ h(x)F(x) \leq \sum_{k \geq 1 }^{} \alpha_k (x) \min 
\left\{\frac{1}{\sqrt{\lambda_{k-1 }  } } \ , 
\ \frac{\sqrt{\lambda_{k-1 }  } }{\lambda_n }\right\} \leq
\frac{1 }{\sqrt{\lambda_n } } \ , $$ 

\noindent
what means the validity of (1.4.2). That is why $F \in 
B^0_w \ $ and implication (1.4.6) is proved.
}}

\bigskip
\bigskip
{\textsc{ 1.7.3. Proof of Corollary 1.1. \ }}
{\normalsize{ Correctness of (1.4.7) was proved in 1.2.  
Implications (1.4.8b)$\Rightarrow$(1.4.8a),
(1.4.8c)$\Rightarrow$(1.4.8b),
(1.4.9c)$\Rightarrow$(1.4.9b)
are evident and (1.4.9b)$\Rightarrow$(1.4.9a) follows 
from (1.3.4).

{\textit{(1.4.8a)$\Rightarrow$(1.4.8c). }}
Since $C^0_w \ $ is a normed 
space then by (1.4.7) ${\mathbb{R}} = \overline{S}_w = 
\overline{S}_h \ $ (see (1.2.3)). Besides that by Theorem 1.1
$T (C^0_w ) = C^0_w \!\restriction_{S_h } =  B^0_w \ $.
Then assumption $\overline{S}_h \setminus S_h \neq 
\emptyset \ $ together with upper semicontinuity $h \ $ and
(1.4.6) leads to a contradiction. That' why 
${\mathbb{R}} = \overline{S}_h = S_h \ $ and $B^0_w = C^0_w 
\subset C({\mathbb{R}} ) \ $.

{\textit{(1.4.8c)$\Rightarrow$(1.4.8d). }}
If $\exists R > 0 \ $:
$\inf_{x \in [-R, R] } h(x) = 0 \ $, then (1.4.6) yields
$S_h \neq {\mathbb{R}} \ $ or $B^0_w \setminus C({\mathbb{R} }
) \neq \emptyset \ $ in the case $S_h = {\mathbb{R}} \ $.

{\textit{(1.4.8d)$\Rightarrow$(1.4.8c). }}
If for any $R > 0 \ $:
$1/ \lambda (R) := \inf_{x \in I_R } h(x) > 0 \ $, then
$S_h = {\mathbb{R}} \ $ and $I_R \subseteq E_{\lambda (R)}
(h) \ $ $\forall R > 0 \ $. Therefore for arbitrary $f \in 
B^0_w \ $: $f \in C(I_R ) \ $ $\forall R > 0 \ $, and hence,
$B^0_w \subset C({\mathbb{R}}) \ $.

{\textit{(1.4.9a)$\Rightarrow$(1.4.9c). }}
Due to (1.3.4) and (1.7.1.1):
$B^0_w = C^0_w (\overline{S}_h )\!\restriction_{S_h } \ $. An 
upper semicontinuity of $h \ $ and (1.4.6) give 
$\overline{S}_h = S_h \ $ and hence, $B^0_w = C^0_w 
(\overline{S}_h )\subseteq C(\overline{S}_h ) \ $.

{\textit{(1.4.9c)$\Rightarrow$(1.4.9d). }}
If $\exists R > 0 \ $:
$S_h \cap I_R \neq \emptyset \ $  and  $\inf_{x \in 
S_h \cap  I_R } h(x) = 0 \ $, then $S_h \neq 
\overline{S}_h \ $ or
$B^0_w \setminus C(\overline{S}_h ) \neq \emptyset \ $ in 
the case $S_h = \overline{S}_h \ $ by the property (1.4.6).

{\textit{(1.4.9d)$\Rightarrow$(1.4.9c). }}
Let $1/ \lambda (R) := 
\inf_{x \in I_R \cap S_h } h(x) > 0 \ $ $\forall R \geq R_0 
\ $, where $S_h \cap I_{R_0 } \neq \emptyset \ $. Then for 
such
values of $R \ $: $I_R \cap S_h \subseteq E_{\lambda (R) } 
(h) \ $ and by the closure of $E_{\lambda (R) } (h) \ $:
$I_R \cap \overline{S}_h \subseteq E_{\lambda (R_1 ) } (h)
\subseteq  S_h \ $ $\forall R_1 > R \geq R_0 \ $. Therefore
$S_h = \overline{S}_h \ $ and for any $f \in B^0_w \ $:
$f \in C(I_R \cap \overline{S}_h ) \ $ $\forall R > R_0 \ $.
That is why $B^0_w \subseteq C(\overline{S}_h ) \ $ and 
Corollary 1.1 is proved.
}}

\bigskip
\bigskip
{\textsc{ 1.7.4. Proof of Corollary 1.2. \ }}
{\normalsize{

{\it Necessity.} \ If some $f: S_w \to {\mathbb{R}} \ $
satisfies (1.4.10) then for $\varepsilon = \frac{1}{p} \ $,
$p \geq 1 \ $, we obtain from there fundamental sequence 
$\left\{m_{1/p } \right\}_{p \geq 1 } \subseteq C^0_w \equiv 
C^0_h \ $, $h := M_w \ $, being mapped by transformation 
(1.4.4) onto the fundamental sequence in Banach space 
$B^0_w 
\ $ whose limit $g \in B^0_w \ $ due to (1.2.3) satisfies:
$g\!\restriction_{S_w } = f \ $.

{\it Sufficiency.} \ Assume that $\exists g \in B^0_w $:
$g\!\restriction_{S_w } = f \ $. Then by Theorem 1.1 the set
${\mathcal{M}}\!\restriction_{S_h } \equiv T({\mathcal{M}})
\ $ will be dense in Banach space $B^0_w \ $. That is why 
those elements of ${\mathcal{M}}\!\restriction_{S_h } \ $
which approximate $g \ $ in $B^0_w \ $ in view of (1.2.3)
will approximate $f \ $ on $S_w \subseteq S_h \ $ in the 
sense of (1.4.10). Corollary 1.2 is proved.
}}

\bigskip
\bigskip
{\textsc{ 1.7.5. Proof of Lemma 1.1. \ }}
{\normalsize{
It is easy to verify that the closure $\overline{K} \ $ in 
$X \ $ will be a normal cone  and $(\overline{K} )^* = 
K^* \ $. That is why we may consider that $K \ $ is a 
closed cone.

Following well-known scheme of [37, I, Ex.2a] let us examine
any subspace $Y \subset X \ $ which is an algebraic 
complementary subspace to $N:= \left\{x \in X \ | \ p(x) = 0 \ 
\right\} \ $ and for arbitrary $x \in X \ $ in its 
representation $x = n + y \ $, $n \in N \ $, $y \in Y \ $, 
denote $P_Y x := y \ $. It follows from the closure of $K \ $
and an obvious equality
$$ p(x + n ) = p(x) \ \ \forall \  n \in N \ \ \forall x \in X
\ , \eqno{(1.7.5.1)}$$

\noindent
that 
$$ P_Y (K) = K \cap Y \ , \ \ K = N + K \cap Y \ , 
\eqno{(1.7.5.2)} $$

\noindent
and cone $K \cap Y \ $ is a normal one in the normed space
$\left(Y , p \right) \ $. Thus, by M. Krein's lemma (see[37])
$$ Y^* = (Y\cap K)^* - (Y\cap K)^* \ . \eqno{(1.7.5.3)}$$

\noindent
Now for any $L \in X^* \ $ equality $\left\| L \right\|_{X^*}  
:= \sup \{ \left| L(x) \right| 
\  \left| \ x \in X , \, p(x) \leq 1  \right. \} \ $ 
implies $L(x) = 0 \ $ $\forall \ x \in N \ $ and so defined 
by formula $l(y) := L(y) \ $ $\forall \ y \in Y \ $ 
functional $l \ $ will be an element of $Y^* \ $. According 
to (1.7.5.3)
$l = l_1 - l_2 \ $, $l_i \in Y^* \ $, $l_i (Y \cap K) 
\geq 0 \ $, $i \in \left\{1, 2 \right\} \ $. Extending 
each functional $l_i \ $ onto the whole space $X \ $ by 
formula
$L_i (x) := l_i (P_Y x ) \ $ with the help of (1.7.5.1) we 
get
$L_i \in X^* \ $ and $0 \leq L_i (Y \cap K) = 
L_i (N + Y \cap K) \stackrel{(1.7.5.2)}{=} L_i (K ) \ $, 
$i \in \left\{1, 2\right\} \ $. That is why 
$L = L_1 - L_2 $, $L_1 , L_2 \in K^* \ $, as was to 
be proved.
}}

 \bigskip
\bigskip
{\textsc{ 1.7.6. Proof of Theorem 1.3. \ }}
{\normalsize{ Since $C^0_w \equiv C^0_{M_w } \ $ then it 
is sufficient to prove the statement of theorem only in the 
case when function $w \ $ is upper semicontinuous 
on ${\mathbb{R}} \ $.

Let $L \in (C^0_w )^* \ $ and $K \ $ be a cone of all 
nonnegative on the real line functions from $C^0_w \ $, which
is a normal one in the seminormed space $C^0_w \ $. Using 
Lemma 1.1 we can find such $L_+ , L_- \in (C^0_w )^* \ $, 
that $$ L = L_+ - L_- \ ,\ L_+ (K) \geq 0 , \ L_- (K) 
\geq 0 \ . \eqno{(1.7.6.1)}$$

\noindent
Formula (1.2.2) allows us to find measures $\mu_L , \mu_L^+ 
\mu_L^-  \in {\mathcal{M}} ({\mathbb{R}} )  \ $, relevant 
to the functionals $L, L_+ , L_- \ $. Exploiting regularity 
of the measures in ${\mathcal{M}} ({\mathbb{R}} )  \ $ and 
density of all compactly supported continuous functions in 
the spaces $L_1 (\nu ) \ $, $\nu \in \left\{\mu_L^+ , 
\mu_L^- , \mu_L^+ + \mu_L^- \ \right\} \ $ (see 1.1), it is 
easy to verify that $\mu_L = \mu_L^+  - \mu_L^- \ $ and 
measures
$ \mu_L^+ , \mu_L^- \ $ are positive.

Consider now an arbitrary measure $\nu \in  \left\{\mu_L^+ , 
\mu_L^-  \right\} \ $ and  corresponding functional 
$L_{\nu }  \in  \left\{L^+ , L^-  \right\} \ $, that for 
any
$  f \in C^0 ({\mathbb{R}} )$:
$$L_{\nu } (f) = \int\limits_{{\mathbb{R}} }^{} f(x)  
d \nu(x)   ; \ \ \ \ \ 
\left\| L_{\nu } \right\|  := \sup \{ \left| L_{\nu }(f) 
\right| \  \left|  f \in C^0_w , \, ||f||_w  \leq 1  \right. \} 
< \infty  . \eqno{(1.7.6.2)}$$

\noindent
Since for arbitrary $\varepsilon > 0 \ $ function $1 / 
(\varepsilon + w(x)) \ $ is lower semicontinuous then using 
the known fact from [17, I, Th.1.4] we get a nondecreasing 
sequence of positive and continuous on the whole real axis 
functions
$\varphi_n^{\varepsilon } (x) $, $n \geq 1 \ $: $\lim_{n \to
\infty} \varphi_n^{\varepsilon } (x) = 1 / (\varepsilon + 
w(x)) \ $ $\forall x \in {\mathbb{R}} \ $. Setting
$$ w_n^{\varepsilon } (x) := e^{- \frac{x^2 }{n } } 
\varphi_n^{\varepsilon } (x) \ , \ n \geq 1 , \ x \in 
{\mathbb{R}} \ , \eqno{(1.7.6.3)} $$

\noindent
and taking into account $\left\|\varphi_n^{\varepsilon } 
\right\|_{C({\mathbb{R}})} \leq 1/\varepsilon \ $, we 
obtain
$w_n^{\varepsilon } \in C^0 ({\mathbb{R}}) \ $, 
$\left\|w_n^{\varepsilon } \right\|_w \leq 1 \ $,
$$ 0 < w_n^{\varepsilon } (x) \leq w_{n + 1}^{\varepsilon }  
(x) \leq \frac{1}{\varepsilon + w(x) } \ \forall \ n \geq 1 
\ ; \ \lim\limits_{n \to \infty} w_n^{\varepsilon } (x) = 
\frac{1}{\varepsilon + w(x) } \ \forall x \in {\mathbb{R}} 
\  .   \eqno{(1.7.6.4)}$$

\noindent
By (1.7.6.2) $\int_{{\mathbb{R}}}^{} w_n^{\varepsilon } (x)
\ d \nu (x)  \leq \left\|L_{\nu } \right\| \ $$\forall n 
\geq 1 \ $$\forall \varepsilon > 0 \ $, and according to 
Beppo-Levi theorem $\nu ({\mathbb{R}} \setminus S_w ) = 0 
\ $, $1/ w \in L_1 (\nu ) \ $, $\left\|1/w \right\|_{L_1 
(\nu )} \leq \left\|L_{\nu }\right\| \ $. That is why measure 
$$\rho (A ) := \int\limits_{A}^{} \frac{1}{w(x)} \ 
d \nu (x) \ \ \ \forall A \in {\mathcal{B}} ({\mathbb{R}}) 
\ \eqno{(1.7.6.5)}$$

\noindent
will be positive measure in ${\mathcal{M}} ({\mathbb{R}}) \ $,
$\rho ({\mathbb{R}} \setminus S_w ) = 0 \ $ and $\left\|\rho 
\right\| \leq \left\|L_{\nu } \right\| \ $. An evident 
inequality $\nu (A) \leq \rho (A) \ $$\forall A \in 
{\mathcal{B}} ({\mathbb{R}}) \ $ due to Radon-Nikodym 
theorem means that there exists such $\alpha \in L_1 (\rho ) 
\ $ that  $$ \nu (A ) = \int\limits_{A}^{} \alpha (x) \ d 
\rho (x) \ \ \ 
\forall \ A \in {\mathcal{B}} ({\mathbb{R}}) \ , 
\eqno{(1.7.6.6)}$$

\noindent
and also $0 \leq \alpha (x) \leq 1 \ $ almost everywhere with 
respect to measure $\rho \ $. Using  changes of variables 
theorem [6, V.3], (1.7.6.5), (1.7.6.6) we get $\nu (A) = 
\int_{A}^{} \frac{\alpha (x)}{w(x)} d \nu (x) \ $ $\forall 
\ A \in {\mathcal{B}} ({\mathbb{R}}) \ $, from where $\alpha 
(x) = w(x) \ $ almost everywhere with respect to measure 
$\nu \ $, and by mutual absolute continuity of the measures 
$\nu \ $ and $\rho \ $: $\alpha (x) = w(x) \ $ almost 
everywhere with respect to measure $\rho  \ $. Therefore 
$\nu (A) = 
\int_{A}^{} w(x) d \rho  (x) \ $ $\forall \ A \in 
{\mathcal{B}} ({\mathbb{R}}) \ $ and due to (1.7.6.2):
$ L_{\nu } (f) = \int\limits_{{\mathbb{R}} }^{} w(x)f(x) \ 
d \rho (x) \ $ $\forall \ f \in C^0 ({\mathbb{R}} )  \ $. 
That equality in view of density $C^0 ({\mathbb{R}}) \ $ in 
the seminormed space $C^0_w \ $ and according to the 
continuity of both its sides can be extended to the whole 
$C^0_w \ $:
$$ L_{\nu } (f) = \int\limits_{{\mathbb{R}} }^{} w(x)f(x) 
\ d \rho (x) \ \ \ \forall \ f \in C^0_w  \ ; \ \ \ 
\rho ({\mathbb{R}} \setminus S_w ) = 0 \ . \eqno{(1.7.6.7)}$$

Denoting constructed measures $\rho \ $ by $\mu^+ \ $ and
$\mu^- \ $ when $\nu \ $ equals to $\mu_L^+ \ $ and $\mu^-_L 
\ $, respectively, and setted $\mu := \mu^+  - \mu^- \ $,
we will get the required representation (1.6.1) taking into 
account  $|\mu | ({\mathbb{R}} \setminus S_w ) = 0 \ $.

Since the inverse statement of the theorem and inequality
$\left\|L\right\| \leq |\mu | (S_w ) \ $ are evident to 
finish the proof one need to show only that $\left\|L\right\| 
\geq |\mu | (S_w ) \ $.

Taking a Hahn expansion  ${\mathbb{R}} = {\mathbb{R}}^+_{\mu} 
\sqcup {\mathbb{R}}^-_{\mu } \ $ with respect to measure 
$\mu \ $  (see 1.1)  and any $R > 0 \ $ we rename 
introduced in (1.1.1) functions by:
$$ \kappa^+_{n, R } := \psi_n [{ I_R \cap {\mathbb{R}}^+_{\mu}
, \  \mu }] \ ; \ \ \kappa^-_{n, R } := \psi_n [{I_R \cap 
{\mathbb{R}}^-_{\mu}  , \   \mu }] \ , \ n \geq 1 \ . $$

\noindent
Then in view of (1.7.6.3),(1.7.6.4): 
$w_m^{\varepsilon } \cdot (\kappa^+_{n, R } - \kappa^-_{n, R }
) \in C^0 ({\mathbb{R}}) \ $, 
$\left\|w_m^{\varepsilon } \cdot (\kappa^+_{n, R } - 
\kappa^-_{n, R } ) \right\|_w \leq 1 \  $ $\forall \ n,m 
\geq 1 \ $, $R, \varepsilon > 0 \ $, and by definition of 
the norm
(see (1.7.6.2) and 1.1):
$$ \left\|L\right\| \geq \int\limits_{{\mathbb{R}}}^{} w(x)
w_m^{\varepsilon } (x) (\kappa^+_{n, R } (x) - 
\kappa^-_{n, R }  (x) ) \ d \mu (x) \ \ \forall \ n,m 
\geq 1 \ , R, \varepsilon > 0 \ . \eqno{(1.7.6.8)}$$

\noindent
Passages to the limit in (1.7.6.8) as $n \to \infty \ $ with
regard to (1.1.1) and then as $m \to \infty \ $ and 
$\varepsilon \downarrow 0 \ $ using Beppo-Levi theorem,
give us (see 1.1): $\left\|L \right\| \geq |\mu| (I_R 
\cap S_w ) \ $ $\forall \ R > 0 \ $, i.e.  $\left\|L 
\right\| \geq |\mu| (S_w ) \ $. Theorem 1.3 is proved. }}

\bigskip
\bigskip
\bigskip
\bigskip

\begin{center}
{\bf{ \large{CHAPTER II. }}}  {{ 
\large{Hamburger criterion of the polynomial density
in $C^0_w \ $  and $ L_p (\mu ) \ $, $1 \leq p < \infty \ $
}}}
\end{center}

\bigskip
\bigskip
{\bf 2.1. Notations. \ } 
{\small{ Let  $C^* ({\mathbb{R}} ) \ $ denote the collection 
of all nonnegative upper semicontinuous 
on the whole real line functions $w \ $ satisfying condition 
$ \left\| x^n w \right\|_{C({\mathbb{R}} )} < + \infty \ \ \ 
\forall \ n \in {\mathbb{Z}}_0 \ $ and ${\mathcal{M}}^* 
({\mathbb{R}}) \ $ -- the set of all positive measures 
$\mu \in {\mathcal{M}}({\mathbb{R}}) \ $
which have all finite moments $\int_{{\mathbb{R}}}^{} |x|^n 
\ d \mu (x) < \infty \ $ $\forall n \in {\mathbb{Z}}_0 \ $ 
and unbounded support ${\mathrm{supp}} \mu := \left\{x \in 
{\mathbb{R}} \ | \  \mu (x - \delta , x + \delta ) > 0 \ 
\right.$ \\ $\left.  \forall \ \delta > 0 \ \right\}\ $.

In order to abridge notations in this chapter introduce
${\mathbb{R}}^* := \left\{*\right\} \cup [1, + \infty ) \ $
and for $\mu \in C^* ({\mathbb{R}}) \ $ rename introduced in
(1.2.1) $C^0_{\mu } = \left(C^0_{\mu } , \left\|\cdot 
\right\|_{\mu }\right)\ $ by $L_* (\mu) := \left(L_* (\mu) 
, \left\|\cdot \right\|_{L_* (\mu) }\right)$. That is why
consideration $L_{\alpha } (\mu ) \ $ for $1 \leq \alpha 
< \infty \ $ will mean that $\mu \in {\mathcal{M}}^* 
({\mathbb{R}}) \ $, but for $\alpha = * \ $ it will 
signify under our stipulation that $\mu \in C^* 
({\mathbb{R}}) \ $.
For every $\alpha \in {\mathbb{R}}^* \ $ complex spaces in 
contrast to the real ones $L_{\alpha } (\mu ) \ $ will be 
denoted by $L^c_{\alpha } (\mu ) \ $. As well as in Chapter I:
\ $ \ S_{\mu } = \left\{x \in {\mathbb{R}} | \ \mu (x) > 0 \ 
\right\} \ $ for $\mu \in C^* ({\mathbb{R}}) \ $.

Denote for $\alpha \in {\mathbb{R}}^* \ $ and $z \in 
{\mathbb{C}} \ $:
$$ M_n^{\alpha } (\mu , z) := \sup 
\left\{|p(z)| \ \left| \ \left\|p\right\|_{L_{\alpha } 
(\mu ) }  \ \leq 1 , \ p \in {\mathcal{P}}_n [{\mathbb{C}}] \ 
\right. \right\} \ 
, \ n \in \ {\mathbb{Z}}_0 \ , \eqno{(2.1.1)} $$
$$ {\rho }_n^{\alpha } (\mu , z ) := \inf
\left\{\left\| p \right\|_{L_{\alpha } (\mu ) }\ \Bigl| \ |p(z)| 
\ = 1 , \ p \in {\mathcal{P}}_n [{\mathbb{R}}] \ \Bigr.\right\} \ 
, \ n \in \ {\mathbb{Z}}_0 \ , \eqno{(2.1.2)} $$
$$ M_{\alpha } (\mu , z) := \lim\limits_{n \to \infty }
 M_n^{\alpha } (\mu , z) \ ; \ \ \ 
 {\rho }_{\alpha } (\mu , z ) :=
 \lim\limits_{n \to \infty } {\rho }_n^{\alpha } (\mu , z ) \ 
 . \eqno{(2.1.3)} $$

\noindent
It easy to verify that
$$ \frac{1}{{\rho }_n^{\alpha } (\mu , z ) } = \sup 
\left\{ |p(z)|  \ \Bigl| \ \left\| p \right\|_{L_{\alpha } 
(\mu ) }   \leq 1 ,  p \in {\mathcal{P}}_n [{\mathbb{R}}]   
\right\}  
,  \frac{1}{{\rho }_n^{\alpha } (\mu , z ) } \leq  
M_n^{\alpha } (\mu , z) \leq \frac{2}{{\rho }_n^{\alpha } 
(\mu , z ) }  . \eqno{(2.1.4)}$$

\noindent
Introduce 
$$ \left\{ \begin{array}{lll}
d \mu_{\alpha } (x) := \frac{1}{(1 + |x|)^{\alpha } } d 
\mu (x) \ ;  &  d \mu^{(2)}_{\alpha } (x) := 
\frac{|x|^{\alpha }}{(1 + |x|)^{\alpha } } d \mu (x) \ ; &
1 \leq \alpha < \infty , \ \mu \in {\mathcal{M}}^* 
({\mathbb{R}}) \ ; \\
\mu_* (x) :=  \frac{1}{1 + |x| } \mu (x) \ ;   &
\mu^{(2)}_{* } (x) := \frac{|x|}{1 + |x| }  \mu (x) \ ;& 
 \mu \in C^* ({\mathbb{R}}) \ . \end{array}\right.
 \eqno{(2.1.5)} $$

\noindent
Restricting the polynomial class in (2.1.1) and (2.1.4) to 
the vanishing at zero polynomials we get for $z \in 
{\mathbb{C}} \ $,
$n \geq 1 \ $ and $\alpha \in {\mathbb{R}}^* \ $:
$$ \left\{\begin{array}{ll}
 M_n^{\alpha } (\mu_{\alpha } , z) \ \geq \  |z| \
 M_{n-1}^{\alpha } (\mu^{(2)}_{\alpha } , z) \ ; & 
 {\rho }_{n-1}^{\alpha } (\mu^{(2)}_{\alpha } , z ) \ 
 \geq \  |z| \ {\rho }_{n}^{\alpha } (\mu_{\alpha } , z ) \ ; 
 \\
M_{\alpha } (\mu_{\alpha } , z) \ \geq \  |z| \ M_{\alpha } 
(\mu^{(2)}_{\alpha } , z) \ ; & {\rho }_{\alpha } 
(\mu^{(2)}_{\alpha } , z ) \ \geq \ |z| \ {\rho }_{\alpha } 
(\mu_{\alpha } , z ) \ . 
  \end{array}\right. \eqno{(2.1.6)} $$
}}

\bigskip
\bigskip
{\bf 2.2. Background. \ } 
Functions $\rho_n (z) := \frac{1}{M_n^2 (\mu , z ) } \ $,
$n \in {\mathbb{Z}}_0 \ $, were introduced by \\ H. Hamburger [15]
in connection with the investigation of an indeterminate 
moment problem. These functions were used by M. Riesz [35] to 
obtain the criterion of the polynomial density in 
$L_2 (\mu ) \ $.
For $\alpha = * \ $ and discrete set $S_{\mu} \ $ function 
$M_* (\mu , z) \ $ was introduced by T. Holl [19] and for an
arbitrary $\mu \in C^* ({\mathbb{R}}) \ $ - by S. Mergelyan
in [32]. B.Ja.Levin [30] generalized these results in the 
following statement a simpler proof of which was found 
recently by Ch. Berg [8].

\bigskip
\bigskip
{\textsc{Proposition 2.1.([30]) }} \ {\textit{Let $\alpha \in 
{\mathbb{R}}^* \ $. If ${\mathcal{P}} [{\mathbb{C}}] \ $ is 
dense in $L^c_{\alpha } (\mu ) \ $ then $M_{\alpha } 
(\mu_{\alpha } ; z ) = \infty \ $ $\forall \ z \in 
{\mathbb{C}} \setminus {\mathrm{supp}} \mu \ $. 
If $\exists z \in {\mathbb{C}} \setminus {\mathrm{supp}} 
\mu \ $: 
$M_{\alpha } (\mu_{\alpha } ; z ) = \infty \ $, then 
${\mathcal{P}} [{\mathbb{C}}] \ $ is dense in $L^c_{\alpha } 
(\mu ) \ $. }}

\bigskip
Denote by ${\mathrm{Close}}_{L_{\alpha } (\mu )} A \ $ the 
closure of $A \subseteq L_{\alpha } (\mu ) \ $ in the space
$L_{\alpha } (\mu ) \ $. It is known and it can be easily 
seen from the Proposition 2.1 and (2.1.4) that ${\mathcal{P}}
[{\mathbb{C}}] \ $ is dense in $L^c_{\alpha } (\mu ) \ $ 
if and only if ${\mathcal{P}} [{\mathbb{R}}] \ $ is dense 
in $L_{\alpha } (\mu ) \ $. That is why everywhere below we 
will examine only real case and use the following statement.

\bigskip
\bigskip
{\textsc{Proposition 2.2. }}([30;  2, Th.2.3.2]) {\textit{ 
\ Let $\alpha \in {\mathbb{R}}^* \ $. The following 
statements are equivalent: }}
$$\begin{array}{ll}
(2.2.1a)\  {\mathrm{Close}}_{L_{\alpha } (\mu )} {\mathcal{P}} 
= L_{\alpha } (\mu ); & (2.2.1d) \  \exists z \in 
{\mathbb{C}} \setminus {\mathrm{supp}} \mu \ : 
{\rho }_{\alpha } (\mu_{\alpha } , z ) = 0 ;  \\
(2.2.1b) \ {\mathrm{Close}}_{L^c_{\alpha } (\mu )} 
{\mathcal{P}}[{\mathbb{C}}] = L^c_{\alpha } (\mu );  
& (2.2.1e) \  {\rho }_{\alpha } (\mu_{\alpha } , z ) = 0 \ 
\ \forall z \in {\mathbb{C}} \setminus {\mathrm{supp}} \mu \ ;
\\    (2.2.1c) \ 
{\small{ \frac{1}{x + i } 
}}
\in {\mathrm{Close}}_{L^c_{\alpha } (\mu )} {\mathcal{P}}
[{\mathbb{C}}] \ ; &  \ 
(2.2.1g) \ 
{\small{
 \frac{1}{1+ x^2 }, \frac{x}{1 + x^2 }  
}} \in
{\mathrm{Close}}_{L_{\alpha } (\mu )} {\mathcal{P}} \ . 
\end{array} $$

\bigskip
H.Hamburger in [15] established another criterion of the 
indeterminacy of a moment problem a simpler proof of which 
was given by M. Riesz [36]. This criterion can be formulated 
as follows:
$$ {\mathrm{Close}}_{L_{2 } (\mu )} {\mathcal{P}} \neq 
L_{2 } (\mu ) \ \Leftrightarrow \ 
\rho_2 (\mu_2 , 0 ) > 0 \ \mbox{and} \ \rho_2 (\mu^{(2)}_2 , 
0 ) > 0 \ . \eqno{(2.2.2)} $$

\noindent
Succeeding Berg's proof [8] of the Proposition 2.1 and using 
Theorem 1.3 we will extend here criterion (2.2.2) to all 
spaces $L_{\alpha } (\mu ) \ $, $\alpha \in {\mathbb{R}}^* \ $,
designated
$$ \rho_n^{\alpha } (\mu ) := \rho_n^{\alpha } (\mu , 0 ) \ , 
\ \ n \in {\mathbb{Z}}_0 \ ; \ \ \ 
\ \rho_{\alpha } (\mu ) := \rho_{\alpha } (\mu , 0 ) \ ,
\ \ \alpha \in {\mathbb{R}}^* \ . \eqno{(2.2.3)}$$

\bigskip
\bigskip
{\bf 2.3.  Main Theorem. \ } 

Arbitrary change of zeros of some polynomial $p \in 
{\mathcal{P}}[{\mathbb{C}}] \ $ by the complex conjugate ones
gives the polynomial set $\pi (p) \ $ containing only one
polynomial $p^* \in \pi (p) \ $ all zeros of which lie in 
the lower complex halfplane ${\mathbb{C}}^- := \left\{z \in 
{\mathbb{C }} \ | \ {\mathrm{Im}} z \leq 0 \right\} \ $. 
It is evident that $|q(x)| = |p(x)| \  $ $\forall x \in 
{\mathbb{R}} \ $$\forall \ q \in \pi (p) \ $, and therefore
for any $\alpha \in {\mathbb{R}}^* \ $: $\left\|q\right\|_{
L_{\alpha } (\mu)} = \left\| p
\right\|_{L_{\alpha } (\mu)} \ $ $\forall \ q \in \pi (p) \ $. 
Besides that for arbitrary $a \in {\mathbb{R}} \ $
and $y \geq 0 \ $: $|p^* (a + i y ) | \geq |q (a + i y )
| \ $ $\forall \ q \in \pi (p) \ $, and $|p^* (a + i y )
| \ $ is a nondecreasing function of $y \geq 0 \ $. 
That is why for any
$\alpha \in {\mathbb{R}}^* \ $, $n \in {\mathbb{Z}}_0 \ $,
$a \in {\mathbb{R}} \ $ and $y \geq 0 \ $:
$$ M_n^{\alpha } (\mu , a + i y ) = 
\sup \left\{|p^* (a + i y )| \ \left| \ \left\|p
\right\|_{L^c_{\alpha } (\mu) } \leq 1 \ , \ p \in 
{\mathcal{P}}_n [{\mathbb{C}}] \right. \right\} \ ,  
\eqno{(2.3.1)} $$

\noindent
and $M_n^{\alpha } (\mu , a + i y ) \ $ is a nondecreasing 
function of $y \geq 0 \ $. Thus, an obvious equality 
$M_n^{\alpha } (\mu , z ) = M_n^{\alpha } (\mu , 
\overline{z} ) \ $$\forall z \in {\mathbb{C}} \ $ 
implies validity of the following statement.

\bigskip
\bigskip
{\textsc{Proposition 2.3. }} \ {\textit{For arbitrary 
$\alpha \in {\mathbb{R}}^* \ $, $a \in {\mathbb{R}} \ $ 
and $n \in {\mathbb{Z}}_0 \ $ functions $M_n^{\alpha } 
(\mu , a + i y )
\ $, $M_{\alpha } (\mu , a + i y ) \ $ of the variable 
$y \in {\mathbb{R}} \ $ are even on 
${\mathbb{R}} \ $ and   nondecreasing on $[0, + \infty ) \ $ 
 and, in particular,
$$ \frac{1}{{\rho }_{\alpha } (\mu , x ) } \leq  M_{\alpha } 
(\mu , x ) \leq M_{\alpha } (\mu , x + i y ) \ \ \
\forall \ \  x, y \in {\mathbb{R}} , \alpha \in {\mathbb{R}}^*
\ , \eqno{(2.3.2)} $$ }}

\noindent
{\textit{where $\frac{1}{0} := + \infty \ $ and $+ \infty 
\leq + \infty \ $. }}

\bigskip
The following criterion of the polynomial density in 
$L_{\alpha } (\mu ) \ $, $\alpha \in {\mathbb{R}}^* \ $, 
is the main result of that Chapter.

\bigskip
\bigskip
{\textsc{Theorem 2.1. }} ({\emph{Hamburger local criterion}}) \ 
{\textit{Let  $\alpha \in {\mathbb{R}}^* \ $. Algebraic 
polynomials  ${\mathcal{P}} \ $ are not dense in 
$L_{\alpha } (\mu ) \ $ if and only if }}
$$ \rho_{\alpha }  (\mu_{\alpha } ) > 0  \ \  \ \mbox{and} \ 
\ \  \rho_{\alpha } (\mu^{(2)}_{\alpha } ) > 0 \ , 
\eqno{(2.3.3)}$$

\noindent
{\textit{where $\mu_{\alpha } \ $, $\mu^{(2)}_{\alpha } \ $ 
are defined in (2.1.5) and (see 2.1) }}
$$ \rho_{\alpha }  ( \nu) = \inf \left\{\left\|p
\right\|_{L_{\alpha } (\nu ) } \ \left| \ p(0) = 1 , 
\ p \in {\mathcal{P}}  \right. \right\} \ ,	\ \ \ 
\nu \in \left\{ \mu_{\alpha } , \  \mu^{(2)}_{\alpha } \ 
\right\} \ , \eqno{(2.3.4)} $$

\bigskip
\bigskip
{\bf 2.4.  Proof of Theorem 2.1. \ } 

\bigskip
{\textsc{2.4.1.}}{\textit{ Sufficiency.}}
{\normalsize{ \ Let $\alpha \in {\mathbb{R}}^* \ $ and 
(2.3.3) is valid.
Then $\mu^{(2)}_{\alpha } \not\equiv 0 \ $ and consequently,
$S_{\mu} \setminus \left\{0 \right\} \neq \emptyset \ $, if 
$\alpha = * \ $, and ${\mathrm{supp}} \mu \setminus 
\left\{0\right\} \neq \emptyset \ $, if $\alpha 
\in [1, + \infty ) \ $. Due to (2.1.4)
$$ |p(0)| \leq \frac{\left\|p\right\|_{L_{\alpha } 
(\mu_{\alpha } ) } }{\rho_{\alpha } (\mu_{\alpha } ) } \ ;
\ \ \ \ |p(0)| \leq \frac{\left\|p\right\|_{L_{\alpha } 
(\mu^{(2)}_{\alpha } ) } }{\rho_{\alpha } (\mu^{(2)}_{\alpha } 
) } \ \ \ \    \forall \ p \in {\mathcal{P}} \ .
\eqno{(2.4.1)} $$

\noindent
By Hahn-Banach theorem, (2.4.1) and Theorem 1.3 there exist
such $f_{\alpha } \in L_{\beta } (\mu_{\alpha }) \setminus 
\left\{0 \right\} \ $, $g_{\alpha } \in L_{\beta } 
(\mu^{(2)}_{\alpha }) \setminus  \left\{0 \right\} \ $, if 
$\alpha \in [1, + \infty ) \ $ where $\alpha \ $, $\beta \ $ 
are dual exponents, and  such $\kappa , \gamma  \in
{\mathcal{M}} ({\mathbb{R}}) \ $, $|\kappa | (S_{\mu }) > 0 
\ $, $|\kappa | ({\mathbb{R}} \setminus S_{\mu} ) = 0 \ $,
$|\gamma | (S_{\mu} \setminus \left\{0\right\} ) > 0 \ $,
$|\gamma | ({\mathbb{R}} \setminus ( \ S_{\mu} \setminus 
\left\{0\right\} ) \ ) = 0 \ $, if $\alpha = * \ $, that
respectively to the considered cases:
$$ p(0) = \int_{{\mathbb{R}}}^{} p(t) f_{\alpha } (t)
\ d \mu_{\alpha } (t) \ = \ 
\int_{{\mathbb{R}}}^{} p(t) g_{\alpha } (t)
\ d \mu^{(2)}_{\alpha } (t) \ \ \  \forall \ p \in 
{\mathcal{P}} \ ,  \eqno{(2.4.2)} $$
$$ p(0) = \int_{{\mathbb{R}}}^{} \mu_* (t) p (t)
\ d \kappa  (t) \ = \ 
\int_{{\mathbb{R}}}^{}\mu^{(2)}_* (t) p(t) 
\ d \gamma (t) \ \ \  \forall \ p \in {\mathcal{P}} \ .  
\eqno{(2.4.3)} $$

\noindent
Consider at first the case $\alpha = * \ $. According to 
the notations (2.1.5) equality (2.4.3) can be rewritten in 
the following way:
$$ p(0) = \int_{{\mathbb{R}}}^{} \mu (t)p(t) \ d \kappa_1 (t)
\ = \ \int_{{\mathbb{R}}}^{} \mu (t) p(t) \ d \gamma^{(2)}_1
(t) \ \ \ \ \forall \ p \in {\mathcal{P}} \ , \eqno{(2.4.4)}$$

\noindent
where, obviously, $\kappa_1 , \gamma^{(2)}_1 \in {\mathcal{M}}
({\mathbb{R}})  \ $. Applying left equality (2.4.4) to the 
polynomials vanishing at zero we get:
$$ 0 = \int_{{\mathbb{R}}}^{} \mu (t)p(t) \ d \ 
\widetilde{\kappa } (t) \ \ \ \forall \ p \in {\mathcal{P}} 
\ , \ \ 
d \ \widetilde{\kappa } (t)  :=  \frac{t}{1 + |t| } 
d \kappa (t)
 \ . $$
 
\noindent
If $|\kappa | (S_{\mu} \setminus \left\{0\right\} ) > 0 \ $ 
then by Theorem 1.3 $\widetilde{\kappa } \in L_* (\mu )^* 
\setminus \left\{0\right\} \ $ and therefore $
{\mathrm{Close }}_{L_* (\mu ) } {\mathcal{P}} \neq L_* (\mu )
\ $. If  $|\kappa | (S_{\mu} \setminus \left\{0\right\} ) 
= 0 \ $, then $0 \in S_{\mu } \ $ and  $|\kappa | 
({\mathbb{R}} \setminus \left\{0\right\} ) = 0 \ $ but by 
(2.4.3)
$d \kappa = (1/ \mu (0) ) \delta_0  \ $, where $ \delta_0  \ $
- Dirac's measure [13, 4.4.1] at the point $0 \ $. That is why 
it follows from (2.4.4) that
$$ 0 = \int_{{\mathbb{R}}}^{} \mu (t)p(t) \ d ( \kappa_1 (t) - 
\gamma_1^{(2)} (t) ) \ \ \ \forall \ p \in {\mathcal{P}} \ ,
$$

\noindent
and in addition $ ( \kappa_1 - \gamma_1^{(2)} ) 
(\left\{0\right\}) = 1/ \mu (0) \ $. This means due to 
Theorem 1.3 and $0 \in S_{\mu } \ $ that $\kappa_1 - 
\gamma_1^{(2)} \in L_* (\mu )^* \setminus \left\{0\right\} 
\  $ and hence, ${\mathrm{Close}}_{L_* (\mu ) } 
{\mathcal{P}} \neq  L_* (\mu ) \ $.

Let now $\alpha \in [1, + \infty ) \ $. Applying left 
equality (2.4.2) to the vanishing at zero polynomials we get
$$ \int_{{\mathbb{R}}}^{} p(t) \frac{t f_{\alpha } (t)}{
(1 + |t|)^{\alpha } } \ d \mu (t) \ = \ 0 \ \ \forall \ p 
\in {\mathcal{P}} \ . $$

\noindent
If (see 1.1)  $\mu (S_{|f_{\alpha } | } \setminus 
\left\{0\right\} ) > 0 \  $ then $\frac{t f_{\alpha } (t)}{
(1 + |t|)^{\alpha } } \in L_{\beta } (\mu ) \setminus 
\left\{0\right\} \ $ and consequently, ${\mathrm{Close}}_{
L_{\alpha } (\mu ) } {\mathcal{P}} \neq  L_{\alpha } (\mu ) 
\ $. But if  $\mu (S_{|f_{\alpha } | } \setminus 
\left\{0\right\} ) = 0 \  $ then $\mu (\left\{0\right\}) > 
0 \ $, $0 \in S_{|f_{\alpha } | } \ $
and by (2.4.2): $\mu (\left\{0\right\}) = 1/ f_{\alpha } (0) 
> 0 \  $. Equalities (2.4.2) yield:
$$ 0 = \int_{{\mathbb{R}}}^{} p(t) \varphi_{\alpha } (t)
\ d \mu (t)  \ \ \forall \ p \in {\mathcal{P}} \ , 
\ \ \varphi_{\alpha } (t) = \frac{ f_{\alpha } (t) - 
|t|^{\alpha } g_{\alpha } (t) }{
(1 + |t|)^{\alpha } }  \ . $$

\noindent
It is easy to verify that $\varphi_{\alpha } \in L_{\beta }
(\mu) \ $ and if $d \mu_0 := d \mu - \frac{1}{f_{\alpha } (0) }
\cdot \delta_0 \  $ then for arbitrary $\varepsilon > 0 \ $:
$$ \left\|\varphi_{\alpha } \right\|_{L_{\beta } (\mu ) 
}^{\beta } \ \geq \  \int\limits_{- \varepsilon }^{+ 
\varepsilon } |\varphi_{\alpha } (t)|^{\beta } \ d \mu_0 (t) 
\ + \ \frac{1}{f_{\alpha } (0) } |\varphi_{\alpha } 
(0)|^{\beta } \ \geq \ f_{\alpha } (0)^{\beta - 1}\ > \ 0 
\ , $$

\noindent
if $\alpha > 1 \ $, and $\left\|\varphi_1 
\right\|_{L_{\infty} (\mu ) } \geq f_{\alpha } (0) > 0 \ $,
if $\alpha = 1 \ $, i.e. $\varphi_{\alpha } \in L_{\beta } 
(\mu) \setminus \left\{0\right\} \ $ and hence, 
${\mathrm{Close}}_{L_{\alpha  } (\mu)}  {\mathcal{P}}  
\neq  L_{\alpha  } (\mu) \ $.
}}

\bigskip
{\textsc{2.4.2.}}{\textit{ Necessity.}}
{\normalsize{ \ Let $\alpha \in {\mathbb{R}}^* \ $ and 
${\mathrm{Close}}_{L_{\alpha  } (\mu)}  {\mathcal{P}}  
\neq  L_{\alpha  } (\mu) \ $. Then by Hahn-Banach theorem 
and Theorem 1.3 for $\alpha \in [1, + \infty ) \ $ there 
exists such $ g_{\alpha } \in L_{\beta } (\mu) \setminus 
\left\{0\right\}  \ $, $1/\alpha  + 1/\beta  = 1 \ $, that: 

$$\int\limits_{{\mathbb{R}}}^{} p(t)g_{\alpha } (t) \ d 
\mu (t) \ = \ 0 \ \ \forall \ p \in {\mathcal{P}} \ ,
\eqno{(2.4.5)} $$

\noindent
and if $\alpha = * \ $ then $\exists \gamma \in 
{\mathcal{M}} ({\mathbb{R}}) \ $, $|\gamma | (S_{\mu } ) > 
0 \ $, $|\gamma | ({\mathbb{R}} \setminus S_{\mu } ) = 0
\ $:
$$ \int\limits_{{\mathbb{R}}}^{} \mu (t) p(t) \ d \gamma  
(t) \ = \ 0 \ \ \forall \ p \in {\mathcal{P}} \ .
\eqno{(2.4.6)} $$

\noindent
Under these conditions function

$$ \varphi_{\alpha } (z) := \int\limits_{{\mathbb{R}}}^{} 
\frac{g_{\alpha } (t) }{t - z} \ d \mu (t), \ \mbox{if} \ 
1 \leq \alpha < + \infty \ ; \ 
\varphi_{* } (z) := \int\limits_{{\mathbb{R}}}^{} 
\frac{\mu (t) }{t - z} \ d \gamma  (t), \ \mbox{if} \ 
\alpha = *  \ ; \eqno{(2.4.7)} $$

\noindent
is analytic on ${\mathbb{C}} \setminus {\mathbb{R}} \ $ 
and not identically zero. Thus, $\exists \lambda_{\alpha } 
\in [1,2]  \ $: $\varphi_{\alpha } (i \lambda_{\alpha } ) 
\neq 0 \ $. Besides that it is easy to derive from (2.4.5) 
and (2.4.6) that for any $z \in {\mathbb{C}} \setminus 
{\mathbb{R}} \ $ and $p \in {\mathcal{P}}[{\mathbb{C}}] \ $:
$$ p(z)\varphi_{\alpha } (z) := 
\int\limits_{{\mathbb{R}}}^{} 
\frac{p(t)g_{\alpha } (t) }{t - z} \ d \mu (t),
1 \leq \alpha < + \infty \ ; \ 
p(z) \varphi_{* } (z) := \int\limits_{{\mathbb{R}}}^{} 
\frac{p(t) \mu (t) }{t - z} \ d \gamma  (t),
\alpha = *  \ ; \eqno{(2.4.8)} $$

\noindent
Setted in (2.4.8) $z = i \lambda_{\alpha } \ $ we get for
$\alpha \in [1, + \infty ) \ $:
$$ |p(i \lambda_{\alpha } ) | \leq 
\frac{\left\|g_{\alpha }\right\|_{L_{\beta } (\mu) } }{
|\varphi_{\alpha } (i \lambda_{\alpha } ) |}
\left[\int\limits_{{\mathbb{R}}}^{} \frac{|p(t)|^{\alpha 
}}{|t - i \lambda_{\alpha } |^{\alpha }} \ d \mu (t) 
\right]^{1 / \alpha } \leq \frac{\sqrt{2}\left\|g_{\alpha 
}\right\|_{L_{\beta } (\mu) } }{
|\varphi_{\alpha } (i \lambda_{\alpha } ) |} 
\left\|p\right\|_{L_{\alpha } (\mu_{\alpha })} \ , 
\eqno{(2.4.9)}$$

\noindent
and for $\alpha = * \ $:
$$ |p(i \lambda_{\alpha } ) | \leq 
\frac{\left\|p\right\|_{L_{* } (\mu_* ) } }{
|\varphi_{* } (i \lambda_{* } ) |}
\int\limits_{{\mathbb{R}}}^{}  \frac{1 + |t|}{|t - i   
\lambda_{* } | }  \ d  |\gamma | (t)  \leq \frac{ \sqrt{2} 
\left\|\gamma \right\| }{
|\varphi_{* } (i \lambda_{* } ) |} \left\|p
\right\|_{L_{* } (\mu_{* })} \ . 
\eqno{(2.4.10)}$$

\noindent
Thus, $M_{\alpha } (\mu_{\alpha } , i \lambda_{\alpha } ) 
< \infty \ $ and by (2.3.2): $ 1 / \rho_{\alpha } 
(\mu_{\alpha } ) \leq  M_{\alpha } (\mu_{\alpha } , 
i \lambda_{\alpha } ) < \infty \ $, what together with 
arised from (2.3.2) and (2.1.6) inequality
$$ \frac{1}{\rho_{\alpha } (\mu^{(2)}_{\alpha } ) } \leq 
M_{\alpha } (\mu^{(2)}_{\alpha } , i \lambda_{\alpha } )\leq
 \frac{1}{|i \lambda_{\alpha } | } M_{\alpha } 
 (\mu_{\alpha } , i \lambda_{\alpha } ) \ < \ \infty \ , $$
 
\noindent
gives validity of the inequalities (2.2.2). Theorem 2.1 
is proved. }}

\bigskip
\bigskip
\bigskip
\bigskip
\begin{center}
{\bf{ \large{CHAPTER III. }}}  {{ 
\large{ Hamburger and Krein classes of entire 
functions }}}
\end{center}

\bigskip
\bigskip
{\bf 3.1. Notations and Definitions. \ } 
{\small{ A function $f : {\mathbb{C}} \to {\mathbb{C}} \ $ is
said to be of {\textit{exponential type}} if $|f(z)| \leq C 
e^{\sigma |z|} \ $
$\forall z \in {\mathbb{C}} \ $ for some $\sigma , C > 0 \ $,
and of {\textit{minimal exponential}} type if 
$$ \forall \varepsilon > 0 \ \ \exists \ 
C_{\varepsilon } > 0 \ : \ \ |f(z)| \leq C_{\varepsilon } 
e^{\varepsilon  |z|} 
\ \ \ \forall \ z \in {\mathbb{C}} \ . \eqno{(3.1.1)}$$

\noindent
Let ${\mathcal{E}} , {\mathcal{E}}_1 , {\mathcal{E}}_0 \ $
denote the sets of all entire functions, entire functions of 
exponential type and entire functions of minimal exponential 
type, respectively; $\Lambda_{f} \ $ -- the set of all zeros 
$f \in {\mathcal{E}} \ $; ${\mathrm{Close }}_{{\mathcal{E}}}
A \ $ --  the closure $A \subseteq {\mathcal{E}} \ $ with respect
to topology $\tau_{{\mathcal{E}}} \  $ of the uniform 
convergence on all compact subsets of ${\mathbb{C}} \ $;
${\mathrm{Close }}_{L_{\alpha } (\mu ) } A \ $ -- the closure 
$A \subseteq L_{\alpha } (\mu ) \ $ in the space $L_{\alpha } 
(\mu ) \ $, $\alpha \in {\mathbb{R}}^* \ $ (see 2.1); 
${\mathrm{deg}} P \ $ -- degree of the polynomial $P \in 
{\mathcal{P}}[{\mathbb{C}}] \ $; $a\vee b := \max 
\left\{a, b \right\} \ $, $a, b \in {\mathbb{R}} \ $; 
${\mathrm{co}} A \ $ -- convex hull of $A \subseteq 
{\mathbb{C}} \ $; ${\mathrm{card}} B \in {\mathbb{Z}}_0 
\cup  \left\{\infty \right\} \ $ -- number of elemens in 
the set $B \ $. 
Function $f \in {\mathcal{E}} \ $ is said to be real if 
$f({\mathbb{R}}) \subseteq {\mathbb{R}} \ $. For $n \in 
{\mathbb{Z}}_0 \ $
and $X \in \left\{{\mathcal{P}}_n , {\mathcal{P}} , 
{\mathcal{E}} , {\mathcal{E}}_1 , {\mathcal{E}}_0 \right\} 
\  $ let $X({\mathbb{R}}) \ $ denote the set of real 
functions from $X \ $ with real zeros only and 
$X^* ({\mathbb{R}}) \ $
 -- the set of real functions $f \in X \ $ all zeros of which
 are real, simple and $f(0) = 1 \ $. The sets of real 
 functions from ${\mathcal{P}} , {\mathcal{E}} \ $ with 
 only real and simple zeros will be denoted by 
 ${\mathcal{P}}_s ({\mathbb{R}}), {\mathcal{E}}_s 
 ({\mathbb{R}})\ $, respectively.

 Let us remind that [30, VIII.1] the set $A \subseteq 
 {\mathcal{E}} \ $ is said to be {\textit{normal}} 
 if ${\mathrm{Close }}_{{\mathcal{E}}} A \ $ is a compact 
 set with respect to the topology $\tau_{{\mathcal{E}}} \  $ 
 or, what is the same [26, IV.41, II.20; 37, I.6.1], if any 
 sequence in $A \ $ contains the convergent subsequence with 
 respect to the same topology $\tau_{{\mathcal{E}}} \  $.
 
\bigskip
{\textsc{3.1.1. Cartwright class.}} 
The set of entire functions $f \in {\mathcal{E}}_1 \ $ 
satisfying inequality
$$ \int\limits_{{\mathbb{R}} }^{} \frac{\log^+ 
|f(t)| }{1 + t^2 } \ d t \ < \infty \  , \ \ \log^+ x := 0 
\vee \log x , \ 
x \geq 0 \ , \eqno{(3.1.2)} $$
 
\noindent
is called the Cartwright class and will be denoted by 
${\mathrm{Cartwright}} \ $. Each $f \in {\mathrm{Cartwright}}
\ $ is an element of so-called (A) class of entire functions,
i.e. $\sum_{\lambda \in \Lambda_f \setminus \left\{0\right\}
}^{} \left|{\mathrm{Im}} \frac{1}{\lambda }\right| < \infty 
\ $, and satisfies [29, V.4, Th.7] the stronger inequality:
$$ \int\limits_{{\mathbb{R}} }^{} \frac{\left|\log |f(t)| 
\right| }{1 + t^2 } \ d t \ < \infty \ \ . \eqno{(3.1.3)} $$ 

\noindent
It is known also [29, V.6, Th.13] that $f \in 
{\mathrm{Cartwright}} \cap {\mathcal{E}} ({\mathbb{R}}) \ $
has the following representation:
$$ f(z) = c \cdot z^m \ 
\lim\limits_{R \to + \infty } \prod\limits_{\lambda \in 
(-R, R ) \cap \left( \Lambda_f  \setminus \left\{0\right\} 
\right) }^{}
\left(1 - \frac{z}{\lambda }\right) \ , 
\ z \in {\mathbb{C}} , \ c \in {\mathbb{R}} \setminus 
\left\{0\right\} , \ m \in {\mathbb{Z}}_0 \ , \eqno{(3.1.4)}$$

\noindent 
and there exist the finite limits: $ \lim\limits_{R \to 
+ \infty }
\frac{{\mathrm{card}} \left( \Lambda_f \cap [0,R]\right) }{R} 
\ = \ 
\lim\limits_{R \to + \infty }
\frac{{\mathrm{card}} \left( \Lambda_f \cap [-R, 0] \right) 
}{R} \  $, 
$$ \delta_f := \lim\limits_{ R \to + \infty } \delta_f (R)
 \ ; \ \ \ \delta_f (R) := \sum_{\lambda \in (-R, R ) \cap 
 \left( \Lambda_f  \setminus \left\{0\right\} \right) }^{} 
 \frac{1}{\lambda } \ . \eqno{(3.1.5)} $$ 

\noindent
It is worth to remind that Lindelof and Hadamard's [29, I] 
theorems for any $f \in {\mathcal{E}}_0 ({\mathbb{R}}) \ $
give an existence of the finite limit (3.1.5), equality 
$\delta_f = 0 \ $ and validity of the representation (3.1.4).

\bigskip
{\textsc{3.1.2. Krein class.}}  According to [2, III] 
and [25] (see also [8, 9, 7, 39]) function $f \in 
{\mathcal{E}}_s ({\mathbb{R}}) \ $ is said to be a function 
of Krein class ${\mathcal{K}} \ $ if its reciprocal can be 
represented as a series of 
simple fractions:

$$ \frac{1}{f(z)} \ = \ A + \frac{B}{z} + \sum\limits_{n 
\geq 1 } \ A_n \ 
\frac{z}{\lambda_n (z - \lambda_n ) } \ , 
\ z \in {\mathbb{C}}\setminus \left\{\lambda_n \right\}_{n 
\geq 1 } \ ; $$

\noindent
where $ \lambda_n \neq 0 , A,B, A_n , \lambda_n \in 
{\mathbb{R }} \ 
\forall \ n \geq 1 \ $ and $\sum\limits_{n \geq 1 } 
\frac{|A_n |}{{
\lambda_n }^2  } \ < \ \infty \ $. 

\bigskip
{\textsc{3.1.3. Related to the Krein class definitions.}} 
For every $f \in {\mathcal{E}}_s ({\mathbb{R}}) \ $ define
$$ d_f \ := \ \inf \ \ \left\{q \in {\mathbb{Z}} \ | \ 
\sum\limits_{\lambda \in \Lambda_f \setminus \{ 0 \} }  
\frac{1}{|f^{\prime } (\lambda ) | \cdot |\lambda |^{q + 1}  
} \ \ < \ \ \infty \ \ 
\right\}  \  ,  \eqno{(3.1.6)}  $$

\noindent
considering $\inf \emptyset := + \infty \ $. If  $f \in 
{\mathcal{E}}_s ({\mathbb{R}}) \ $ and $d_f < + \infty \ $,
then for every $p \in {\mathbb{Z}}_0 \ $, $p \geq d_f \ $
one can introduce entire function:
$$ \Delta_f^p (z) :=  \frac{1}{f(z)} \ - \ \frac{\chi_{
\Lambda_f } (0)  }{f^{\prime } (0 ) 
\cdot z } \ - \ 
 \sum\limits_{\lambda \in \Lambda_f \setminus \left\{0\right\}
 }  \frac{z^p }{ \lambda^p f^{\prime } (\lambda  ) 
 (z - \lambda ) } \ ,  z \in {\mathbb{C}} \ , 
 \eqno{(3.1.7)} $$

\noindent
and meromorphic function:

$$ m_f^p (z) :=  \frac{\chi_{\Lambda_f } (0)  }{
f^{\prime } (0 ) \cdot z } \ + \ \sum\limits_{ k =0 }^{
p - 1 } \frac{z^k}{ k ! } \left(\frac{1}{f(z)} -  
\frac{\chi_{\Lambda_f } (0)  }{f^{\prime } (0 ) \cdot z }  
\right)^{(k)} (0) +
 \sum\limits_{\lambda \in \Lambda_f \setminus 
 \left\{0\right\}    }
\frac{z^p }{ \lambda^p f^{\prime } (\lambda  ) (z 
- \lambda ) } \ , z \in {\mathbb{C}} \ ,  \eqno{(3.1.8)} $$

\noindent
where $\sum\limits_{ k =0 }^{- 1 } := 0 \ $.

\bigskip
{\textsc{3.1.4. Hamburger class.}}  Hamburger in [16] defined
the class ${\mathcal{H}} \ $ of entire functions $f \in 
{\mathcal{E}}_s ({\mathbb{R}}) \ $, satisfying the following 
two conditions:

\medskip
$ (3.1.9a) \  \frac{1}{f(z)} \ = \ \sum\limits_{ \lambda \in 
\Lambda_f }^{}  
\frac{1}{f^{\prime }(\lambda ) (z - \lambda ) } \ , \ z 
\in {\mathbb{C}}\setminus \Lambda_f  \ ; $

\medskip
$ (3.1.9b) \ \ \sum\limits_{\lambda \in \Lambda_f } 
\frac{|\lambda |^n }{|f^{
\prime }(\lambda ) | } \ < \ \ \infty \ \ \forall \ n 
\in {\mathbb{Z}}_0  \ . $

\bigskip
{\textsc{3.1.5. Laguerre-P\'olya class.}}
We will consider below a certain subclass of the well-known
second Laguerre-P\'olya class [21, p.336; 18, III.3, Def.3.1;
29, VIII; 27]:
$${\cal LP}_{II}:=\left\{be^{az - c^2 z^2}\prod\limits_{
n\geq 1} (1-{z\over{\lambda_n}}) e^{z\over{\lambda_n}} \ | 
\  a, c \in {\mathbb{R}}, \; b, \lambda_n \in {\mathbb{R}} 
\setminus \left\{0 \right\} 
\; \forall n \geq 1 ,\;
\sum\limits_{n\geq 1}{1\over{\lambda_n^2 }} <\infty\right\},
$$

\noindent
namely, 
$${\cal LP}_{II}^0:=\left\{be^{az}\prod\limits_{n\geq 1} (
1-{z\over{\lambda_n}}) e^{z\over{\lambda_n}} \ | \  a \in {
\mathbb{R}}, \; b, \lambda_n \in {\mathbb{R}} \setminus 
\left\{0 \right\} 
\; \forall n \geq 1 ,\;
\sum\limits_{n\geq 1}{1\over{\lambda_n^2 }} <\infty\right\}.
\eqno{(3.1.10)}$$

}}

\bigskip
\bigskip
{\bf 3.2. Normal polynomial families in $ L_p (\mu )  \ $.\ }

\medskip
{\textsc{3.2.1. }} It follows from the definition of 
$M_{\alpha } (\mu , z ) \ $ that for every $\alpha \in 
{\mathbb{R}}^* \ $ the following implication holds:
$$ p \in {\mathcal{P}}[{\mathbb{C}}] \ , \ \ 
\left\|p \right\|_{L_{\alpha }^c (\mu ) } \leq 1 \ 
\Rightarrow \ |p(z)| \leq M_{\alpha } (\mu , z ) \ 
\forall \ z \in {\mathbb{C }} \ . \eqno{(3.2.1)} $$

\noindent
Inequality (3.2.1) become essential when 
${\mathrm{Close}}_{L_{\alpha } (\mu )} {\mathcal{P}} 
\neq L_{\alpha } (\mu ) \ $. It is known that in 
this case (2.4.8), (2.4.9) and (2.4.10) imply the uniform 
boundedness
of $M_{\alpha } (\mu , z ) \ $ on any segment of the form
$i a + b \cdot I \ $, $a \in {\mathbb{R}} \setminus 
\left\{0\right\} \ $, $b > 0 \ $, which also does not 
include any zero of the defined in (2.4.7) function 
$\varphi_{\alpha } (z) \ $. In view of the Proposition 2.3 
this means the uniform boundedness of $M_{\alpha } 
(\mu , z ) \ $ on any compact subset of the complex plane 
${\mathbb{C}} \ $. Thus, by virtue of Vitali's classical 
compactness theorem and (3.2.1) we deduce that each subset

$$ {\mathcal{P}}_{\alpha } (\mu ) := \left\{p \in 
{\mathcal{P}} \ | \ \left\|p\right\|_{L_{\alpha } (\mu )} 
\leq 1 \  \right\} \ , \ {\mathcal{P}}^c_{\alpha } (\mu ) 
:= \left\{p \in {\mathcal{P}}[{\mathbb{C}}] \ | \ 
\left\|p\right\|_{L^c_{\alpha } (\mu )} \leq 1 \  \right\} 
\ , \eqno{(3.2.2)} $$

\noindent
is normal and so their closures $ {\mathcal{E}}_{\alpha } 
(\mu ) := {\mathrm{Close}}_{{\mathcal{E}}} {\mathcal{P}}_{
\alpha } (\mu) \ $, $ {\mathcal{E}}^c_{\alpha } (\mu ) := 
{\mathrm{Close}}_{{\mathcal{E}}} {\mathcal{P}}^c_{\alpha } 
(\mu) \ $ are compact sets in the topology $\tau_{{
\mathcal{E}} } \  $ and moreover
$$ |f(z)| \leq M_{\alpha } (\mu , z ) \ \forall \ z \in 
{\mathbb{C}} \ \ \forall \ f \in {\mathcal{E}}^c_{\alpha } 
(\mu ) \ . \eqno{(3.2.3)}$$

Inequality (3.2.1) for $\alpha = 2 \ $ was indicated in 
[38, Th.2.5] where also the problem about the complete 
description of the set ${\mathcal{E}}^c_{\alpha } (\mu ) \ $ 
was raised. In addition, known \\ M. Riesz's theorems [2, 
Th.2.4.1,2.4.3] assert that the function $M_{\alpha } 
(\mu , z ) \ $ for $\alpha = 2 \ $ is of minimal exponential 
type and inequality (3.1.2) holds for $f(t) = M_{\alpha } 
(\mu , t ) \ $ there. These two properties of $M_{\alpha } 
(\mu , z ) \ $ were proved for $\alpha = * \ $ by S. Mergelyan
in [32] and for arbitrary $\alpha \in {\mathbb{R}}^* \ $ - 
by B.Ja.Levin in [30].

Observe, that for arbitrary $\alpha \in {\mathbb{R}}^* \ $ 
condition (3.1.2) for $f(t) = M_{\alpha } (\mu , t ) \ $
in view of the evident lower bound:

$$ M_{\alpha } (\mu , z ) \ \geq \ \frac{1}{\left\|\mu 
\right\|_{\alpha } } \ \ \ \forall \ z \in {\mathbb{C}} \ , \ \ 
\left\|\mu\right\|_{\alpha } :=
\left\{\begin{array}{ll}
|\mu| ({\mathbb{R}} ) \ , &  \ 1 \leq \alpha < \infty ; \\
\left\|\mu\right\|_{C ({\mathbb{R}})} \ ,   & \alpha = * ;  
  \end{array}\right. \eqno{(3.2.4)} $$

\noindent
is equivalent to the condition (3.1.3) with the same 
$f(t) \ $. In the next item 3.2.2 we will prove the 
following statement using proofs from [32, item 13] and 
[2, Th.2.4.3].
  
\bigskip
\bigskip
{\textsc{Proposition 3.1. }}([30, 32, 36, 2]) 
{\textit{Let $\alpha \in {\mathbb{R}}^* \ $ and 
${\mathrm{Close}}_{L_{\alpha } (\mu )} {\mathcal{P}} 
\neq L_{\alpha } (\mu ) \ $. Then function
$M_{\alpha } (\mu , \cdot  ) : {\mathbb{C}} \to 
(0, +\infty ) \ $ is of minimal exponential type and 
inequality (3.1.3) holds for $f(t) = M_{\alpha } (\mu , t ) 
\ $.}}

\bigskip
That is why Proposition 3.1 together with inequality (3.2.3) 
shows that for arbitrary  $\alpha \in {\mathbb{R}}^* \ $
 incompleteness of ${\mathcal{P}} \ $ in
$L_{\alpha } (\mu) \ $ means that

$$ {\mathcal{E}}^c_{\alpha } (\mu ) , \ 
{\mathcal{E}}_{\alpha } (\mu ) \subseteq \ {\mathcal{E}}_0 
\cap {\mathrm{Cartwright}} \ \eqno{(3.2.5)}$$

\noindent
and therefore in that case polynomials can approximate in 
the space $L_{\alpha } (\mu) \ $ only those functions which 
from their domain of definition in  $L_{\alpha } (\mu) \ $ 
can be extended into the whole complex plane as an entire 
function of minimal exponential type from the Cartwright 
class.

\bigskip
{\textsc{3.2.2. Proof of Proposition 3.1. }}

{\normalsize{

\medskip
Let $\alpha \in {\mathbb{R}}^* \ $. Since $M_{\alpha } 
(\lambda \cdot \mu , z ) = \frac{1}{\lambda }M_{\alpha } 
(\mu , z ) \ $$\forall \lambda > 0 \ $, then without loss 
of generality one may consider $\left\|\mu\right\|_{\alpha } 
= 1 \ $. By Proposition 2.3 and (3.2.4) to prove inequality 
(3.1.1) it is sufficient to show that:
$$ \forall \varepsilon > 0 \ \exists C_{\varepsilon } > 0 :
\ \log |M_{\alpha } ( \mu , z )| \leq \ C_{\varepsilon } + 
\varepsilon |y| \ \ \forall \ z = x + i y \in {\mathcal{A}}
 \ , \eqno{(3.2.6)}$$
 
\noindent
where ${\mathcal{A}} := \left\{z \in {\mathbb{C}} \ | \ y 
\geq 2 + |x| \ \right\} \ $. But defined in (2.4.7) function
$\varphi_{\alpha } (z) \  $ is uniformly bounded in 
${\mathrm{Im}} z \geq 1 \ $ and so by known corollary of 
Jensen theorem [22, IV.D, VI.C]:
$$ \int\limits_{{\mathbb{R}} }^{} \frac{\left| \ \log 
|\psi (t)| \ \right| }{1 + t^2 } \ d t \ < \infty \ \ ,
\eqno{(3.2.7)} $$  

\noindent
for $\psi (t) = \varphi_{\alpha } (1 + i t) \ $, $t \in 
{\mathbb{R}} \ $. 
Using the similar to (2.4.9) and (2.4.10) estimates one can
easily obtain an existence of such constant 
$C_{\alpha } > 0 \ $ that:
$$ |p(1 + i t ) | \leq \Phi (t) \ \ \ \forall \ t \in 
{\mathbb{R}} \ \ \ \forall \ p \in {\mathcal{P}}_{
\alpha }^{c} (\mu) \ , \eqno{(3.2.8)}$$

\noindent
where:
$$ \Phi (t) := \frac{C_{\alpha } }{|\varphi (1 + i t )|} \ 
 \ \  \forall \  t \in {\mathbb{R}} \ . \eqno{(3.2.9)}$$
 
\noindent
Since (3.2.7) is valid and for $\psi = \Phi  \ $, then
using the Poisson formula [29, V.2, Th.4] for the harmonic 
in ${\mathrm{Im}} z \geq 1 \ $ function $\log p^* (z) \ $ 
(polynomial $p^* \in \pi (p) \subseteq {\mathcal{P}}_{
\alpha }^c (\mu) \ $ has been defined at the beginning of 
the section 2.3.) we will get from the formulas (2.1.4), 
(2.3.1) for 
$z = x + iy \in {\mathcal{A}} \ $ the following inequalities:
$$ \log |M_{\alpha } (\mu , z ) | = 
\sup\limits_{p \in {\mathcal{P}}_{\alpha }^c (\mu) }
\log |p^* (z)| = \sup\limits_{p \in {\mathcal{P}}_{\alpha 
}^c (\mu) } \frac{y - 1 }{\pi } \int\limits_{{\mathbb{R}}}^{}
\frac{\log |p^* (i + t)|}{(t - x )^2 + (y -1 )^2 } \ d t 
\leq $$
$$ \leq \frac{y - 1 }{\pi } \int\limits_{{\mathbb{R}}}^{}
\frac{\gamma (t) }{(t - x )^2 + (y -1 )^2 } \ d t  \ , $$

\noindent
where $\gamma (t) := |\log \Phi (t)| \ $, $t \in 
{\mathbb{R}} \ $. Since for all $z \in {\mathcal{A}} \ $:
$(t - x )^2 + (y -1 )^2 \geq \frac{1}{2} (1 + t^2 ) \ $ 
$\forall \ t \in {\mathbb{R}} \ $, then validity 
for arbitrary $T > 0 \ $ and $z \in {\mathcal{A}} \ $ of 
the following relations (cf. [2, Th.2.4.3, Proof]):
$$ \int\limits_{|t| \geq T }^{} \frac{\gamma (t)}{(t - x )^2 
+ (y -1 )^2 } \ d t \leq \ 2 \int\limits_{|t| \geq T }^{} 
\frac{\gamma (t)}{1 + t^2 } \ d t ; \ 
\lim\limits_{y \to + \infty }  \int\limits_{- T }^{T}
\frac{\gamma (t)}{(t - x )^2 + (y -1 )^2 } \ d t \ = \ 0 
\ , $$

\noindent
implies a correctness of (3.2.6).

To prove (3.1.3) for $f(t) = M_{\alpha } (\mu , t ) \ $, 
observe, at first, that due to (2.1.4):
$$ \frac{1}{4} M_{\alpha } (\mu , x )^2 \leq \sup\limits_{p 
\in {\mathcal{P}}_{\alpha } (\mu) } (1 + p(x)^2 ) = 
\sup\limits_{p \in {\mathcal{P}}_{\alpha } (\mu) } |p_- 
(x)|^2
\ \ \forall \ x \in {\mathbb{R}} \ , $$

\noindent
where the polynomial $p_- \in {\mathcal{P}} [{\mathbb{C}}] 
\ $ contains all zeros of the polynomial $1 + p(z)^2 \ $ 
lying
in ${\mathbb{H}}_- := \left\{z \in {\mathbb{C}} \ | \ 
{\mathrm{Im}}z < 0 \ \right\} \ $ and $1 + p(x)^2 = |p_- 
(x)|^2 \ $ $\forall \ x \in {\mathbb{R}} \ $. But 
$\left\|p_- \right\|_{L_{\alpha }^c (\mu ) } \leq \sqrt{2} 
\left\|1 + |p|\right\|_{L_{\alpha } (\mu ) } \leq 2\sqrt{2} 
\ $ and so
$$ 1 \leq M_{\alpha } (\mu , x ) \leq 2 \sup\limits_{q \in 
{\mathcal{P}}^-_{\alpha } (\mu) } |q(x)| \ \ \ \forall 
\ x \in {\mathbb{R}} \ , \eqno{(3.2.10)}$$

\noindent
where ${\mathcal{P}}^-_{ \alpha } (\mu) := \left\{ q \in 
{\mathcal{P}} [{\mathbb{C}}] \ | \ \Lambda_q  \subset 
{\mathbb{H}}_- , \left\| q  \right\|_{ L_{\alpha }^c (\mu ) }  
\leq 2 \sqrt{2} , |q(x)| \geq 1 \ \forall \ x  \in 
{\mathbb{R}} \ \right\} \ $. Since \\ $ {\mathcal{P}}^-_{\alpha } 
(\mu) \subseteq 2 \sqrt{2} {\mathcal{P}}^c_{\alpha } (\mu) \ $, 
then denoting for the introduced in (3.2.9) function $\Phi \ $:
$\tau (t) := |\log 2\sqrt{2} \Phi (t)| \ $, $t \in 
{\mathbb{R}} \ $, we can state that for arbitrary $q \in 
{\mathcal{P}}^-_{\alpha } (\mu) \ $ function $\log |q(z)| \ $
is a harmonic function in $0 \leq {\mathrm{Im}} z \leq 1 \ $
and in view of (3.2.8) it satisfies inequality: $\log |q(1 + 
i t ) | \leq \tau (t) \ $$\forall \ t \in {\mathbb{R}} \ $. 
But (3.2.7) is true for $\psi = 2\sqrt{2} \Phi \ $ and that is 
why
application of the Poisson formula to $\log |q(z)| \ $ 
gives possibility to continue inequality (3.2.10) as follows:
$$ \log \frac{ M_{\alpha } (\mu , x )}{2} \leq \sup\limits_{
q \in {\mathcal{P}}^-_{\alpha } (\mu) } 
\frac{1}{\pi } \int\limits_{{\mathbb{R}}}^{} \frac{\log |q(1 
+ i t ) | }{1 + (t - x )^2 } \ d t \leq 
\frac{1}{\pi } \int\limits_{{\mathbb{R}}}^{} \frac{\tau (t)
}{1 + (t - x )^2 } \ d t , \ \ \forall \ x \in {\mathbb{R}} 
\ . \eqno{(3.2.11)} $$

\noindent
Using Fubini theorem and equality $\int_{{\mathbb{R}}}^{} 
(1 + x^2 )^{-1 } [1 + (t - x )^2 ]^{-1 } \ d x  = 
2 \pi (t^2 + 4 )^{-1 } \ $, we get from (3.2.11) the required
inequality:
$$ \int\limits_{{\mathbb{R}}}^{} \frac{\log | M_{\alpha }
(\mu , x ) | }{1 + x^2 } \ d x \leq \pi \log 2 + 4 
\int\limits_{{\mathbb{R}}}^{} \frac{\tau (t)}{ 1 + t^2 } \ 
d t \ < \ \infty \ . $$

\noindent
Proposition 3.1 is proved. 
}}

\bigskip
\bigskip
{\bf 3.3. New version of M.G.Krein's Theorem and its 
corollaries.\ }

\bigskip
{\textsc{3.3.1.  Setting of a problem and M. Krein's theorem. 
}} \ \

Let $f \in {\mathcal{E}}^* ({\mathbb{R}}) \ $ and 
$\Lambda_f = \left\{\lambda_n \right\}_{n \geq 1 }  \ $. 
By Mittag-Leffler theorem [31, v.II] there exists such 
sequence $\left\{p_n \right\}_{n \geq 1 } \subset 
{\mathbb{Z}}_0 \ $
that 
$$ \sum_{n \geq 1 }^{} \frac{|z|^{p_n } }{|\lambda_n |^{1 + 
p_n } |f^{\prime } (\lambda_n )| } \ < \ \infty \ \ \ 
\forall \ z \in {\mathbb{C}} \ \eqno{(3.3.1)}$$

\noindent
and function
$$ \frac{1}{f(z)} - \sum_{n \geq 1 }^{}
\frac{z^{p_n } }{\lambda_n^{ p_n } f^{\prime } (\lambda_n )
(z - \lambda_n ) } \ $$

\noindent
is an entire function, where for $\lambda , z \in 
{\mathbb{C}} \ $ and positive integer $p \ $: 
$(z / \lambda )^p (z - \lambda )^{-1 } = (z - \lambda )^{-1 } 
+ \frac{1}{\lambda } + \frac{z}{\lambda^2 } + \ldots + 
\frac{z^{p-1} }{\lambda^{p}  } \ $. Assumptiom 
(see (3.1.6)) $d_f < + \infty \ $ gives possibility to set 
in (3.3.1) $p_n = p \geq 0 \vee d_f \ $
$\forall \ n \geq 1 \ $ and consider an entire function 
$\Delta_f^p (z) \ $ defined by (3.1.7). Just such assumption 
about the entire functions from more wide class $(A) \ $ 
(see 3.1.) have been made by M. Krein in [24]. But everywhere 
below we will consider M. Krein's results only on the set of 
real entire functions all zeros of which are real.  

So, for $f \in {\mathcal{E}}^* ({\mathbb{R}}) \ $ M. Krein in 
[24] made an assumption $d_f < + \infty \ $ and 
considered the problem of the description of all those 
functions $f \in {\mathcal{E}}^* ({\mathbb{R}}) \ $ entire 
function $\Delta_f^p (z) \ $ of which for some $p \geq 0 
\vee d_f \ $ is a polynomial.

M. Krein in [24] proved theorem which is described detally in 
[29], has a self-contained proof in [23] and has been 
discussed also in [34, 25, 9, 11]. We will use the following its 
version given in [9, Th.6.1].

\bigskip
\bigskip
{\textsc{Proposition $3.2. \ $ }} 
{\textit{Let $f \in {\mathcal{E}}_s ({\mathbb{R}} ) \ $ and 
the following relations hold: }}
$$ \sum_{\lambda \in \Lambda_f }^{} \frac{1}{|f^{\prime } 
(\lambda ) | } < \infty \ ; \ \ 
\frac{1}{f(z)} = \sum_{\lambda \in \Lambda_f }^{}
\frac{1}{f^{\prime } (\lambda ) (z - \lambda ) } \ , 
\ \forall \ z \in \ {\mathbb{C}} \setminus  \Lambda_f \ . $$

\noindent
{\textit{Then $f \in {\mathrm{Cartwright } } \ $. }}

\bigskip
Such form of M. Krein's theorem requires some additional 
comments. It was proved in [9, Lemma 6.3] that if  $f \in 
{\mathcal{E}}_s ({\mathbb{R}} ) \ $, $d_f < + \infty \ $ and 
for some  $p \geq 0 \vee d_f \ $ entire function $\Delta_f^p 
(z) \ $ is a polynomial then for any polynomial $Q \in 
{\mathcal{P}}_s ({\mathbb{R}}) \ $, satisfying $\Lambda_Q 
\cap  \Lambda_f = \emptyset \ $, ${\mathrm{deg}} Q > p \vee 
{\mathrm{deg}} \Delta_f^p \ $, function $g := Q \cdot f \ $
has the property $\Delta_g^0 (z) \equiv 0 \ $. That is why 
under these  conditions using Proposition $3.2 \ $ we get 
also $f \in {\mathrm{Cartwright } } \ $. So, taking into 
account that remark from the paper [9] and also indicated 
in 3.1 possibility to substitute inequality (3.1.2)
by (3.1.3) we can reformulate Proposition $3.2 \ $ as follows.

\bigskip
\bigskip
{\bf{ M. Krein Theorem.}}([24]) 
{\textit{Let $f \in {\mathcal{E}}_s ({\mathbb{R}} ) \ $ and 
$d_f < \infty \ $. If for some $p \geq 0 \vee d_f \ $ entire 
function $\Delta_f^p (z) \ $ is a polynomial then the 
following two properties hold:}}
$$\begin{array}{ll}
(3.3.2a) \ \ \ \ \  f \in {\mathcal{E}}_1 ({\mathbb{R}}) \  ; 
& \ \ \ \ \   
 (3.3.2b)\ \ \ \ \  \int\limits_{{\mathbb{R}} }^{} \frac{
 \left|\log |f(t)| \right| }{1 + t^2 } \ d t \ < \infty \ \ . 
  \end{array} $$

\noindent
 Note the following evident properties of the quantity $d_f 
 \  $ for $f \in {\mathcal{E}}_s ({\mathbb{R}} ) \ $:

\medskip
(3.3.3a) $d_{f\cdot Q} = d_f - d_Q \ \ \ \forall \ Q \in 
{\mathcal{P}}_s ({\mathbb{R}}) : \ \Lambda_Q \cap  \Lambda_f 
= \emptyset \ ; \ $

\medskip
(3.3.3b)$ d_{\frac{f}{Q}} = d_f + d_Q \ \ \ \forall \ Q 
\in {\mathcal{P}}_s ({\mathbb{R}}) : \ \Lambda_Q \subset   
\Lambda_f \ ; \ $

\medskip
\noindent
In the paper [9] the following Akhiezer's [2, III.11] remark 
was considered: if $f \in {\mathcal{E}}_s ({\mathbb{R}} ) \ $, 
$d_f \leq 1 \ $ and $\Delta_f^1 (z) \ $ is a polynomial then 
${\mathrm{deg}} \Delta_f^1 = 0 \  $. In the Corollary 6.4 
from [9] the more general result was proved: 
if $f \in {\mathcal{E}}_s ({\mathbb{R}} ) \ $, $d_f < \infty 
\ $ and $\Delta_f^p \in {\mathcal{P}} \ $ for some $p \geq 1 
\vee d_f \ $, then ${\mathrm{deg}} \Delta_f^p \leq p - 1 \ $,
and for $p = 0 \leq d_f \ $: $\Delta_f^0 (z) \equiv 0 \ $. 
Using this fact, Phragmen-Lindelof principle and possibility 
to differentiate series (3.1.7) it is easy to derive the 
validity of the following statement given here without proof.

\bigskip
\bigskip
{\textsc{Proposition 3.3. \  }} 
{\textit{Let $f \in {\mathcal{E}}_s ({\mathbb{R}} ) \ $,   
$d_f < + \infty \ $ and for some $p_0 \geq 0 \vee d_f \ $ 
entire function $\Delta_f^{p_0} (z) \ $ is a polynomial. 
Then (see (3.1.8))
}}

\medskip
(3.3.4a) $\frac{1}{f_a (z)} = m_{f_a}^p (z) \ \ \ 
\forall z \in {\mathbb{C}} 
\setminus \Lambda_{f_a } \ \ \ \forall p \geq 0  \vee d_f  \ 
\ \ \forall a \in {\mathbb{R}} \ $;

\medskip
(3.3.4b) $\frac{1}{f(z)Q(z)} = m_{f\cdot  Q}^p (z) \ 
\ \ \forall z \in {\mathbb{C}} \setminus \Lambda_{f \cdot Q } 
\ \ \ \forall p \geq 0 \vee (d_f - {\mathrm{deg}} Q )  \ 
\ \ \forall \ 
Q \in {\mathcal{P}}_s ({\mathbb{R}}) : \ \Lambda_Q \cap  
\Lambda_f = \emptyset \    $;

\medskip
(3.3.4c) $ \frac{Q(z)}{f(z)} = m_{\frac{f}{ Q} }^p (z) \ 
\ \ \forall z \in {\mathbb{C}} \setminus 
\Lambda_{ \frac{f}{ Q}} 
\ \ \ \forall p \geq 0 \vee (d_f + {\mathrm{deg}} Q )
\ \ \ \forall \ Q \in {\mathcal{P}}_s ({\mathbb{R}}) : \ 
\Lambda_Q \subset   \Lambda_f \ . \ $

 \medskip
\noindent
{\textit{where}} $f_a (z) := f(z + a) \ $, $a \in 
{\mathbb{R}} \ $, $z \in {\mathbb{C}} \ $.

\bigskip
That is why the following property of functions $f \in 
{\mathcal{E}}_s ({\mathbb{R}} ) \ $: $d_f < \infty \ $ and 
$\exists p \geq 0 \vee d_f \ $: $\Delta_f^p \in {\mathcal{P}} 
\ $, is invariant with respect to the translation $f (z + a) 
\ $, $a \in {\mathbb{R}} \ $, multiplication and division on 
the pointed out in (3.3.4b) and (3.3.4c) polynomials. 
              
\bigskip
{\textsc{3.3.2.Main results. }} \ \

Continuing consideration of the simple properties of the 
functions from the class
${\mathcal{E}} ({\mathbb{R}}) \cap {\mathrm{Cartwright}} \ $,
being started in [9, Th.6.2], we establish the following 
statement.

\bigskip
\bigskip
{\textsc{Lemma 3.1. \ }} {\textit{Let $f \in {\mathcal{E}}_1 
({\mathbb{R}}) \ $ and }}
$$ \sum_{\lambda \in \Lambda_f \setminus \left\{0\right\} }^{}
\frac{1}{|\lambda | } < \infty \ . \eqno{(3.3.5)} $$

\noindent
{\textit{Then the following statements hold:}}

\medskip
(3.3.5a) \ {\textit{if $d_f < +\infty \ $ and ${\mathrm{co}} 
\Lambda_f = {\mathbb{R}} \ $ then $f \in {\mathcal{E}}_0 
({\mathbb{R}}) \ $; }}

\medskip
(3.3.5b) \ {\textit{if $ f \in {\mathrm{Cartwright}} \ $ then 
$f \in {\mathcal{E}}_0 ({\mathbb{R}}) \ $; }}

\medskip
(3.3.5c) \ {\textit{if $ f \in {\mathrm{Cartwright}} \ $ and 
${\mathrm{co}} \Lambda_f \neq  {\mathbb{R}} \ $ then
$\sum\limits_{\lambda \in \Lambda_f \setminus 
\left\{0\right\} }^{}  \frac{\log^+ |\lambda | }{|\lambda | } 
 < \infty \ $. }}

\bigskip
{\textsc{Proof of Lemma 3.1. \ }}
{\normalsize{
By Hadamard's theorem $f(z) = e^{- a z } P(z) \ $, $a \in 
{\mathbb{R}} \ $, $z \in {\mathbb{C}} \ $, $P(z) = 
\prod_{\lambda \in \Lambda_f }^{} (1 - \frac{z}{\lambda } ) 
\ $ where without loss of generality we assume $0 \notin 
\Lambda_f \ $. If $a \neq 0 \ $ then change of variables 
$z \ $ by $-z \ $ allows us to consider only the case 
$a > 0 \ $.
It is known [42, 8.6.4., Ex.8.15] that $P \in {\mathcal{E}}_0 
\ $ and therefore $\exists C > 0 $: $|P(z)| \leq C 
e^{\frac{a}{2} x } \ $ $\forall \ z = x + i y $, $y \in I \ $,
$x \geq 0 \ $. Then $|f(x)|, |f^{\prime } (x)| \leq  
C e^{ - \frac{a}{2} (x - 1) } \ $ $\forall \ x \geq 1 \ $. 
Such inequalities lead to a contradiction with $d_f < 
\infty \ $, if ${\mathrm{co}} \Lambda_f = {\mathbb{R}} \ $ 
and with inequality (3.1.3),
if $ f \in {\mathrm{Cartwright}} \ $. So, $f \equiv P \in 
{\mathcal{E}}_0 ({\mathbb{R}}) \ $ and (3.3.5a), (3.3.5b) 
are proved. Prove at last (3.3.5c). Change of variables 
$f(b \pm  x) \ $, $b \in {\mathbb{R}} \ $, allows us to assume
$\Lambda_f = \left\{\lambda_n \right\}_{n \geq 1 } 
\subset [1, + \infty ) \ $. Then according to (3.3.2b) 
for arbitrary
$N \geq 1 \ $:
$$  \int\limits_{{\mathbb{R}} }^{} \frac{\left|\log |f(x)| 
\right| }{1 + x^2 } \ d x \ \geq \ \int\limits_{- \infty }^{0}
\frac{\log P(x)  }{1 + x^2 } \ d x \ \geq \ \sum_{n = 1 }^{N }
\int\limits_{0 }^{ \infty } \frac{\log (1 + \frac{x}{\lambda_k
} ) }{1 + x^2 } \ d x \ \geq \ (\log 2 ) \cdot  
\sum_{n = 1 }^{N } \frac{\log \lambda_k }{\lambda_k  } \ , $$

\noindent
what was to be proved. }} $\Box$

\bigskip
Since the Lindelof's theorem [29, I] implies validity of 
(3.3.5) for $f \in {\mathcal{E}}_1 ({\mathbb{R}}) \ $ 
with 
${\mathrm{co}} \Lambda_f \neq  {\mathbb{R}} \ $, then due 
to Lemma 3.1 we have validity of the following 
implication:
$$ f \in {\mathcal{E}}({\mathbb{R}}) \cap 
{\mathrm{Cartwright}} , \ {\mathrm{co}} \Lambda_f \neq  
{\mathbb{R}} \ \Rightarrow \ f \in {\mathcal{E}}_0 
({\mathbb{R}}) \ . \eqno{(3.3.6)} $$

\noindent
The following statement represents another version of 
M.G. Krein's theorem.

\bigskip
\bigskip
{\textsc{Theorem 3.1. \ }}  {\textit{Let $f \ $ be 
non-constant real entire function with only real and 
simple zeros and}}
$$ d_f \ := \ \inf \ \ \left\{q \in {\mathbb{Z}} \ | \ 
\sum\limits_{\lambda \in \Lambda_f \setminus \{ 0 \} }  
\frac{1}{|f^{\prime } (\lambda ) | \cdot |\lambda |^{q + 1}  } 
\ \ < \ \ \infty \ \ 
\right\}  \ < + \infty .  \eqno{(3.3.7)}  $$

\noindent
{\textit{The following statements are equivalent:}}

\medskip
(3.3.8a) \ {\textit{There exists such integer  $p \geq 0 
\vee d_f \ $
that entire function}}
$$ \frac{1}{f(z)} \ - \ \frac{\chi_{\Lambda_f } (0)  }{f^{
\prime } (0 ) 
\cdot z } \ - \ 
 \sum\limits_{\lambda \in \Lambda_f \setminus \left\{
 0\right\}    }
\frac{z^p }{ \lambda^p f^{\prime } (\lambda  ) (z - \lambda ) 
} \ $$

\noindent
{\textit{is a polynomial.}}

\medskip
(3.3.8b) \ {\textit{Function $f \ $ is an entire function of 
exponential type and  }}
$$ \int\limits_{{\mathbb{R}} }^{} \frac{\left| \ \log |f(t)| 
 \ \right| }{1 + t^2 } \ d t \ < \infty \ \ . $$

\medskip
(3.3.8c) \ {\textit{If ${\mathrm{co}} \Lambda_f = {\mathbb{R}}
\ $, then $f \ $ is an entire function of exponential type,
but if ${\mathrm{co}} \Lambda_f \neq  {\mathbb{R}} \ $,
then $f \ $ is an entire function of minimal exponential 
type.}}

\bigskip
Implication $(3.3.8a) \ \Rightarrow \ (3.3.8b) \ $ coincides 
identically with M. Krein's theorem. Implication 
$(3.3.8b) \ \Rightarrow \ (3.3.8a) \ $ has been proved by 
L.de Branges in his famous paper [12, Lemma 2], where one need
to take $G \in {\mathcal{P}} \ $ and observe that 
$|F(iy)| \ $ tends to infinity faster than any exponential 
function. Implication $(3.3.8b) \ \Rightarrow \ (3.3.8c) \ $
follows from (3.3.6). Implication $(3.3.8c) \ 
\Rightarrow \ (3.3.8a) \ $ for the entire functions of 
minimal exponential type was proved in the master's thesis
of Henrik L. Pederson at University of Copenhangen and can 
be found in [9] as Theorem 6.6. In view of Lemma 3.1 it is
remained to prove only those part of $(3.3.8c) \ 
\Rightarrow \ (3.3.8a) \ $  
where $f \in {\mathcal{E}}_1 ({\mathbb{R}}) \setminus 
{\mathcal{E}}_0 ({\mathbb{R}}) \  $, ${\mathrm{co}} 
\Lambda_f = {\mathbb{R}} \ $ and so by (3.3.5a) $ 
\sum_{\lambda \in \Lambda_f \setminus \left\{0\right\} }^{} 
\frac{1}{|\lambda | } = \infty \  $. This part can be easily
derived from the following theorem which will be proved in 
3.5.

\bigskip
\bigskip
{\textsc{Theorem 3.2. \ }} {\textit{Let real entire 
transcendental function $f \ $ has only real zeros and 
taking into account their multiplicity
$  \left\{ \lambda_k \  \left|  \ Q < k < P	\ 
\right. \right\} := \Lambda_f \setminus \left\{0\right\}\ $,
$ P, Q  \in { \mathbb{Z} } \cup  \{ \pm \infty \}
\ \ $, $\lambda_k \leq \lambda_{k+1 } \   $
$\forall  \ Q < k < P - 1 \ $. Let also exist such increasing 
sequences of positive real numbers $R_n , r_n \  $, $n \geq 1 \ $,
that $R_n , r_n \ \to + \infty \ $, $n  \to  \infty \ $, and
$$ f(z) = z^m \frac{f^{(m)} (0) }{m ! } \cdot \lim\limits_{n
\to \infty } \prod\limits_{\lambda \in (\Lambda_f \setminus 
\left\{0\right\})\cap (- r_n ,
R_n ) }^{} \left(1 - \frac{z}{\lambda } \right) \ \ \forall \ 
z \in {\mathbb{C}} \ , \eqno{(3.3.9)}$$
}}

\noindent
{\textit{where $m \in {\mathbb{Z}}_0 \ $. }}

{\textit{Then $f \in {\mathcal{LP}}_{II}^0 \  $ and there 
exist such sequences of integers $p_N , q_N \ $:
$Q < q_N < p_N < P \ $, $N \geq 1 \ $, that the polynomial 
divisors of the function $f \ $ which have the following 
form: }}

$$ P_N (z) := \frac{f^{(m)} (0) }{m ! } \cdot  z^m  \cdot 
\prod\limits_{k = q_N }^{p_N } \left(1 - \frac{z}{\lambda_k } 
\right) \ \eqno{(3.3.10)}$$

\noindent
{\textit{converge to $f(z) \ $ uniformly on any compact 
subset of $ {\mathbb{C}} \ $ and satisfy conditions:}}

\medskip
(3.3.11a)\  $(-N, N ) \cap \Lambda_f \subseteq \left\{
\lambda_k \right\}_{k = q_N}^{p_N } \  $;

\medskip
(3.3.11b) \  $|P_N (x) | \geq \frac{1}{e} \cdot |f(x)| \ 
 \ \ \forall \ x \in [\lambda_{q_N} ,  \lambda_{p_N }  ] \ $;

\medskip
(3.3.11c) \ $|P^{(m_k )}_N (\lambda_k ) | \geq \frac{1}{e} 
\cdot |f^{(m_k )} (\lambda_k )| \ \  \forall \ q_N \leq k 
\leq p_N  \ $,

\medskip
\noindent
{\textit{where $m_k \geq 1 \ $ denotes the multiplicity of zero 
$\lambda_k \in \Lambda_f  \setminus \left\{0\right\} \ $ 
$\forall  \ Q < k < P \ $ (in terms of the set $\Lambda_f \ $
this means that $m_k \ $ is a  number of the 
equal to $\lambda_k \ $ elements in  $\Lambda_f \ $).}}

\bigskip
\bigskip
{\textsc{Remark 3.1.}}({\emph{Sense of the condition (3.3.8c)}})
Representation (3.3.9) means in particular that function 
$f \ $ can be obtained not only as a limit of some 
sequence of real polynomials with real zeros  
but as a limit of its polynomial divisors. Consider an 
arbitrary  $f \in {\mathcal{E}}_s ({\mathbb{R}}) \cap 
{\mathcal{E}}_1 \ $ with $d_f < + \infty \ $ and 
clarify in what cases that function cannot be represented
as a limit of its polynomial divisors.  If $f \in 
{\mathcal{E}}_0 \ $ then (3.3.9) is a corollary of 
Lindelof's theorem. Let  $f \in 
{\mathcal{E}}_1  \setminus {\mathcal{E}}_0  \  $.

If $ \sum\limits_{\lambda \in \Lambda_f \setminus 
\left\{0\right\} }^{}  
\frac{1}{|\lambda | } = \infty \  $ then once more by 
Lindelof's theorem ${\mathrm{co}} \Lambda_f = {\mathbb{R}} \ $,
$ \sum\limits_{\lambda \in \Lambda_f , \lambda > 0  }^{}  
\frac{1}{|\lambda | } = \sum\limits_{\lambda \in \Lambda_f ,
\lambda < 0  }^{} \frac{1}{|\lambda | } = + \infty \  $
and  using the Hadamard's theorem we get for some
$a \in {\mathbb{R}} \ $, $m \in {\mathbb{Z}}_0 \ $ and any 
$R, r > 0 \ $:

$$ \frac{m ! f(z)}{f^{(m)} (0) z^m } = 
e^{- a z } \cdot \prod\limits_{\lambda \in \Lambda_f 
\setminus \left\{0\right\} }^{} \left(1 - \frac{z}{\lambda } 
\right)e^{\frac{z}{\lambda }} = 
e^{(\delta_f (r, R) - a ) z} \cdot f_{r, R } (z) \cdot 
\prod\limits_{\lambda \in (\Lambda_f \setminus 
\left\{0\right\})\cap (-r , R ) }^{} \left(1 - 
\frac{z}{\lambda } \right)  \ , 
$$

\noindent
where $ f_{r, R } (z) := \prod\limits_{\lambda \in \Lambda_f 
\setminus  (-r , R ) }^{} \left(1 - \frac{z}{\lambda } 
\right)e^{ \frac{z}{\lambda }} \ $, $\delta_f (r, R) := 
\sum\limits_{\lambda \in (\Lambda_f \setminus 
\left\{0\right\})\cap (-r , R ) }^{} \frac{1}{\lambda } 
\ $. Choosing two sequences $r_n , R_n \ $, $n \geq 1 \ $, 
so that $r_n , R_n \to + \infty \ $,
$\delta_f (r_n, R_n ) \to a \ $, $n \to \infty \ $, we 
get  representation (3.3.9).

If now  $ \sum_{\lambda \in \Lambda_f \setminus 
\left\{0\right\} }^{}  
\frac{1}{|\lambda | } < \infty \  $ then by (3.3.5a) 
${\mathrm{co}} \Lambda_f \neq  {\mathbb{R}} \ $ and by 
virtue of Hadamard's theorem $f (z) = e^{a z } f_0 (z) \ $,
$f_0 \in {\mathcal{E}}_0 \ $, $a \in {\mathbb{R}} \setminus 
\left\{0 \right\} \ $. Our condition $d_f < + \infty \ $ for
functions of such kind indicates only that $a > 0 \ $, if
$\sup \Lambda_f = + \infty \ $, and $a < 0 \ $, if 
$\inf \Lambda_f = - \infty \ $. But in both cases $f \ $ of
this kind cannot be represented as a limit of some its 
polynomial divisors. Just that class of entire functions 
has been excluded by condition (3.3.8c). 

That is why any $f \in {\mathcal{E}}_s ({\mathbb{R}}) \ $ which
satisfies (3.3.8c) and $d_f < + \infty \ $ can be 
represented in the form of (3.3.9). $\Box $

\bigskip
\bigskip
{\textsc{Proof of Theorem 3.1. \ }}
As it was noted above it is remained to prove only 
implication $(3.3.8c) \ \Rightarrow \ (3.3.8a) \ $
where in view of Remark 3.1 one can apply to the 
considered function $f \ $ the Theorem 3.2. But 
firstly we multiply $f \ $ on the polynomial $Q \ $ 
satisfying (3.3.3a) in order to
obtain $d_g \leq -1 \ $ for $g := f \cdot Q \ $. 
Approximating $g \ $ by polynomials $P_N \ $ from 
Theorem 3.2
we will get by (3.3.11c) for arbitrary $z \in {\mathbb{C}} 
\setminus \Lambda_g \ $ and $R > 0 \ $:

$$ \left|\frac{1}{P_N (z) } - \sum\limits_{\lambda \in 
\Lambda_{P_N } \cap (-R , R ) }^{} \frac{1}{P_N^{\prime } 
(\lambda ) (z - \lambda ) } \right| = \left| \sum\limits_{
\lambda \in \Lambda_{P_N } \setminus  (-R , R ) }^{} 
\frac{1}{P_N^{\prime } (\lambda ) (z - \lambda ) } \right| 
\leq $$

$$ \leq e \left( \sum\limits_{\lambda \in \Lambda_{g } 
\setminus  (-R , R ) }^{} \frac{1}{|g^{\prime } (\lambda )|} 
\right) \cdot  \sup_{ \lambda \in \Lambda_{g } \setminus  
(-R , R )}  \frac{1}{|z - \lambda  | }  \ 
\eqno{(3.3.12)}$$

\noindent
and passing to the limit as $N \to \infty \ $, we obtain 
for every $z \in {\mathbb{C}} \setminus  \Lambda_g \ $:

$$ \left|\frac{1}{g (z) } - \sum\limits_{\lambda \in 
\Lambda_{g } \cap (-R , R ) }^{} \frac{1}{g^{\prime } 
(\lambda ) (z - \lambda ) } \right| \ \to \ 0 , \ R \ \to \ 
+ \infty \ , \eqno{(3.3.13)}$$

\noindent
from where with the help of Proposition 3.3 we derive the 
required property (3.3.8a) for the function $f \ $. $\Box$

\bigskip
\bigskip
Theorems 3.1, 3.2 and Remark 3.1 give possibility to 
characterize Hamburger and Krein classes of entire 
functions (see 3.1) in terms of the behavior of the 
{\it{ derivative numbers}}  $\left\{ f^{\prime } (\lambda )
\right\}_{\lambda \in \Lambda_f } \ $ 
of the entire function $f \ $.

\bigskip
\bigskip
{\textsc{ Corollary 3.1. \ }}

\medskip
1. {\textit{Entire function $f(z) \ $ belongs to the Krein class
${\mathcal{K}} \ $ {\it if and only if } it has the 
following properties:  }}

\bigskip
(3.3.14a) $f \ $ {\textit{is a real function with only real and 
simple  zeros $\Lambda_f \  $}}
 \ ; 

\bigskip
(3.3.14b) {\textit{if ${\mathrm{co}} \Lambda_f = {\mathbb{R}} \  $ 
then $f \ $ is of 
exponential type, but if ${\mathrm{co}} \Lambda_f \neq  
{\mathbb{R}} \  $ 
then $f \ $ is of  minimal exponential type \ ; }}

\bigskip
(3.3.14c) $ \sum\limits_{\lambda  \in \Lambda_f } \ 
\frac{1}{(1 + {\lambda }^2 ) |f^{\prime } (\lambda )| }
\ < \ \ \infty \ \ $.

\bigskip
\bigskip
2.{\textit{ Entire function $f(z) \ $ belongs to the Hamburger 
class
${\mathcal{H}} \ $ {\it if and only if } it has the following 
properties:	}}

\bigskip
(3.3.15a) $f \ $ {\textit{is a real function with only real 
and simple 
zeros $\Lambda_f \  $ \ ; }}

\bigskip
(3.3.15b) {\textit{if ${\mathrm{co}} \Lambda_f = {\mathbb{R}} \  $ 
then $f \ $ is of 
exponential type,  but if ${\mathrm{co}} \Lambda_f \neq  
{\mathbb{R}} \  $ 
then $f \ $ is of minimal exponential type \ ; }}

\bigskip
(3.3.15c) $ \lim\limits_{\stackrel{|\lambda | \to \infty }{
\lambda  \in \Lambda_f 
}} \ \frac{|\lambda |^n }{ |f^{\prime } (\lambda )| }
\ = \ \ 0 \ \ \ \forall \ n \in {\mathbb{Z}}_0 \ $. \hfill 
$\bigtriangleup \ $

\bigskip
\bigskip
It should be noted that entire functions satisfying 
conditions (3.3.14a) and (3.3.14b) form a sufficiently 
large subset of the second Laguerre-P\'olya class 
${\mathcal{LP}}_{II} \ $ of entire functions  and each 
of them can be represented as follows:

$$ f(z) = \left\{ \begin{array}{lll}
c \cdot z^q \cdot \prod\limits_{k \geq 1}^{}
\left(1 - \frac{z}{\lambda_k }\right) \ ; &
\sum\limits_{k \geq 1 }^{} \frac{1}{|\lambda_k |} < \infty \ ,
&\mbox{if} \ {\mathrm{co}} \left\{\lambda_k \right\}_{
k \geq 1 } \neq {\mathbb{R}} \ ;   \\
c \cdot z^q \cdot e^{a z} \cdot  \prod\limits_{k \geq 1}^{}
\left(1 - \frac{z}{\lambda_k }\right) e^{\frac{z}{\lambda_k } 
} \ ; & a \in {\mathbb{R}} , \ \sum\limits_{k \geq 1 }^{} 
\frac{1}{\lambda_k^2 } < \infty \ ,   &
\mbox{if} \ {\mathrm{co}} \left\{\lambda_k \right\}_{
k \geq 1 } = {\mathbb{R}} \ ;   \end{array}\right. 
\eqno{(3.3.16)} $$

\noindent
where $c, \lambda_k \in {\mathbb{R}}\setminus 
\left\{0\right\} \ $ $\forall \ k \geq 1 \ $,
$\lambda_k \neq \lambda_m \ $, if $k \neq  m \ $, $k, 
m \geq 1 \ $, $q \in {\mathbb{Z}}_0 \ $. But if we 
impose on the function $f \ $ from that class  
only one condition on their derivative
numbers: $d_f < + \infty \ $, or, what is the same,

$$ \overline{\lim\limits_{\stackrel{|\lambda | \to 
\infty }{\lambda \in \Lambda_f }} } \ \  \frac{\log^+ 
\frac{1}{|f^{\prime }(\lambda )| }}{ \log |\lambda | } \ < 
+ \infty \ , \eqno{(3.3.17)}$$

\noindent
then by Theorem 3.1 we can conclude that this function $f \ $
will be an element of \\ Cartwright class, what in the case
${\mathrm{co}} \left\{\lambda_k \right\}_{k \geq 1 } 
\neq {\mathbb{R}} \ $ means by virtue of Lemma 3.1 that 
$\sum_{\lambda \in \Lambda_f \setminus \left\{0\right\} }^{} 
\frac{\log^+ |\lambda_k |}{|\lambda_k | } < \infty \ $ and 
in the case ${\mathrm{co}} \left\{\lambda_k \right\}_{k 
\geq 1 } = {\mathbb{R}} \ $ gives an existence of the limit
$\delta_f := \lim\limits_{ R \to + \infty } \sum_{\lambda 
\in (-R, R ) \cap \left( \Lambda_f  \setminus
\left\{0\right\} \right) }^{} \frac{1}{\lambda } \  $, 
equality $a + \delta_f = 0 \ $ and also an existence and 
equality of two finite limits: $ \lim\limits_{R \to + 
\infty }
\frac{{\mathrm{card}} \left( \Lambda_f \cap [0,R]\right) }{
R} \ = \ 
\lim\limits_{R \to + \infty }
\frac{{\mathrm{card}} \left( \Lambda_f \cap [-R, 0] 
\right) }{R} \  $ (see [29, V.4, Th.11]). 

\bigskip
Besides that Theorem 3.2 and Remark 3.1 allow us to establish 
direct and inverse polynomial approximation  theorem for 
the entire functions from  Hamburger and Krein classes. 
It should be noted here that since both classes 
${\mathcal{H}} \ $ and ${\mathcal{K}} \ $ are subsets of 
the second Laguerre-P\'olya class ${\mathcal{LP}}_{II} \ $ then 
every function from these classes can be approximated 
[18, III, Th.3.2] by real polynomials with real zeros only 
with respect to the topology $\tau_{{\mathcal{E}}} \ $.

\bigskip
\bigskip
{\textsc{ Corollary 3.2. \ }}

\medskip
1. {\textit{Entire function $f(z) \ $ belongs to the Krein 
class
${\mathcal{K}} \ $ {\it if and only if } there exists the 
sequence of real 
polynomials $\left\{P_n \right\}_{n \geq 1 } \ $ with only 
real and simple 
zeros $\left\{\Lambda_{P_n }\right\}_{n \geq 1 } \ $ which 
uniformly on any 
compact subset of the complex plane converges to the 
function $f \ $ and for 
some does not depending on $n \ $ constant $C \ > \ 0 \ $ 
the following 
inequality holds:}}

 $$ \sum\limits_{\lambda  \in \Lambda_{P_n } } \ 
\frac{1}{(1 + {\lambda }^2 ) |P_n^{\prime } (\lambda )| }
\ \leq \ \ C \ \ \forall \ n \geq \ 1 \ .\eqno{(3.3.18)} $$

\medskip 
2. {\textit{Entire function $f(z) \ $ belongs to the 
Hamburger class
${\mathcal{H}} \ $ {\it if and only if } there exists 
the sequence of real 
polynomials $\left\{P_n \right\}_{n \geq 1 } \ $  with 
only real and simple 
zeros $\left\{\Lambda_{P_n }\right\}_{n \geq 1 } \ $ which 
uniformly on any 
compact subset of the complex plane converges to the 
function $f \ $ and for 
some does not depending on $n \ $ function
$ w : {\mathbb{R}} \ \to \ (0, + \infty ) \ $, 
$ \sup\limits_{x \in {\mathbb{R}} } |x|^n w (x)  < 
+ \infty \ \ $ 
$\forall \ n \in {\mathbb{Z}}_0 \ $, the following 
inequality holds: }}

$$ 
\left|P^{\prime }_n (\lambda )\right| \ \geq \ \ 
\frac{1}{w(\lambda )} \ \  
\ \ \forall \ \lambda \in  \Lambda_{P_n } \ \ \ 
\forall \ n \geq 1 \ . \eqno{(3.3.19)} $$

\medskip
\noindent
{\textit{Moreover, in both items the approximating 
polynomial sequence  can be chosen as the subset of all
polynomial
divisors of the function $f(z) \ $. }}

\bigskip
\bigskip
\bigskip
{\textsc{ Proof of Corollary 3.2. \ }}
Necessity  follows easily from the Theorems 3.1, 3.2 and 
Remark 3.1.

{\textit{Sufficiency. }} 1. Multiplying each polynomial 
$P_n \ $, $n \geq 1 \ $, on the polynomial of the second 
degree
$Q_2 \in {\mathcal{P}}_s ({\mathbb{R}}) \ $ satisfying
$\Lambda_{Q_2 } \cap \Lambda_f = \emptyset \  $ we by 
Hurwitz's
theorem obtain $\Lambda_{Q_2 } \cap \Lambda_{P_n } = 
\emptyset \  $ for sufficiently large $n \ $. Resulting 
polynomial sequence converges to $g := Q_2 \cdot f \ $,  
where
$d_g \leq -1 \ $. Performing similar to (3.3.12)  
estimate with $R > 2 |z| \ $ we get after passage to the 
limit as $n \to \infty \ $  relation (3.3.13) which by 
Proposition 3.3 yields $f \in {\mathcal{K}} \ $.

2. Here one need to use P\'olya-Laguerre [29, VIII.1, Th.3] 
theorem according to which for some does not depending on 
$n \ $ constant $M > 0 \ $: $\sum\limits_{\lambda \in 
\Lambda_{P_n } }^{} \frac{1}{1 + \lambda^2 } \leq M \ $ 
$\forall \ n \geq 1 \ $.  We can perform a similar to 
(3.3.12) estimate from where taking into account the 
following corollary of (3.3.19):

$$ \sum\limits_{\lambda \in \Lambda_{P_n } \setminus  
(-R , R ) }^{} \frac{1}{|{P_n }^{\prime } (\lambda )|} 
\leq 
\sum\limits_{\lambda \in \Lambda_{P_n } \setminus  
(-R , R ) }^{} \frac{(1 + \lambda^2 ) w(\lambda ) }{1 + 
\lambda^2 }
\leq M \left\|(1 + x^2 ) w \right\|_{C ({\mathbb{R}})} \ ,
$$

\noindent
one can easily obtain for $R > 2|z| \ $ (3.3.13) and then 
by
(3.3.19) $f \in {\mathcal{H}} \ $. $\Box$

\bigskip
\bigskip
{\bf 3.4. Strictly normal polynomial families.\ }

\bigskip
{\textsc{3.4.1. Main results.}}
Recall (see 3.1) that ${\mathcal{P}}^* ({\mathbb{R}}) \ $ 
denotes the set of real polynomials $P \ $ with only real 
and simple zeros and $P(0) = 1 \ $.

Let $G \subseteq {\mathcal{P}}^* ({\mathbb{R}}) \ $ be a 
normal family of polynomials (see 3.1). Making use a proof 
by contradiction it is easy to derive succeeding the proof 
of P\'olya-Laguerre [18, III, Th.3.3] theorem that 
$ \lambda_1 (G) := \sup_{P \in G } |P^{\prime } (0)| = 
\sup_{P \in G } \left|\sum_{\lambda \in \Lambda_P }^{} 
\frac{1}{\lambda } \right| < \infty \ $ and
$ \lambda_2 (G) := \sup_{P \in G } \left(P^{\prime } (0)^2 
- P^{\prime \prime } (0) \right) = \sup_{P \in G } \sum_{
\lambda \in \Lambda_P }^{} \frac{1}{\lambda^2 } < \infty \ $.
Conversely, if these quantities are finite then an obvious 
inequality $\left|(1 - z) e^{z} \right| \leq e^{\frac{1}{2} 
|z|^2 } \ $$\forall \ z \in {\mathbb{C}} \ $
implies $|P(z)| \leq \exp \left(M_1 |z| + \frac{1}{2} M_2 
|z|^2 \right) \ $$\forall \ z \in {\mathbb{C}} \ $, i.e.
by  Vitali's classical compactness theorem $G \ $ is a 
normal set. Thus, we have proved the following statement 
(see also [29, VIII.1]).

\bigskip
\bigskip
{\textsc{ Proposition 3.4. \ }} {\textit{Polynomial family 
$G \subset {\mathcal{P}}^* ({\mathbb{R}}) \ $ is a normal 
set if and only if the following two conditions hold:}}

\medskip
(3.4.1a) \ $ \lambda_1 (G) := \sup\limits_{P \in G } 
\left|\sum\limits_{\lambda \in \Lambda_P }^{} \frac{1}{
\lambda } \right| < \infty \ $; 
\  \ 
(3.4.1b) \ $ \lambda_2 (G) :=  \sup\limits_{P \in G } 
\sum\limits_{\lambda \in \Lambda_P }^{} \frac{1}{\lambda^2 } 
< \infty \ $.

\bigskip
By well-known P\'olya-Laguerre theorem [18, III, Th.3.3] the 
closure of normal set  $G \subset {\mathcal{P}}^* ({
\mathbb{R}}) \ $ denoted as ${\mathrm{Close}}_{{\mathcal{E}}}
G \ $ is a compact subset of the second Laguerre-P\'olya class of 
the entire functions ${\mathcal{LP}}_{II} \ $ [18, III, Def.3.1].

\bigskip
\bigskip
{\textsc{ Definition 3.1. \ }} 
{\it{Normal family 
$G \subset {\mathcal{P}}^* ({\mathbb{R}}) \ $ 
is said to be a {\textit{strictly normal polynomial 
family}} if for any convergent with respect to the 
topology 
$\tau_{{\mathcal{E}}} \ $ sequence $\left\{P_n \right\}_{
n \geq 1 } \subseteq G \ $ satisfying $\lim\limits_{n \to 
\infty } {\mathrm{deg}} P_n = \infty \ $ an entire function 
$\lim\limits_{n  \to \infty } P_n (z)\ $ is  transcendental. }}  

\bigskip
In terms of the introduced by P. Painleve [26, II.29] notion 
of upper limit of the sequence $\left\{A_n \right\}_{n \geq 
1 } \ $ of subsets of some topological space with topology 
$\tau \ $: 
$${\mathrm{Ls}}_{n \to \infty } A_n := 
\bigcap\limits_{n \geq 1 }^{} {\mathrm{Close}}_{\tau } 
\left(\bigcup\limits_{k \geq n }^{} A_k \right) \ , $$

\noindent
where ${\mathrm{Close}}_{\tau } A \ $ denotes the closure 
of $A \ $ in the considered topological space, Definition 
3.1 means that $G \subset {\mathcal{P}}^* ({\mathbb{R}}) \ $
is a strictly normal polynomial family if and only if 
$G \cap {\mathrm{Ls}}_{n \to \infty } \left(G_{n + 1 } 
\setminus G_n \right) \ = \ \emptyset \ $, or, what is the 
same,
$$G \cap  \left[\bigcap\limits_{n \geq 1 }^{} {\mathrm{
Close}}_{{\mathcal{E}} } \left( G \setminus G_n \right) 
\right] = \emptyset \ , \ \ 
G_n := \left\{P \in G \ \left| \ {\mathrm{deg}} P \leq n 
\right.\right\} , \ \ n \geq 1 \ . \eqno{(3.4.2)}$$

\noindent
It is easy to verify that the closure ${\mathrm{Close}}_{{
\mathcal{E}} } G \ $ of the strictly normal polynomial 
family $G \ $ aside from the compactness property have 
one more  characteristic one: the set $G \ $ being 
considered as a subset of the topological space 
${\mathrm{Close}}_{{\mathcal{E}} } G \ $ with induced  
topology (from the whole space of all entire functions 
with topology $\tau_{{\mathcal{E}}} \ $) is an open set, 
i.e.
$$ {\mathrm{Close}}_{{\mathcal{E}} }  \left(\left({\mathrm{
Close}}_{{\mathcal{E}} } G  \right)\setminus G \right) \ = 
\ \left({\mathrm{Close}}_{{\mathcal{E}} } G \right)
\setminus G   \ . \eqno{(3.4.3)}$$

\noindent
In other words (3.4.3) means that any convergent sequence 
of the transcendental entire functions $\left\{f_n \right\}_{
n \geq 1  } \subseteq {\mathrm{Close}}_{{\mathcal{E}} } G
\ $
can have in capacity of its limit only also transcendental 
entire function or, what is the same, the set $\left({
\mathrm{Close}}_{{\mathcal{E}} } G  \right)\setminus G \ $ 
of all transcendental entire functions from  ${\mathrm{
Close}}_{{\mathcal{E}} } G \ $ is a closed and hence, 
compact set. That is why the closure of any strictly normal 
polynomial set $G \ $ with respect to topology $\tau_{{
\mathcal{E}}} \ $ generates at once two compact sets: 
${\mathrm{Close}}_{{\mathcal{E}} } G \ $ and $\left({
\mathrm{Close}}_{{\mathcal{E}} } G  \right)\setminus G \ $.

In contrast to normality criterion of the Proposition 3.4
we will be interested here in those sufficient conditions 
for the normality and strictly normality of the polynomial 
set $G \subset {\mathcal{P}}^* ({\mathbb{R}}) \ $ which can 
be formulated in terms of derivative numbers $\left\{P^{
\prime } (\lambda )\right\}_{\lambda \in \Lambda_P } \ $ 
of the polynomials from that set and which would give 
possibility to exclude condition of the type (3.4.1a) 
at all and to make condition of the type (3.4.1b) a 
little weaker.

\bigskip
\bigskip
{\textsc{ Lemma 3.2. \ }}
{\textit{For arbitrary finite constants $\alpha , \beta , 
\gamma , \delta_{\alpha } ,  \delta_{\beta } > 0 \ $ 
the set}}
$$\left\{P \in {\mathcal{P}}^* ({\mathbb{R}}) \ \left| \ 
\right.  \sum\limits_{\lambda \in \Lambda_P }^{} e^{- \alpha 
|\lambda |} \leq \delta_{\alpha } \ ; \ \ 
|P^{\prime } (\lambda )| \geq \delta_{\beta } e^{- 
\beta  |\lambda |} |\lambda |^{-1 -\gamma } \ \forall \ 
\lambda \in \Lambda_P \  \right\} \eqno{(3.4.4)}$$ 

\noindent
{\textit{is normal with respect to topology of the 
uniform convergence on all compact subsets of the complex 
plane (see 3.1). }}

\bigskip
The sequence of polynomials $\left\{1 - n x \right\}_{n 
\geq 1 } \ $ shows that the set (1) for $\gamma = 0 \ $  
is not  normal. Denote by $C_{+}^* ({\mathbb{R}}) \ $ the 
family of all positive functions from $C^* ({\mathbb{R}}) \ $
(see 2.1), i.e. the set of all upper semicontinuous functions
$ \mu : {\mathbb{R}} \to (0, + \infty )  \ $, satisfying 
conditions: $\left\|x^n \cdot \mu \right\|_{C 
({\mathbb{R}}) } < \infty \ $ $\forall \ n \in 
{\mathbb{Z}}_0 \ $.

\bigskip
\bigskip
{\textsc{ Theorem 3.3. \ }} {\textit{For any $\mu \in 
C_{+}^* ({\mathbb{R}}) \ $ and arbitrary finite constants
$\alpha , \gamma , \delta_{\alpha } > 0 \ $ the set}}
$$\left\{P \in {\mathcal{P}}^* ({\mathbb{R}}) \ \left| \ 
\right.  \sum\limits_{\lambda \in \Lambda_P }^{} e^{- 
\alpha |\lambda |} \leq \delta_{\alpha } \ ; \ \ 
|P^{\prime } (\lambda )| \geq  \frac{1}{\mu (\lambda ) |
\lambda |^{1 + \gamma } } \ \forall \ \lambda \in \Lambda_P \ 
\right\} \eqno{(3.4.5)}$$  

\noindent
{\textit{is a strictly normal polynomial set (see 
Definition 3.1).}}

\bigskip
\bigskip
Let ${\mathcal{H}}^* := \left\{f \in {\mathcal{H}} \ 
\left| \  f (0) = 1 \right.\right\} \ $ (${\mathcal{P}}^* 
({\mathbb{R}}) \subset {\mathcal{H}}^* \ $) and for $\mu 
\in C_{+}^* ({\mathbb{R}}) \ $, $\gamma , C_{\gamma } \in 
(0, + \infty ) \ $, denote
$$ {\mathcal{P}}^*_{{\mathcal{H}}} (C_{\gamma } , \mu ) :=
\left\{P \in {\mathcal{P}}^* ({\mathbb{R}}) \ \left| \ 
\right.  \sum\limits_{\lambda \in \Lambda_P }^{} e^{-  
|\lambda |} \leq C_{\gamma  } \ ; \ \ 
|P^{\prime } (\lambda )| \geq  \frac{1}{\mu (\lambda ) 
|\lambda |^{1 + \gamma } } \ \forall \ \lambda \in 
\Lambda_P \  \right\} \ . \eqno{(3.4.6)}$$  

\noindent
The set (3.4.6) by Theorem 3.3 is a strictly normal 
polynomial set. In addition, by virtue of Corollary 3.2:
$ {\mathcal{H}}^* (C_{\gamma } , \mu ) := {\mathrm{
Close}}_{\mathcal{E}} {\mathcal{P}}^*_{{\mathcal{H}}} 
(C_{\gamma } , \mu )
\subseteq {\mathcal{H}}^* \ $ and so:
$$ {\mathcal{H}}^* (C_{\gamma } , \mu ) = \left\{f \in 
{\mathcal{H}}^*  \ \left| \ \right.  \sum\limits_{\lambda 
\in \Lambda_f }^{} e^{-  |\lambda |} \leq C_{\gamma  } \ ; \ \ 
|f^{\prime } (\lambda )| \geq  \frac{1}{\mu (\lambda ) 
|\lambda |^{1 + \gamma } } \ \forall \ \lambda \in 
\Lambda_f \  \right\} \ . \eqno{(3.4.7)}$$ 

\noindent
The set (3.4.7) being a compact subset of the Hamburger 
class of entire functions ${\mathcal{H}} \ $ possesses 
due to (3.4.3) the following property.

\bigskip
\bigskip
{\textsc{ Corollary 3.3. \ }} {\textit{Let ${\mathcal{H}
}_{\infty} \ $
denote the set of transcendental entire functions from 
the Hamburger class ${\mathcal{H}} \ $ and for $\mu \in 
C_{+}^* ({\mathbb{R}}) \ $, $\gamma , C_{\gamma } \in (0, 
+ \infty ) \ $: }}
$$ {\mathcal{H}}^*_{\infty} (C_{\gamma } , \mu ) :=
{\mathcal{H}}_{\infty} \cap {\mathcal{H}}^* (C_{\gamma } , 
\mu ) \ . \eqno{(3.4.8)} $$

\noindent
{\textit{Then the set ${\mathcal{H}}^*_{\infty} (C_{\gamma } ,
\mu ) \ $ is a compact subset of ${\mathcal{H}}_{\infty} \ $ 
and }}
$$  {\mathcal{H}}^*_{\infty} (C_{\gamma } , \mu ) =
\left\{f \in {\mathcal{H}}_{\infty }  \left|  \right.  
f(0)= 1  ,  \sum\limits_{\lambda \in \Lambda_f }^{} e^{-  
|\lambda |} \leq C_{\gamma  }  ;  \ 
|f^{\prime } (\lambda )| \geq  \frac{1}{\mu (\lambda ) 
|\lambda |^{1 + \gamma } } \ \forall  \lambda \in \Lambda_f 
\  \right\}  . \eqno{(3.4.9)}$$ 

\bigskip
It should be noted at last that the statements of Lemma 3.2, 
Theorem 3.3 and Corollary 3.3 remain valid if we substitute
the conditions on the derivative numbers of the entire 
functions $f \ $ in (3.4.4, 5, 6, 7, 9) by the following 
ones: 
$$ \sum\limits_{\lambda \in \Lambda_f  }^{} \frac{1}{\mu 
(\lambda ) |\lambda |^{\beta } |f^{\prime } (\lambda )|^{
\gamma } } \ \leq \ C(\beta , \gamma )\ , \ \ \ \ \beta > 
\gamma
> 0 \ . \eqno{(3.4.10)}$$

\noindent
Such substitution is possible because for any $\mu \in 
C_{+}^* ({\mathbb{R}}) \ $: $(\mu )^{\alpha } \in C_{+}^* 
({\mathbb{R}}) \ $ $\forall \ \alpha \ > 0 \ $, and an 
arbitrary subset of normal or strictly normal polynomial 
set possesses also the corresponding property.

\bigskip
{\textsc{3.4.2. Proof of Lemma 3.2.}} 
{\normalsize{
Let $M \geq \alpha + \beta + 1 \ $ and $P \ $ be an 
arbitrary polynomial of the defined by (3.4.4) polynomial 
family.
Since for $z \in  {\mathbb{C}} \setminus \left\{ \alpha_n  
\right\}_{n \geq 1 } \ $, $\alpha_n := \frac{\pi }{M } (n - 
\frac{1}{2} )  \ $, $n \in {\mathbb{Z}} \ $, the following 
equality holds:
$$ \frac{1}{\cosh (M z ) } = 1 + \sum\limits_{n \geq 1 }^{}
\frac{(-1)^n }{ M \cdot \alpha_n } \left(\frac{z}{z + i 
\alpha_n } + \frac{z}{z - i \alpha_n } \right) \ , $$

\noindent
then, denoting $\left\{\lambda_k \right\}_{k =1 }^{n} := 
\Lambda_P \ $, $0 < |\lambda_1 | \leq |\lambda_2 | \leq 
\ldots \leq |\lambda_m | \leq 1 \leq |\lambda_{m + 1 } | 
\leq \ldots \leq |\lambda_N  | \ $, $0 \leq  m \leq N \ $, 
($m = 0 \ $, if $|\lambda_1 | > 1 \ $) we obtain for $z \in  
{\mathbb{C}} \setminus \left( \Lambda_P \cup  \left\{ 
\alpha_n  \right\}_{n \geq 1 } \right) \ $:
$$ \Phi (z) := \frac{1}{P(z) \cdot \cosh (M z ) } = 
\sum\limits_{k = 1}^{N} \frac{1}{P^{\prime } (\lambda_k ) 
\cosh (M \lambda_k ) } \frac{1}{z - \lambda_k } + \frac{i}{M}
\sum\limits_{n \geq 1 }^{} (-1)^n \times  $$
$$ \times \left[\frac{1}{P(i \alpha_n  
 ) (z - i \alpha_n ) } + \frac{1}{ P(- i \alpha_n  
 ) (z + i \alpha_n ) } \right] \ . \eqno{(3.4.11)} $$

\noindent
Differentiating equality (3.4.11) we have $\Phi^{\prime } (0)
= - P^{\prime } (0) \ $, $\Phi^{\prime \prime } (0) = 2 
P^{\prime  } (0)^2 - P^{\prime \prime } (0) - M^2 \ $ and
{\small{
$$ P^{\prime } (0) \ = \ \sum\limits_{k = 1 }^{N } \frac{1}{
\lambda_k^2  \cdot P^{\prime } (\lambda_k  ) \cosh (M 
\lambda_k ) } + \frac{i}{M } \sum\limits_{n \geq 1 }^{} 
\frac{(-1)^n }{\alpha_n^2  } \left(\frac{1}{P (- i \alpha_n ) 
} - \frac{1}{P(i \alpha_n ) } \right) \ , $$
$$ P^{\prime \prime  } (0) + M^2 - 2 P^{\prime }(0)^2  =  
\sum\limits_{k = 1 }^{N } \frac{1}{
\lambda_k^3  \cdot P^{\prime } (\lambda_k  ) \cosh (M 
\lambda_k ) } + \frac{1}{M } \sum\limits_{n \geq 1 }^{} 
\frac{(-1)^n }{\alpha_n^3  } \left[  \frac{1}{P(i \alpha_n )} 
+ \frac{1}{P (- i \alpha_n ) } \right]  , $$
}}

\noindent
from where taking into account $|P(i \lambda )| \geq 1 \ $
$\forall \ \lambda \ \in {\mathbb{R}} \ $ and 
$ \sum\limits_{n \geq 1 }^{} \frac{1}{(n - 0,5)^3 } < 
\frac{\pi^3}{2} \ $, we derive

$$ |P^{\prime } (0)| \ \leq \  M + \sum\limits_{k = 1 }^{N } 
\frac{1}{ \lambda_k^2  \cdot |P^{\prime } (\lambda_k  )| 
\cosh (M \lambda_k ) }  \ , \eqno{(3.4.12a)}$$
$$ 2 P^{\prime }(0)^2 - P^{\prime \prime  } (0)  \ \leq  \ 
2 M^2 + \sum\limits_{k = 1 }^{N } \frac{1}{
|\lambda_k |^3  \cdot | P^{\prime } (\lambda_k  )| \cosh 
(M \lambda_k ) } \ . \eqno{(3.4.12b)}$$

If $m = 0 \ $, then equalities (3.4.12a), (3.4.12b) and 
conditions (3.4.4) give estimates of $|P^{\prime } (0)| = 
\left|\sum\limits_{k = 1}^{N} \frac{1}{\lambda_k } \right| 
\ $ and $P^{\prime } (0)^2 - P^{\prime \prime } (0) = 
\sum\limits_{k = 1}^{N} \frac{1}{\lambda_k^2 }\ $ depending 
on five constants of Lemma 3.2 only.

If $m \geq 1 \ $ then (3.4.12a) and (3.4.4) yield:
$$ \left|\sum\limits_{k = 1}^{N} \frac{1}{\lambda_k } 
\right|
\leq M + \frac{2}{\delta_{\beta } } \sum\limits_{k = 1 }^{
m   } \frac{1}{|\lambda_k | } |\lambda_k |^{\gamma } e^{- (
1 + \alpha ) |\lambda_k | } + 
\frac{2}{\delta_{\beta } } \sum\limits_{k = m + 1 }^{N   }  
|\lambda_k |^{\gamma - 1 } e^{- (1 + \alpha ) |\lambda_k | }
\leq $$
$$\leq M + 2 \frac{\delta_{\alpha } }{\delta_{\beta } } 
\frac{1}{|\lambda_1 | } + 2 \frac{\delta_{\alpha } }{\delta_{
\beta } } \left(\frac{\gamma }{e } \right)^{\gamma } $$

\noindent
and therefore, using inequality $\log (1 + x ) \leq x \  $
$\forall \ x > -1 \ $, we get:

$$ \frac{\delta_{\beta} e^{- \beta }  }{ |\lambda_1 |^{1 + 
\gamma } } \leq |P^{\prime } (\lambda_1 ) | = 
\frac{1}{|\lambda_1 | } \prod\limits_{k = 2}^{N } 
\left(1 - \frac{\lambda_1 }{\lambda_k } \right) \leq \frac{e
}{|\lambda_1 | } e^{- \lambda_1 \left(\sum\limits_{k = 1}^{N} 
\frac{1}{ \lambda_k } \right) } \leq $$
$$ \leq \frac{e}{|\lambda_1 | } {\mathrm{exp}}{ \left[ M + 2 
\frac{\delta_{\alpha } }{\delta_{\beta } } \left(1 + \left(
\frac{\gamma }{e } \right)^{\gamma } \right) \right] }
 \ . $$

\noindent
That is why there exists such constant $\delta > 0 \ $ 
depending on constants (3.4.4) only that $|\lambda_1 | 
\geq \delta \ $. Since among all zeros of $P \ $ zero  $
\lambda_1 \ $ has a minimal absolute value then it follows 
from
(3.4.12a), (3.4.12b) and conditions (3.4.4) that there exist
the estimates of $\left|\sum\limits_{k = 1}^{N} \frac{1}{
\lambda_k } \right| \ $ and $ \sum\limits_{k = 1}^{N} 
\frac{1}{\lambda_k^2 }\ $ depending on  $\alpha , \beta , 
\gamma , \delta_{\alpha } ,  \delta_{\beta } > 0 \ $ 
only. This means that conditions of the Proposition 3.4 are 
fulfilled and so the set (3.4.4) is normal. }}

\bigskip
{\textsc{3.4.3. Proof of Theorem 3.3.}} 
{\normalsize{
Denote defined in (3.4.5) set by $G \ $. Since $\mu \ $ is 
uniformly bounded on the whole real axis then by Lemma 3.2 
$G \ $ is normal and due to Proposition 3.4 : $ M := 
\lambda_1 (G) \vee \lambda_2 (G) < + \infty \ $. 

Consider an arbitrary convergent to some entire function 
$f \ $ polynomial sequence $\left\{P_n \right\}_{n \geq 1 } 
\subseteq G \ $, with $\lim_{n \to \infty }{\mathrm{deg}} 
P_n = \infty \ $. Denote $\left\{\lambda_k^{(n)} \right\}_{
k=1}^{r_n } := \Lambda_{P_n } \ $, $\frac{1}{M} \leq |
\lambda_1^{(n)} | \leq |\lambda_2^{(n)} | \leq \ldots \leq |
\lambda_{r_n}^{(n)} | < \infty \  $ $\forall \ n \geq 1 \ $. 
By Hurwitz's theorem for arbitrary $p \geq 1 \ $ the sequence 
$\left\{ \lambda_p^{(n)} \right\}_{ n \geq n_p} \ $,$ \ 
r_{n_p} \geq p \ $, has finite or infinite limit. Let for $ 
n \geq n_p \ $:
$$ P_{n, p } (x) := \frac{P_n (x)}{ \Delta_{n, p} (x) } \ , 
\ \ \Delta_{n, p} (x) := \left( 1 - \frac{x}{\lambda_1^{(n)}} 
\right) \left(1 - \frac{x}{\lambda_2^{(n)}} \right)
\cdot \ldots \cdot \left(1 - \frac{x}{\lambda_{p-1}^{(n)}} 
\right) \ , $$

\noindent
where $ p \geq 2 , \ \Delta_{1, n} (x) \equiv 1 $. Then for 
all $p \leq k \leq r_n \ $, $n \geq n_p \ $:
$$|P^{ \prime } (\lambda_k^{(n)} ) | = 
\left| \Delta_{p, n } (\lambda_k^{(n)} ) \right| \cdot 
\left| P^{\prime }_{n , p } (\lambda_k^{(n)} ) \right| 
\leq  \left(1 + M |\lambda_k^{(n)} | \right)^{p - 1}\cdot 
\left| P^{\prime }_{n , p } (\lambda_k^{(n)} ) \right| \ . $$

\noindent
Using (3.4.5) and decomposition formula of $P_{n, p} (
z)^{-1} \ $ on the simple fractions we get:
$$ 1 \leq \sum\limits_{k = p}^{r_n } \frac{1}{| P^{\prime }_{n,
p } (\lambda_k^{(n)} ) | } \frac{1}{|\lambda_k^{(n)} |} \leq 
\sum\limits_{k = p }^{r_n} \mu (\lambda_k^{(n)} ) |
\lambda_k^{(n)}|^{\gamma } \leq   \ \frac{M}{|\lambda_p^{(
n)}|} \left\|x^{1 + \gamma } (1 + x^2)^p  \mu \right\|_{C 
({\mathbb{R}})} \ . $$

\noindent
That is why for arbitrary $ p \geq 1  \ $ the sequence $
\left\{ \lambda_p^{(n)} \right\}_{ n \geq n_p  } \ $ has a 
finite limit $\lambda_p \in \Lambda_f \ $ and since by 
(3.4.5) function $f \ $ has only simple zeros then all 
numbers  $\left\{\lambda_p \right\}_{p \geq 1 } \ $ are 
distinct and hence, $f \ $ is a transcendental entire 
function.}}

\bigskip
\bigskip
{\bf 3.5. Proof of Theorem 3.2.\ }
It is easy to verify that existence of the indicated by  
theorem polynomials is invariant with respect to the  
change of  variables of the kind $f(a \pm x ) \ $, 
$a \in {\mathbb{R}} \ $, and inequality (3.3.11c) can be 
obtained from (3.3.11b) by  division it on 
$\left(1 - \frac{x}{\lambda_k }\right)^{ m_k} \ $. 
Thus, to prove the theorem it is sufficient to show 
the validity of (3.3.11a), (3.3.11b) under conditions $0 
\notin \Lambda_f  \ $, $f(0) = 1 \ $ and  if ${\mathrm{co}} 
\Lambda_f \neq {\mathbb{R}} \ $ then without loss of 
generality $\Lambda_f \subset (0, + \infty ) \ $ .
That is why in the case ${\mathrm{co}} \Lambda_f \neq {
\mathbb{R}} \ $ (3.3.9) implies $\sum\limits_{\lambda \in 
\Lambda_f }^{} \frac{1}{\lambda } < \infty \ $, $f \in {
\mathcal{E}}_0 ({\mathbb{R}}) \ $ and for $N > \min_{} 
\Lambda_f \ $ polynomials $P_N (x):= \prod\limits_{\lambda 
\in \Lambda_f \cap (0 , N )  }^{} \left(1 - \frac{x}{
\lambda } \right) \ $ satisfy all conditions of the theorem.

Let now ${\mathrm{co}} \Lambda_f = {\mathbb{R}} \ $. Then 
(3.3.9) for $z = i \ $ means  $\sum\limits_{\lambda \in 
\Lambda_f }^{} \frac{1}{\lambda^2 } < \infty \ $ and then 
it is easy to derive from the (3.3.9) at some $z \in {
\mathbb{R}} \setminus  \left\{0\right\} \ $, $f(\pm z) \neq 0 
\ $, an existence of finite limit  $a := $ $ \lim\limits_{
p \to \infty } \sum\limits_{\lambda \in \Lambda_f \cap (
-r_p , R_p ) }^{} \frac{1}{\lambda } \ $
and the following representation:
$$ f(z) = e^{- a z } \prod\limits_{\lambda \in \Lambda_f}^{}
\left(1 - \frac{z}{\lambda } \right) e^{\frac{z}{\lambda } }
 \ . \eqno{(3.5.1)}$$
 
\noindent
So, $f \in {\mathcal{LP}}_{II}^0 \ $. Rename zeros of 
$f \ $ by: $\left\{\lambda_k \right\}_{k \geq 0 } \cup 
\left\{- \lambda_{- l } \right\}_{l \geq 1 } := \Lambda_f 
\ $, $0 < \lambda_k \leq \lambda_{k + 1 } \  $,
$0 < \lambda_{- k - 1 } \leq \lambda_{-k - 2 } \ $ $\forall 
\ k \in {\mathbb{Z}}_0 \ $, and for arbitrary $n \geq 0 \ $,
$m \geq 1 \ $ set:
$$ S(n , m ) := \lim\limits_{p \to \infty } 
\left(\sum\limits_{k = n + 1}^{n_+ (p) } \frac{1}{\lambda_k 
} \ - \ \sum\limits_{l = m + 1}^{n_- (p) } \frac{1}{
\lambda_{- l } } \right) \ = \ 
a - \sum\limits_{k = 0}^{n } \frac{1}{\lambda_k  } \ + \ 
\sum\limits_{l =  1}^{m } \frac{1}{\lambda_{- l } } \ , 
\eqno{(3.5.2)} $$

\noindent
where $n_+ (p) := \max \left\{ k \in {\mathbb{Z}}_0 \ | \ 
\lambda_k  < R_p \right\} \ $,  
 $n_- (p) := \max \left\{ l \geq 1  \ | \ \lambda_{- l}  
 < r_p \right\} \ $, $p \geq 1 \ $.
Then for
$$ {\mathcal{R}}_{n , m } (z) :=  \lim\limits_{p \to 
\infty }
\left[\prod\limits_{k = n + 1}^{n_+ (p) } \left(1 - 
\frac{z}{\lambda_k } \right) \prod\limits_{l = m + 1}^{n_- 
(p) }  \left(1 + \frac{z}{\lambda_{- l } } \right)\right]
\ , $$

\noindent
using inequality $\log (1 + x ) \leq x \  $
$\forall \ x > -1 \ $, we get the following estimate:
$$ 0 < {\mathcal{R}}_{n , m } (x)  \leq e^{- x S(n,m) }
\ \ \forall \ x \in [- \lambda_{- m } , \lambda_n  ] , \ 
\ n \geq 0 , \ \ m \geq 1 \ . \eqno{(3.5.3)}$$

For fixed  arbitrary $N \geq 1 \ $ let us find now such 
positive integers $p_N , q_N \ $ that polynomial 

$$ P_N (z) := \prod\limits_{k = 0 }^{p_N } \left(1 - 
\frac{z}{\lambda_k } \right)  \prod\limits_{l = 1 }^{q_N } 
\left(1 + \frac{z}{\lambda_{- l } } \right) \ $$ 

\noindent
will satisfy conditions of the theorem.

Choosing subsequences of $\left\{R_N \right\}_{N \geq 1 } \ $
and $\left\{r_N \right\}_{N \geq 1 } \ $ and reindexed  them
we can consider that
$R_N > \lambda_{n_f^+ (N) + 1} \  $, $r_N > \lambda_{n_f^- 
(N) + 1 } \  $ $\forall \ N \geq 1 \ $, where
$n_f^+ (N) := \max \left\{ k \in {\mathbb{Z}}_0 \ | \right. 
\ $ $ \left.  \lambda_k < N \  \right\} \ $, 
$n_f^- (N) := \max \left\{ l \geq 1  \ | \ \lambda_{- l} 
< N \  \right\} \ $. Therefore inequalities $p_N \geq n_+ 
(N) \ $, $q_N \geq n_- (N) \ $ yield  validity of 
(3.3.11a).

Denote $S(q) := S (n_+ (q) , n_- (q) ) \ $, $q \geq 1 \ $,
and observe that $\lim\limits_{q \to \infty } S (q) = 0 \ $.

Let $S(N) = 0 \ $. Setting $p_N := n_+ (N) \ $, $q_N := n_- 
(N) \ $ we get from (3.5.3) validity of (3.3.11b).

Let $S(N) > 0 \ $. Since  function 
$\varphi_+ (n) := S (n_+ (N) +n , n_- (N) ) \  $, $n \in {
\mathbb{Z}}_0 \ $, decreases from $S(N) \ $ to
$- \infty \leq - \sum\limits_{l = 1 + n_- (N) }^{\infty }
\frac{1}{\lambda_{- l } } < 0 \ $, then it is possible to 
find such $r \in {\mathbb{Z}}_0 \ $ that $\varphi_{+ } 
(1 + r )\leq 0 < \varphi_{+ } ( r) \ $, where $\varphi_{+ } 
(1 + r ) = \varphi_{+ } ( r) - 1/ \lambda_{1 + r + n_+ 
(N)}  \ $. It follows from (3.5.3) that
 
 $${\mathcal{R}}_{n_+ (N) + r + 1 , n_- (N) } (x) \leq
 \left\{\begin{array}{ll}
1 ,   & \forall \ x \in [- \lambda_{- n_- (N)} , 0  ] \ ; \\
{\mathrm{exp}}\left(-x  \varphi_{+ } ( r) + \frac{x}{
\lambda_{1 + r + n_+ (N)} } \right) \ \leq \ e \ ,     
&  \forall \ x \in [0 , \lambda_{1 + r + n_+ (N) } ] \ , 
   \end{array}\right.$$
   
\noindent
i.e. (3.3.11b) will be true for $p_N := 1 + r + n_+ (N) \ $
and $q_N := n_- (N) \ $.

Let $S(N) < 0 \ $. Since the function $\varphi_{-} (m) := 
S(n_+ (N) , n_- (N) + m ) \  $, $m \geq 0 \ $, increases from
$S(N) < 0 \ $ to $0 < \sum\limits_{k \geq 1 + n_+ (N)}^{} 
\frac{1}{\lambda_k } \leq + \infty \ $, then one can find
such $r \in {\mathbb{Z}}_0 \ $ that $\varphi_- (r) < 0 \leq 
\varphi_- (r + 1 ) \ $ and, obviuosly, $\varphi_- (r + 1 ) = 
\varphi_- (r) + 1 / \lambda_{- (r + 1 + n_- (N) ) } 
\ $. As well as in the previous case using inequality (3.5.3)
we obtain validity of (3.3.11b) for $p_N := n_+ (N) \ $,
$q_N := 1 + r + n_- (N) \ $.

Observe now that according to our choice $\lim\limits_{N \to 
\infty } S(p_N , q_N ) = 0 \ $ and therefore by (3.5.2) we 
get:

$$ \sum\limits_{k = 0}^{p_N } \frac{1}{\lambda_k  } \ - \ 
\sum\limits_{l =  1}^{q_N } \frac{1}{\lambda_{- l } } \ = \ 
a - S(p_N , q_N ) \to a \ , \ \ N \to \infty \ , $$

\noindent
i.e. due to (3.5.1) for arbitrary $z \in {\mathbb{C}} \ $:

 $$ P_N (z) := e^{- (a - S(p_N , q_N ) ) z }
\cdot \prod\limits_{k = 0 }^{p_N} \left(1 - \frac{z}{
\lambda_{k} }\right) e^{ \frac{z}{\lambda_{k} }}
\cdot \prod\limits_{l = 1 }^{q_N} \left(1 + \frac{z}{
\lambda_{-l} } \right) e^{ \frac{ z}{\lambda_{-l} }}
\ \to \ f(z) \ , \ \ N \to \infty \ , $$

\noindent
and moreover, that convergence is uniform on any
compact subset of ${\mathbb{C}} \ $. Theorem 3.2 is proved.

\newpage
\begin{center}
{\bf{ \large{CHAPTER IV. }}}  {{ 
\large{  Criterion of the polynomial density
in $ L_p (\mu )  \ $ }}}
\end{center}

\bigskip
\bigskip
{\bf 4.1. Representation Theorem. \ } 
Let $w(x)\ $ be a positive and continuous function of real 
$x \ $ such that for each $n = 0, 1, 2, \ldots $, \ 
$x^n w(x) \ $ is bounded on the whole real line.
In 1924 S. Bernstein [10] asked for conditions on $w \ $ 
that
algebraic polynomials ${\mathcal{P}} \ $ are dense in the 
space $C^0_w \ $. 
In 1959 L.de Branges [12] gived the following its solution.

\bigskip
\bigskip
{\textsc{De Branges Theorem.}}([12]) 
{\textit{If $w : {\mathbb{R}} \to (0, + \infty ) \ $, $w 
\in C ({\mathbb{R}}) \ $ and $\left\|x^n w \right\|_{C 
({\mathbb{R}}) } < \infty \ $
$\forall \ n \in {\mathbb{Z}}_0 \ $, then ${\mathrm{Close
}}_{C^0_w }{\mathcal{P}} \neq C^0_w \ $  if and only if 
there exists a real entire function $F \ $ of exponential 
type all whose zeros $\Lambda_F \  $ are real and simple  
and which satisfies:}}

$$  \ \ \int\limits_{{\mathbb{R} } }^{} \frac{ \log^+   | 
F(x) |  }{1 + x^2 } \ d  x \ \ < \ \ + \infty \  ,  \ \ 
 \sum_{\lambda \in \Lambda_F  }^{} \frac{1}{ w(\lambda )
|F^{\prime } (\lambda ) | } \ < \ \infty  .  $$

\bigskip
It should be noted that conditions on $w \ $ in that theorem 
in view of Corollary 1.1 mean that $C^0_w \ $ is a Banach 
space, and  by Theorem 3.1 and Corollary 3.1 conditions on 
$F \ $ signify in fact that $F \in {\mathcal{H}} \ $.

In 1996 M. Sodin and P. Yuditskii [41] found a  simpler 
proof of de Branges theorem and gived its version with 
weakened conditions on $w \ $.

\bigskip
\bigskip
{\textsc{Sodin-Yuditskii Theorem.}}([41]) 
{\textit{Let $w : {\mathbb{R}} \to [0, 1] \ $ is an upper 
semicontinuous function on ${\mathbb{R}} \ $ and $\left\|
x^n w \right\|_{C ({\mathbb{R}}) } < \infty \ $ $\forall \ 
n \in {\mathbb{Z}}_0 \ $. Algebraic polynomials 
${\mathcal{P}} \ $ are not dense in $C^0_w \ $ if and only 
if there exists such $B \in {\mathcal{H}} \cap 
{\mathcal{E}}_0 \ $ with zeros $\Lambda_B \subseteq S_w :=  
$ $\left\{x \in {\mathbb{R}} \ | \ w(x) > 0 \ \right\} \ $ 
that }}

$$   \sum_{\lambda \in \Lambda_B  }^{} \frac{1}{ w(\lambda )
|B^{\prime } (\lambda ) | } \ < \ \infty  . \eqno{(4.1.1)} $$

\bigskip
Observe, that if (4.1.1) is valid for some $B \in 
{\mathcal{H}} \ $ then according to the established  
by Hamburger [16, 2] property of such functions:
$ \sum_{\lambda \in \Lambda_B  }^{} \frac{\lambda^n }{ 
B^{\prime } (\lambda )  } = 0 \ $ $\forall \ n \in 
{\mathbb{Z}}_0 \ $, we obtain that defined by equality 
$d \mu := \sum_{\lambda \in \Lambda_B  }^{} \frac{
\delta_{\lambda } }{B^{\prime } (\lambda ) } \ $ 
measure $\mu \ $, where $\delta_{\lambda } \  $ denotes 
Dirac's measure at the point $\lambda \ $, belongs to 
${\mathcal{M}} ({\mathbb{R}}) \ $ and evaluates (see 
Th. 1.3) on $C^0_w \ $ linear continuous functional $
\int_{{\mathbb{R}}}^{} f(x) d \mu (x) \ $ vanishing at 
all exponential functions. That is why in that case 
${\mathcal{P}} \ $ is not dense in $C^0_w \ $ and 
taking into account indicated in Chapter I coincidence of 
the seminormed spaces
$C^0_h  \ $ and $C^0_{M_h}  \ $ for arbitrary $h : 
{\mathbb{R}} \to [0,1] \ $, we can reformulate 
aforementioned theorems as follows.

\bigskip
\bigskip
{\textsc{Proposition 4.1.}}([12, 41])
{\textit{Let $h : {\mathbb{R}} \to [0,1] \ $, $\left\|x^n h 
\right\|_{C ({\mathbb{R}}) } < \infty \ $ $\forall \ n \in 
{\mathbb{Z}}_0 \ $, $M_h \ $ is an upper Bair function of $h
\  $ and $S_{M_h } := \left\{x \in {\mathbb{R}} \ | \ 
M_h (x) > 0 \ \right\} \ $. Then the following statements 
hold:}}

\medskip
(1){\textit{ if algebraic polynomials ${\mathcal{P}} \ $ are
not dense 
in $C^0_h \ $ then there exists such $B \in {\mathcal{H}} 
\cap {\mathcal{E}}_0 \ $ that
$$ \Lambda_B \subseteq S_{M_h } \ \mbox{and} \ \ 
\sum\limits_{\lambda \in \Lambda_B  }^{} \frac{1}{M_h (
\lambda ) |B^{\prime } (\lambda )| } < \infty \eqno{(4.1.2)}$$ }}

\medskip
(2){\textit{ if there exists satisfying (4.1.2) $B \in
{\mathcal{H}} \ $, then algebraic polynomials ${\mathcal{P}}
\ $ are not dense 
in $C^0_h \ $.}} 

\bigskip
\bigskip
For any positive integer $N \ $ let ${\mathcal{P}}_N^* \ $ 
denote the set of real algebraic polynomials $P \ $ of 
degree $N \ $ with real and simple zeros only and
$P(0) = 1 \ $. Note that the proof of Theorem 3.3 shows 
that the closure of intersection (3.4.5) and 
${\mathcal{P}}_N^* \ $ is a 
compact set in the topology $\tau_{{\mathcal{E}}} \ $.
Using the Theorems 3.2, 3.3 and Corollary 3.2 it is easy to 
get the validity of the following assertion.

\bigskip
\bigskip
{\textsc{Lemma 4.1.}} \ 
{\textit{Let $w : {\mathbb{R}} \to [0, 1] \ $,
$\left\|x^n w \right\|_{C ({\mathbb{R}})} < \infty \ $ 
$\forall n \in {\mathbb{Z}}_0 \ $,  $w \ $ is an upper 
semicontinuous function on ${\mathbb{R}} \ $, 
 $ \sigma := \chi_{S_w } (0 ) \ $ $ \in \left\{0, 1 \right\} \ $
and function $\theta : [0, + \infty ) \to {\mathbb{R}} \ $ 
for some finite constants $C, c, \alpha > 0 \  $ satisfies 
$c e^{- \alpha x } \leq \theta (x) \leq \frac{C}{1 + x } \ $
$\forall x \geq 0 \ $.  Algebraic 
polynomials ${\mathcal{P}} \ $ are dense in $C^0_w \ $ if 
and only if :}}

$$ \lim\limits_{N \to \infty } \ \ \ \min\limits_{P \in 
{\mathcal{P}}_N^* } \ \left( \sum\limits_{\lambda \in 
\Lambda_P }^{}
\frac{ \theta (|\lambda |) }{|\lambda |}
\ + \  \sum\limits_{\lambda \in \Lambda_P }^{} 
\frac{1}{w(\lambda ) |\lambda |^{\sigma } |P^{\prime } 
(\lambda )| } \right) \ = \ + \infty \ \eqno{(4.1.3)}$$

\bigskip
Calling  to mind (see 1.1) that ${\mathcal{B}}({
\mathbb{R}}) \ $ denotes the family of Borel subsets 
of ${\mathbb{R}} \ $
we formulate now (see also Proposition 2.2) the main 
result of this paper.

\bigskip
\bigskip
\bigskip
{\textsc{Theorem 4.1.}} \ 
{\textit{Let $1 \leq p < \infty \ $ and $\mu \ $ be a 
positive Borel measure on ${\mathbb{R}} \ $ with  all 
finite moments: $\int_{{\mathbb{R}}}^{} |x|^n \ d \mu (x) 
< \infty \  $ $\forall \ n \in {\mathbb{Z}}_0 \ $, and 
unbounded support: ${\mathrm{supp}} \mu := 
\left\{ x \in {\mathbb{R}} \ | \ \mu (x - \delta , x + 
\delta ) > 0 \ \ \forall \ \delta > 0 \right\} \ $.  }}

{\textit{Algebraic polynomials ${\mathcal{P}} \ $ are 
dense in the space $L_p ({\mathbb{R}} , d \mu ) \ $ if 
and only if 
the measure $\mu \ $ can be represented in the following 
form:}}

$$ \mu (A) = \int\limits_{A }^{} w(x)^p \ d \nu (x)
\ \ \ \forall \ A \in {\mathcal{B}}({\mathbb{R}}) \ , 
\eqno{(4.1.4)} $$

\noindent
{\textit{where $\nu \ $ is some finite positive Borel 
measure on ${\mathbb{R}} \ $ and $w \ $ is some upper 
semicontinuous on  
${\mathbb{R}} \ $ function $w : {\mathbb{R}} \to [0,1] \ $,
$\left\|x^n w \right\|_{C ({\mathbb{R}})} < \infty \ $ 
$\forall \ n \in {\mathbb{Z}}_0 \ $, for which algebraic 
polynomials ${\mathcal{P}} \ $ are dense in the 
seminormed space $C^0_w \ $, i.e. $w \ $ satisfies (4.1.3). }}

\bigskip
\bigskip
\bigskip
It is interesting to note that by Theorem 1.3 representation
(4.1.4) for the measure $\mu \ $ holds if and only if $L(f)
:= \int_{{\mathbb{R}}}^{} f(x) \ d \mu (x) \ $ is a linear 
continuous functional on the seminormed space $C^0_{w^p } \ $. 
Another equivalent to Theorem 4.1 statements and their 
important corollaries will be given 
in the second part of that paper.

\bigskip
\bigskip
{\bf 4.2.  Preliminary Lemmas. \ }

\bigskip
\bigskip
{\textsc{4.2.1. Formulations.}} \

\bigskip
For arbitrary function $f : {\mathbb{R}} \to [- \infty , + 
\infty ] \ $, denote
$$ {\mathrm{dom}} f := \left\{ x \in {\mathbb{R}} \ | \  
f (x) < + \infty \ \right\} \ , \ \ 
{\mathrm{epi}} f := \left\{ (x , y ) \in {\mathbb{R}}^2 \ | 
\ y \geq f(x) \ \right\} \ . \eqno{(4.2.1)}$$

\bigskip
\bigskip
{\textsc{Lemma 4.2.}} 
{\textit{\ Let $\mu \ $ be a positive Borel measure on $
{\mathbb{R}} \ $ and function $\alpha : {\mathbb{R}} \to 
[0, + \infty ] \ $ is  lower semicontinuous on $
{\mathbb{R}} \ $ with $\mu ({\mathrm{dom}} \alpha ) > 0
\ $. Then there exists such lower semicontinuous on
${\mathbb{R}} \ $ function $\beta \ $ that:}}

\medskip
(4.2.2a) $\beta (x) \geq \alpha (x) \ $ $\forall \ x \in
{\mathbb{R}} \ $; \ \ 
(4.2.2b) $\mu \left(x \in {\mathbb{R}} \ | \ 
\beta (x) \neq  \alpha (x) \ \right) = 0 \ $ \ ;

\medskip
(4.2.2c) $\mu \left( y \in {\mathbb{R}} \ | \ |x - y | +
|\beta(x) - \beta(y) | < \varepsilon \right) > 0 \ $
\ $\forall \ \varepsilon > 0 \ $ \ $\forall \ x \in {
\mathrm{dom}} \beta \ $. 

\bigskip
\bigskip
Denote by ${\mathcal{K_R   }} \ $ the class of entire functions
$f \ $satisfying conditions: 

\bigskip
(3.3.14a) $f \ $ is a real function with only real and simple
zeros $\Lambda_f \  $
 \ ; 

\medskip
(3.3.14b) if ${\mathrm{co}} \Lambda_f = {\mathbb{R}} \  $
then $f \ $ is of 
exponential type, but if ${\mathrm{co}} \Lambda_f \neq  
{\mathbb{R}} \  $ 
then $f \ $ is of  minimal exponential type \ ; 

\medskip
(3.3.17)$ \limsup\limits_{|\lambda | 
\to \infty , \ \lambda \in \Lambda_f  } \ \  \frac{\log^+ 
\frac{1}{|f^{\prime }(\lambda )| }}{ \log |\lambda | } \ 
< + \infty \ . $

\bigskip
\noindent
As it was noted after the Corollary 3.1 every function $f 
\in {\mathcal{K_R   }} \ $ can be restored by its zeros up 
to a constant factor:
$$ f(z) = f^{(q)} (0) \cdot z^q \cdot 
\lim\limits_{ R \to + \infty } \prod\limits_{|\lambda | < R ,
\  \lambda \in \Lambda_f \setminus \left\{0\right\} }^{}
\left(1 - \frac{z}{\lambda }\right) \ , \ \ 
q \in \left\{0, 1 \right\} \ , \ \ z \in {\mathbb{C}} \ . 
\eqno{(4.2.3)} $$

\bigskip
\bigskip
{\textsc{Lemma 4.3.}} \ {\textit{Let $X \in \left\{ 
{\mathcal{K_R}}, {\mathcal{K}},  {\mathcal{H}} \right\} \ $.
For arbitrary $B \in X \ $ with zeros $\left\{ a_n 
\right\}_{n \geq 1} := 
\Lambda_B \ $ there exist such constant $C > 0 \ $ and such 
sequence of real positive numbers $\left\{\delta_n
\right\}_{n \geq 1 } \ $ that for any sequence of real 
numbers $\left\{b_n \right\}_{n \geq 1 } \ $ satisfying
condition:}}
$$ |b_n - a_n | \leq \delta_n \ \ \forall \ n \geq 1 \ 
\eqno{(4.2.4)} $$

\noindent
{\textit{it is possible to find such $D \in X \ $ that 
$\Lambda_D = \left\{b_n \right\}_{n \geq 1 } \ $ and}}
$$ |B^{\prime } (a_n )| \leq C \cdot |D^{\prime } (b_n )|
\ \ \  \forall \ n \geq 1 \ . \eqno{(4.2.5)} $$

\bigskip
\bigskip
{\textsc{4.2.2. Proof of Lemma 4.2.}} \

For arbitrary $ A \subseteq {\mathbb{R}}^2 \ $ denote 
$P[A] := \left\{x \in {\mathbb{R}} \ | \ \exists y \in 
{\mathbb{R}} : (x , y ) \in A \right\} \ $ and let 
$E_{\alpha } \ $ denote the set of those points 
$\overline{x} := (x , y ) \in {\mathrm{epi }} \alpha \ $,
for which one can find such $\varepsilon (\overline{x}) > 
0 \ $ that:
$$ \mu \left( P\left[ \ \left( \ \overline{x} + \varepsilon 
(\overline{x}) V \ \right) \ \cap \  {\mathrm{epi}} \alpha 
\ \right]\right) \ = \ 0 , \eqno{(4.2.6)}$$

\noindent
where $V := \left\{ (x , y) \in {\mathbb{R}}^2 \ | \ |x| + 
|y| < 1 \ \right\} \ $. The known property of separable 
metric spaces [3, I.5, Lemma 2] means that from the open 
covering of the set
$$ G := \bigcup\limits_{\overline{x} \in E_{\alpha } }^{}
\left(\overline{x} + \varepsilon (\overline{x}) V \right)
 \ \eqno{(4.2.7)} $$
 
\noindent
it is possible to extract the countable subcovering:
$G = \bigcup_{n \geq 1 }^{} G_n \ $, $G_n := \overline{x}_n 
+ \varepsilon (\overline{x}_n ) V \ $, $\overline{x}_n \in 
E_{\alpha } \ $ $\forall \ n \geq 1 \ $.
That is why equality (4.2.6) and countable semiadditivity 
of the measure $\mu \ $ imply : $\mu (P [G \cap {
\mathrm{epi}}\alpha  ] ) = $$\mu \left( \bigcup_{n 
\geq 1 }^{} P [G_n \cap {\mathrm{epi}} \alpha  ]\right) 
\leq 0 \ $, i.e.
$$ \mu (P [ \ G \ \cap \ {\mathrm{epi}} \alpha \ ] ) = 0 \ . 
\eqno{(4.2.8)}$$

\noindent
Since every point in $G \cap {\mathrm{epi}} \alpha \ $ 
possesses 
the property (4.2.6) then the set $B := 
({\mathrm{epi}}\alpha )  
\setminus G \  $ contains all those points $\overline{x } 
\in {\mathrm{epi}} \alpha \ $ for which
$$ \mu (P [ \ ( \ \overline{x} + \varepsilon V \ ) \ \cap 
\ {\mathrm{epi}} \alpha \ ] ) > 0 \  \ \ \forall \ \varepsilon 
\ > 0 \ . \eqno{(4.2.9)}$$

\noindent
In addition an evident inclusion $P[(\overline{x} + 
\varepsilon V) \cap {\mathrm{epi}} \alpha ] \subseteq P[(
\overline{x} + (0, t ) + \varepsilon V ) \cap {
\mathrm{epi}} \alpha ] \ $ $\forall \ t \geq 0 \ $, $ 
\varepsilon > 0 \ $, $\overline{x } \in {\mathbb{R}}^2
\  $, yields $B + (0 , t ) \subseteq  B \ $$\forall \ t
\geq 0 \ $. Therefore the set $B \ $ coincides with 
epigraph ${\mathrm{epi}} \beta \ $ of the lower 
semicontinuous on ${\mathbb{R}} \ $ function $\beta \ $
defined by formula:
$$ \beta (x) := \left\{ \begin{array}{ll}
+ \infty \ ,  & x \notin P [B] \ ;   \\
\min \left\{ y \in {\mathbb{R}} \ | \ (x , y) \in 
B \right\} \ ,    & x \in P [B] \ . 
  \end{array}\right. $$
  
\noindent
Validity of (4.2.2a) follows from ${\mathrm{epi}} \beta
\subseteq  {\mathrm{epi}} \alpha \ $. Due to (4.2.8):
$ \mu (P [ {\mathrm{epi}} \alpha \ \setminus \ 
{\mathrm{epi}} \beta ] ) = \mu (P [G \cap {\mathrm{epi}}
\alpha  ] ) = 0 \ $, and so (4.2.2b) is also true.
Since by (4.2.9) for any $\varepsilon > 0 \ $ and 
$\overline{x} \in {\mathrm{epi}} \beta \ $:
$$ 0 <   \mu (P [ (\overline{x} + \varepsilon V) 
\cap {\mathrm{epi}} \alpha  ] ) = 
\mu \left( P[ (\overline{x} + \varepsilon V )  \cap 
{\mathrm{epi}} \beta ] \cup  P[ (\overline{x} + 
\varepsilon V )  \cap ( {\mathrm{epi}} \alpha \setminus 
{\mathrm{epi}} \beta ) ] \right) \ \leq $$
$$ \leq \mu (P [ (\overline{x} + \varepsilon V) \cap 
{\mathrm{epi}} \beta   ] ) \ , $$
 
\noindent
then the property (4.2.2c) is fulfilled. That is why 
$\beta \ $ satisfies all required conditions  and 
Lemma 4.2 is proved.

\bigskip
\bigskip
{\textsc{4.2.3. Proof of Lemma 4.3.}} \
Since the Lemma's statement is invariant with respect 
to the substitution $x \ $ by $x + a$, $a \in {\mathbb{R}}
\ $, we can without loss of generality consider $0 \notin
\Lambda_B \  $. Assuming $0 < |a_1 | \leq |a_2 | \leq
\ldots \leq |a_n |
\leq \ldots \ $, let us set
$$ \rho_k := \min \left\{ 1, |a_k | , \left\{ \ \left| \ 
|a_k | - |a_m | \ \right| \ \ \left| \ m \geq 1 , \ 
|a_m | \neq |a_k | \right. \right\} \ 
\right\} , \ k \geq 1 \ .             \eqno{(4.2.10)} $$

\noindent
An elementary reasonings show that for real constants 
$\alpha , \beta , \rho , \Delta \ $, satisfying
$$ \alpha  \in {\mathbb{R}} \setminus \left\{0\right\} \ ; 
\ \ 0 < \rho < |\alpha  | , \ \ 0 < \Delta \leq \frac{1}{2} 
\rho \ , \ \ |\alpha  - \beta | \leq \Delta \ , 
\eqno{(4.2.11)}$$

\noindent
the following inequality holds:
$$ \left|\left(1 - \frac{x}{\alpha }\right) \cdot \left(1 - 
\frac{x}{\beta } \right)^{-1} \right| \leq \ 1 + 4 
\frac{\Delta }{\rho } \ \ \ \ \forall x \in {\mathbb{R}} 
: \ |x - \alpha | \geq \rho \ . \eqno{(4.2.12)} $$

\noindent
As for every $k \geq 1 \ $ function $B_k (x) := (1 - 
\frac{x}{a_k } )^{-1 } B(x) \ $ is continuous on 
${\mathbb{R}} \ $ and $B^{\prime } (a_k) = - 
a_k B^{\prime } (a_k ) \ $, then there exists such 
$ \alpha_k > 0 \  $ that
$$ |B_k (x)| \geq \frac{1}{2} \ |a_k | \ |B^{\prime } (a_k ) | 
\ \ \ \forall \ x \in {\mathbb{R}} : \ |x - a_k | \leq  
\alpha_k  \ . \eqno{(4.2.13)}$$

\noindent
Let us set
$$ \delta_k := \min \left\{ \alpha_k ,\frac{\rho_k }{ 4 
(1 + a_k^2 ) }  \right\} \ \ \ \forall \ k \geq 1 \ , 
\eqno{(4.2.14)}$$

\noindent
and consider an arbitrary sequence $\left\{b_k \right\}_{
k \geq 1 } \ $ satisfying inequalities (4.2.4). Since \\
$\sum_{k \geq 1 }^{} \frac{1}{|a_k |^{1 + \varepsilon } } 
< \infty \  $ $\forall   \varepsilon > 0 \ $ and $\left|
\frac{1}{b_k} - \frac{1}{a_k} \right| \leq \frac{2}{a_k^2} 
\  $$\ \forall \ k \geq 1 \ $, then due to (4.2.3) one can 
determine an entire function of exponential type by the 
folowing equality:
$$ D (z) := \lim\limits_{R \to + \infty } \prod\limits_{ 
|b_k| < R , \ k \geq 1 }^{} \left(1 - \frac{z}{b_k } \right)
 \ , \eqno{(4.2.15)} $$
 
\noindent
which, obviously, possesses properties (3.3.14a) and (3.3.14b).
It is easy to verify that (4.2.10) and (4.2.14) give
possibility to choose such tending to infinity sequence 
of positive real numbers $R_p \ $, $p \geq 1 \ $, that 
interval $(- R_p , R_p ) \ $ for every $p \geq 1 \ $ will
include the same number  $N_p \ $  of zeros of the 
functions $D (z) \ $
and $B(z) \ $. Therefore the following relation holds:
$$ \frac{B_k (b_k ) }{ (- b_k ) D^{\prime } (b_k ) } = 
\lim\limits_{p \to \infty }
\prod\limits_{m = 1, \ m \neq k }^{N_p } \left(1 - 
\frac{b_k }{b_m }\right)^{-1}  \left(1 - \frac{b_k}{a_m } 
\right) 
\ , \ \ k \geq 1 \ . \eqno{(4.2.16)}$$

\noindent
Applying estimate (4.2.12) we get: 
$\left| \left(1 - \frac{b_k }{b_m }\right)^{-1}  
\left(1 - \frac{b_k}{a_m } \right)   \right| \leq 1 + 
4 \frac{\delta_m }{\rho_m } \leq 1 + \frac{1}{1 + a_m^2 } \ $,
$k , m \geq 1 \ $, $m \neq k \ $, and so by (4.2.13) and
(4.2.14):
$$ C :=  8 \cdot {\mathrm{exp}} \left\{ \sum\limits_{
m \geq 1 }^{} \frac{1}{1 + a_k^2 }\right\}  \geq 4 
\frac{|a_k |}{|b_k|} \frac{|B^{\prime } (b_k )|}{|D^{
\prime }(b_k )|}  \geq  \frac{|B^{\prime } (a_k )|}{|
D^{\prime }(b_k )|}  \ \ \forall \ k \geq 1 \ . 
\eqno{(4.2.17)}$$

\noindent
Since defined in (4.2.17) constant $C \ $ does not depend 
on choice of the sequence $\left\{b_k \right\}_{k \geq 1 }
\ $ then (4.2.17) represents the required inequality (4.2.5),
from where by the Theorem 3.1 and Corollary 3.1 we will have 
$D \in X \ $ for any indicated in Lemma 4.3 choice of the 
class $X \ $. Lemma 4.3 is proved.

\bigskip
\bigskip
{\bf 4.3.  Proof of Theorem 4.1. \ }

 \bigskip
 {\textit{Sufficiency.}} \ Since $\frac{1}{w} \in L_p (\mu)
 \ $ then the density in $L_p (\mu) \ $ of all compactly 
 supported continuous on ${\mathbb{R}} \ $ functions and 
 an evident inequality $\left\|f\right\|_{L_p (\mu) } \leq
 \left\|f\right\|_w \cdot \left\|\frac{1}{w} \right\|_{L_p
 (\mu) } \ $ $\forall \ f \in C^0_w \ $ by virtue of the
 Proposition 2.2 means the density ${\mathcal{P}}\ $ in 
 $L_p (\mu) \ $. 

\bigskip
 {\textit{Necessity.}} \ 
 
 \medskip
  {\textsc{$4.3.1^0. \ $}} \ By (2.2.1c) density of 
  ${\mathcal{P}}\ $ in $L_p (\mu) \ $ is equivalent 
  to the existence of such sequence of polynomials 
  $P_n \in {\mathcal{P}} [{\mathbb{C}}] \ $, $n \geq 1 \ $,
  that
 $$ \alpha_n := \left\| \frac{1}{x + i } - P_n \right\|_{
 L_p (\mu) \ }^p \ \to \ 0 \ , \ \ n \to \infty \ , 
 \eqno{(4.3.1)} $$
 
 \noindent
 where without loss of generality we can assume that 
 $\sum_{n \geq 1}^{} \alpha_n \leq 1 \ $. Then nondecreasing
 sequence of nonnegative continuous on ${\mathbb{R}} \ $ 
 functions
 $$ \varphi_N (x) := \sum\limits_{n = 1 }^{N} \left|
 \frac{1}{x + i } - P_n (x) \right|^p \ , \ \ N \geq 1 \ 
 \eqno{(4.3.2)}$$
 
 \noindent
 satisfies $\left\|\varphi_N \right\|_{L_1 (\mu ) } \leq 1
 \ $ $\forall \ N \geq 1 \ $ and by Beppo-Levi theorem has
 a limit $\varphi \in L_1 (\mu ) \ $: $\left\|\varphi 
 \right\|_{L_1 (\mu )  } \leq 1 \ $. It is easy to see also
 that $\varphi \ $ is a lower semicontinuous function and
 $\mu ({\mathrm{dom}} \varphi ) > 0 \ $.
 
 \medskip
  {\textsc{$4.3.2^0. \ $}}  \ Under conditions of the 
  Theorem 4.1: $0 < s_n := \int_{{\mathbb{R}}}^{} |x|^n \
  d \mu (x)  < \infty \ $ $\forall \  n \in {\mathbb{Z}}_0
  \ $. Therefore the function 
$$ h (x) := 2 s_0 \cdot \sum\limits_{n \geq 0 }^{}
\frac{1}{2^{n + 1} }  \frac{|x|^n }{s_n } \ , \ \ x
\in {\mathbb{R}} \ , $$

\noindent
has the following properties: $h \in C({\mathbb{R}}) \ $, $h(x) \in
[1, + \infty ) \ $ $\forall \ x \in {\mathbb{R}} \ $,
$2 s_0 = \int_{{\mathbb{R}}}^{} h(x) \ d \mu (x) \ $ 
and $\inf_{x \in {\mathbb{R}} } (1 + x^{2n} )^{-1} 
\cdot h(x)  > 0 \ $$\forall \ n \in {\mathbb{Z}}_0 \ $.
Let
$$ a(x) := h(x) + \varphi (x) \ , \ \ x \in {\mathbb{R}} 
\ . \eqno{(4.3.3)}$$

\noindent
Then for $f = a \ $ the following conditions hold:
$$\begin{array}{ll}
(4.3.4a) \ f(x) \geq 1 \ \ \forall \  x \in {\mathbb{R}} 
\ ; & (4.3.4b)\inf\limits_{x \in {\mathbb{R}} } (1 + x^{2n} 
)^{-1} \cdot f(x)  > 0  \ \forall \ n \in {\mathbb{Z}}_0  ;  
\\
(4.3.4c) \ f \ \mbox{{\small{is lower semicontinuous on }}}
{\mathbb{R}} \ ;  & (4.3.4d) \ f \in L_1 (\mu ) \ .  
  \end{array} $$

\medskip
{\textsc{$4.3.3^0. \ $}} \ By (4.3.4a) and (4.3.4d) with 
$f = a \ $, the sequence of positive numbers 
$$t_n := \int\limits_{|x| > n }^{} a(x) \ d \mu (x) \ , 
\ \ 
 \ n \in {\mathbb{Z}}_0 \ , $$

\noindent
tends to zero as $n \to \infty \ $, and therefore one can 
find such subsequence $\left\{ n_k \right\}_{k \in 
{\mathbb{Z}}_0 } \ $, $n_0 := 0 \ $, that 
$\sum_{k \geq 0 }^{} \sqrt{t_{n_k } } < \infty \ $ and 
$t_{n_{k+1 } }  < t_{n_k } \ $ $\forall \ k \in 
{\mathbb{Z}}_0 \ $. Then for the function
$$ \frac{1}{\theta (x) } := \sqrt{t_0 } \left( 
\frac{\chi_{\left\{0 \right\}} (x) }{\sqrt{t_0 } }  + 
\sum\limits_{k \geq 0 }^{} \frac{\chi_{(n_k , n_{k+1} ] }
(|x|) }{\sqrt{t_{n_k } } }  \right) \ $$

\noindent
we have $\theta (x) \to 0 \ $,  $|x| \to + \infty \ $, 
$\theta (x) \ $ is an even lower semicontinuous on 
${\mathbb{R}} \ $ function, which does not increase as 
$x \geq 0 \ $, \ $\theta (x) \in (0 , 1 ] \ $
$\forall \ x \in {\mathbb{R}} \ $ and
$$ \int\limits_{{\mathbb{R}}}^{} \frac{a(x)}{\theta (x)}
\ d \mu(x) = a(0) \cdot \mu (\left\{0 \right\}) + 
\sum\limits_{k \geq 0  }^{} \frac{t_{n_k} - t_{n_{k + 1 
}}}{
\sqrt{t_{n_k } } } \ < \ \infty \ . $$

\noindent
That is why all properties (4.3.4a-d) are valid and for 
$f = \alpha_0 \  $, where $\alpha_0 (x) := \frac{a(x)}{
\theta (x) } \  $, $x \in {\mathbb{R}} \ $.
Applying  Lemma 4.2 to the function $\alpha_0 \  $ we 
obtain the function $\alpha \ $ for which  all conditions 
(4.3.4a-d) with $f = \alpha \ $ will be true and also:

\medskip
(4.3.4e) \ $\alpha (x) \geq \frac{a(x)}{\theta (x)} \ $
$\ \forall \ x \in {\mathbb{R}} \ $ ;

\medskip
(4.3.4g) \ $\mu \left( \ y \in {\mathbb{R}} \ | \ 
|x - y | + |\alpha (x) - \alpha (y)| < \varepsilon \right)
> 0 \ $ $\forall \ \varepsilon > 0 \ $ $\forall \ x \in 
{\mathrm{dom} } \alpha \ $.

\bigskip
{\textsc{$4.3.4^0. \ $}} \ In view of (4.3.1) and (4.3.4a)
with $f = \alpha \ $ we can apply known Riesz's theorem to 
the convergent to zero in the space $L_1 (\mu ) \ $ sequence 
$\frac{1}{\alpha (x) } \left|\frac{1}{x + i } - P_n (x)
\right|^p \ $, $n \geq 1 \ $ (we consider here $1 / + 
\infty := 0 \ $). That is why taking into account $\mu 
({\mathbb{R}} \setminus {\mathrm{dom}} \alpha ) = 0 \ $
we can find such $A \subseteq {\mathbb{R}} \ $,
$\mu (A) = 0 \ $, and such subsequence $\left\{n_k 
\right\}_{k \geq 1 }  \ $ that
$$ \lim\limits_{k \to \infty } \frac{1}{\alpha (x)} 
\left| \frac{1}{x + i } - P_{n_k } (x) \right|^p = 0 \ \ 
\forall \ x \in {\mathbb{R}} \setminus A \ ; \ \ 
{\mathbb{R}} \setminus A \subseteq  {\mathrm{dom}}\alpha \ . 
\eqno{(4.3.5)} $$

\noindent
On the other hand for arbitrary $T > 0 \ $, $x \in 
{\mathrm{dom}}\alpha \ $ and $k \geq 1 \ $ properties 
(4.3.4e), (4.3.3) and (4.3.2) yield:
$$ \frac{1}{\alpha (x)} 
\left| \frac{1}{x + i } - P_{n_k } (x) \right|^p \leq 
\frac{\left| \frac{1}{x + i } - P_{n_k } (x) \right|^p }{ 
a(x)} \theta (x) \leq \theta (T) \ \ \ \forall \ |x| \geq T
\ , $$

\noindent
from where 
$$ \lim\limits_{T \to + \infty } \ \  \sup\limits_{|x| \geq T }
\frac{1}{\alpha (x)} \left| \frac{1}{x + i } - P_{n_k } (x) 
\right|^p = 0 \ . 
\eqno{(4.3.6)}$$

\noindent
By virtue of the Proposition 2.2 established properties 
(4.3.5), (4.3.6) mean that for any countable set $\Lambda 
\subseteq {\mathbb{R}} \setminus A \ $, which has not finite 
limit points the following statement holds:
$$ {\mathcal{P}} \ \ \mbox{{\textit{is dense in }}} \ \ 
C^0_{\beta } \ \ , \ \  \mbox{{\textit{where}}} \ 
\beta (x) := \frac{\chi_\Lambda (x) }{\alpha (x)^{1/ p }}
\ , \ \ x  \in {\mathbb{R}} \ . \eqno{(4.3.7)}$$

\bigskip
\bigskip
{\textsc{$4.3.5^0. \ $ }} \ Let us exhibit that in fact  
more strong than (4.3.7) statement is valid:

$$ {\mathcal{P}} \ \ \mbox{{\textit{is dense in }}} \ \ 
C^0_{w } \ \ , \ \  \mbox{{\textit{where}}} \ 
w (x) := \frac{\chi_{{\mathrm{supp}} \mu } (x) }{\alpha 
(x)^{1/ p }}
\ , \ \ x  \in {\mathbb{R}} \ . \eqno{(4.3.8 )}$$

\noindent
Assuming a contrary we by Proposition 4.1 can find such
$B \in {\mathcal{H}}\cap {\mathcal{E}}_0 \ $, that 
$\Lambda_B \subseteq  {\mathrm{supp}} \mu \cap 
{\mathrm{dom}} \alpha \ $
and 
$$ \sum\limits_{\lambda \in \Lambda_B }^{} \frac{\alpha 
(\lambda )^{1/ p }}{ |B^{\prime } (\lambda )| } \ < \ 
\infty \ . \eqno{(4.3.9)} $$

\noindent
Applying Lemma 4.3 to the function $B \in {\mathcal{H}} \ $,
we can find also corresponding to that function constant 
$C > 0 \ $ and the sequence of positive real numbers 
$\left\{\delta_{\lambda } \right\}_{\lambda \in \Lambda_B }
\ $. Determine now the numbers $b_{\lambda } \ $, 
$\lambda \in \Lambda_B \ $, satisfying $| b_{\lambda } 
- \lambda | \leq \delta_{\lambda } \  $ $\forall \ 
\lambda \in \Lambda_B \  $.

If $\mu (\left\{\lambda \right\}) > 0 \ $, then $\lambda 
\notin A \ $ and let $b_{\lambda } = \lambda \ $ in that 
case.
If $\mu (\left\{\lambda \right\}) = 0 \ $ then choose an 
arbitrary $b_{\lambda } \ $ from the nonempty by (4.3.4g)
set:
$$ \left\{ y \in {\mathbb{R}} \ | \ y \neq \lambda , \ 
|y - \lambda | + |\alpha (y) -  \alpha(\lambda ) | \leq 
\delta^*_{\lambda }  \ \right\} \setminus A \ , 
\eqno{(4.3.10)} $$

\noindent
where $\delta^*_{\lambda } := \min \left\{ \alpha(\lambda ), 
\delta_{\lambda } \right\} \ $. 

Then $\alpha (b_{\lambda } ) \leq 2 \alpha(\lambda ) \ $ 
$\forall \ \lambda \in \Lambda_B  \ $ and constructed by 
such
sequence $\left\{b_{\lambda }\right\}_{\lambda \in 
\Lambda_B } \ $ entire function $D \in {\mathcal{H}} \ $ 
in Lemma 4.3 
will satisfy inequality:
$ |B^{\prime } (\lambda )| \leq C |D^{\prime } (b_{\lambda } )
| \ $ $\forall \ \lambda \in \Lambda_B \ $. In view of 
(4.3.9) this means that
$$ \sum\limits_{\lambda \in \Lambda_B }^{} \frac{\alpha
(b_{\lambda} )^{1/ p }}{ |D^{\prime } (b_{\lambda} )| } 
\ < \ \infty \ . \eqno{(4.3.11)} $$

\noindent
Since the sequence $\left\{b_{\lambda }\right\}_{\lambda 
\in \Lambda_B } \  \subseteq {\mathbb{R}} \setminus A \ $ 
and has not the finite limit points then obtained inequality 
(4.3.11)
contradicts (4.3.7) with $\Lambda = \left\{b_{\lambda }
\right\}_{\lambda \in \Lambda_B } \ $. Thus, statement 
(4.3.8) has been proved.

It remains to observe that defined in (4.3.8) function
$w \ $
in view of (4.3.4a-d) with $f = \alpha \ $ is upper
semicontinuous on ${\mathbb{R}} \ $ and satisfies :
$\left\|x^n w\right\|_{C ({\mathbb{R}})} < \infty \ $
$\forall \ n \in {\mathbb{Z}}_0 \ $, $0 \leq w(x) \leq 
1 \ $
$\forall \ x \in {\mathbb{R}} \ $,  $\frac{1}{w} \in L_p 
(\mu ) \ $. That is why defined by the following equality
$$ \nu (A) := \int\limits_{A}^{} \frac{1}{w(x)^p } \ 
d \mu (x) \ \ \ \ \forall \ A \in {\mathcal{B}} ({\mathbb{R}}) \ ,
$$

\noindent
measure $\nu \ $ will be finite positive Borel measure on 
the real axis. Since the bounded 
function $w(x)^{p} \ $ is Borel we  by  known 
change of variables theorem in the Lebesgue integral 
get the required  representation of the measure $\mu \ $: 
$\mu (A) = \int\limits_{A}^{} w(x)^p  \ d \nu (x) \ $ $ 
\forall \ A \in {\mathcal{B}} ({\mathbb{R}}) \ $. 
Theorem 4.1 is proved. 

\newpage
\begin{center}
{\textsc{References}}
\end{center}
{\small{
 \begin{enumerate}
\item N.I. Akhiezer, {\textit{On the weighted approximation 
of continuous functions on the real axis}}, 
Uspekhi Mat.Nauk  {\bf{11}}(1956), 3-43; AMS Transl.(ser.
2)  {\bf{22}}(1962),
95-137.
\item N.I. Akhiezer, {\textit{The classical moment problem}},
Oliver \& Boyd, Edinburgh,1965, 253 pp.
\item G.P. Akilov, L.V. Kantorovich, {\textit{Functional 
analysis in normed spaces}}, Macmillan, N.Y., 1964, 773 pp.  
\item P.S. Aleksandrov, {\textit{Introduction to the set 
theory and general topology}}, Nauka, Moscow, 1977, 368 pp.
(in Russian).
\item A. Bakan, {\textit{The Moreau-Rockafellar equality for
sublinear functionals}}, Ukr. Math. J., {\bf{41}} (1990),
861-871.
\item Ju.M. Berezanskii, Z.G. Sheftel, G.F. Us, {\textit{
Functional analysis}}, I-II, Birkh\"auser, 1996.
\item C. Berg, {\textit{Indeterminate moment problem and 
the theory of entire functions}}, J.Comput. Appl.Math., 
{\bf{65}}
(1995), 27-55.
\item C. Berg, {\textit{Moment problems and polynomial 
approximation}}, Ann.Fac.Sci.Toulouse, Stieltjes special 
issue(1996), 9-32.
\item C. Berg, H.L. Pedersen, {\textit{Nevanlinna matrices 
of entire functions}}, Math.Nachr., {\bf{171}} \\ (1995), 29-52.
\item S. Bernstein, {\textit{Le probleme de l'approximation 
des fonctions continues sur tout l'axe reel at l'une de ses 
applications}}, Bull.Math. de France.,  {\bf{52}}(1924), 
399-410.
\item A. Borichev, M. Sodin, {\textit{The Hamburger moment 
problem and weighted polynomial approximation on discrete 
subsets of the real line}}, Preprint, Univ. Bordeaux 1 , 
Laboratoire de Math. Pures de Bordeaux, {\bf{71}}(1998),
43 pp. 
\item L. de Branges, {\textit{The Bernstein problem}}, 
Proc. Amer. Math. Soc.,  {\bf{10}} (1959), 825-832.
\item R.E. Edwards, {\textit{Functional analysis}}, Holt, 
Rinehart \& Winston, 1965.
\item I. Hachatryan, {\textit{On weighted approximation of 
entire functions of zero degree by polynomials}}, Zapiski
Kharkovskogo Mat. ob-va, ser.4, vol.XXIX (1963), 129-142
(in Russian).
\item H. Hamburger, {\textit{Uber eine Erweiterung des 
Stieltjesschen Momentenproblems}}, Math. Ann. {\bf{81}}(1920)
235-319; Math.Ann.{\bf{82}}(1921) 120-164; 168-187.
\item H. Hamburger, {\textit{Hermitian transformations of 
deficiency index (1,1), Jacobi matrices and undetermined 
moment problems}}, Amer.J.Math., {\bf{66}}(1944), 489-522.
\item W.K. Hayman, P.B. Kennedy, {\textit{Subharmonic 
functions}}, Acad.Press, 1976.
\item I.I. Hirschman, D.V. Widder, {\textit{The convolution 
transform}}, Princeton Univ. Press, Princeton, NJ, 1955.
\item  T. Holl, {\textit{Sur l'approximation polynomiale des 
fonctions contenues d'une variable}}, Proc. of the 9th Congr.,
Math.Scand,1939.
\item S. Isumi, T. Kawata, {\textit{Quasi-analytic class and 
closure of $\{ t^n \} \ $ in the interval $(- \infty ,  
\infty ) \ $}}, Tohoku Math.J., {\bf{43}}(1937), 267-273. 
\item S. Karlin, {\textit{Total positivity}}, v.I, 
Stanford.Univ.Press, Stanford, Calif., 1968.
\item P. Koosis, {\textit{Introduction to $H_p \ $ spaces}},
Cambr.Univ.Press, 1980.
\item P. Koosis, {\textit{The logarithmic integral}} I,
 Cambridge Univ.Press, Cambridge, 1988.
\item M.G. Krein, {\textit{On the theory of entire functions
of exponential type}}, Izv.Akad.Nauk.SSSR, Ser.Mat.  
{\bf{309}}(1947), 11 (in Russian). 
\item M.G. Krein, {\textit{On the indeterminate case of the 
Sturm-Liouville boundary problem in the interval $(0, + 
\infty) \ $}}, Izv.Akad.Nauk.SSSR., {\bf{16}}(1952), 
293-324 (in Russian).
\item K. Kuratowski, {\textit{Topology}}, I-II, Acad.Press, 
1968.
\item E. Laguerre, {\textit{Sur quelques points de la theorie 
des equations numeriques}}, Acta Math., Vol. 4 (1884), 97-120.
\item D. Leviatan, I. Shevchuk, {\textit{Some positive results 
and counterexamples
in comonotone approximation II}}, J. Approx. Theory (to appear).
\item B.Ja.Levin, {\textit{Disribution of zeros of entire 
functions}}, Transl.Math.Mono., v.5, AMS, Providence, RI, 
1964; rev.ed. 1980, 523 pp.
\item B.Ja.Levin, {\textit{Density of functions, 
quasianalyticity and subharmonic majorants}},  Zap. \\ nauchn. 
seminarov LOMI {\bf{170}}(1989), 102-156;
English transl. in J.Soviet Math. {\bf{63}}(1993),  171-201.
\item A.I. Markushevich {\textit{Theory of functions of a 
complex variable}}, I-III, Chelsea Publ. \\ Comp., N.Y., 1965. 
\item S.N. Mergelyan,  
{\textit{Weighted approximation by polynomials}},
Uspekhi Mat.Nauk {\bf{ 11}}(1956), 107-152(Russian);
English transl.(ser.2) {\bf{ 10}} (1958).
\item I.P. Natanson, {\textit{Theory of functions of a real 
variable}}, F. Ungar Pub. Co., N.Y., 1960.
\item I.V. Ostrovskii, {\textit{On a class of entire 
functions}}, Soviet Math.Dokl.,  {\bf{17}}(1976), 977-981.
\item M. Riesz, {\textit{Sur le probleme des moments et 
le theoreme de Parseval correspondant}}, Acta Litt.Ac.Sci 
Szeged,
{\bf{1}}(1923), 209-225.
\item M. Riesz, {\textit{Sur le probleme des moment. 
Troisieme note}}, Ark.Mat.Astronom.Fys.,  {\bf{17}} (16), 
1923.
\item H.H. Schaefer, {\textit{Topological vector spaces}},
Macmillian., N.Y., 1966.
\item J. Shohat, J.D. Tamarkin, {\textit{The problem of 
moments}}, AMS, Providence, RI, rev.ed., 1950.
\item M.L. Sodin, {\textit{A remark to the definition of 
Nevanlinna matrices}}, Math.Physics,Analysis and Geometry,
 (Kharkov) {\bf{3}}(1996), 412-422.
\item M.L. Sodin, {\textit{Which perturbations of 
quasianalytic weights preserve quasianaliticity ? 
How to use de Branges theorem}}, J. d'Anal.Math.
{\bf{69}}(1996), 293-309.
\item M. Sodin, P. Yuditskii, {\textit{Another approach to 
de Branges' theorem on weighted polynomial approximation}}, 
Israel Math.Conf.Proc.,v.11, Amer.Math.Soc., Providence, RI, 1997, 221-227. 
\item E.C. Titchmarsh, {\textit{The theory of functions}},
Oxford Univ.Press, Oxford, 1939.

\end{enumerate}

}}

\end{document}